\documentclass[dvips,ejs, reqno, preprint]{imsart}

\RequirePackage{amsthm, amsmath, natbib, amsfonts, amssymb}%
\RequirePackage[OT1]{fontenc}%
\usepackage{natbib}%

\numberwithin{equation}{section}%
\theoremstyle{plain}%

\usepackage{graphicx, color}%
\usepackage{algorithm2e}

\definecolor{darkblue}{rgb}{0.0,0.0,0.7}

\RequirePackage
{hyperref}%

\hypersetup{%
  pdfauthor = {St\'ephane Ga\"iffas, Guillaume Lecu\'e},%
  pdftitle = {Optimal rates and adaptation in the single-index model
    using aggregation},%
  pdfcreator = {pdflatex},%
  pdfproducer = {pdflatex}}

\startlocaldefs

\newcommand \cA{{\cal A}}%
\newcommand \cC{{\cal C}}%
\newcommand \cF{{\cal F}}%
\newcommand \cI{{\cal I}}%
\newcommand \cN{{\cal N}}%
\newcommand \cS{{\cal S}}%
\newcommand \cX{{\cal X}}%

\newcommand{\smin}{s_{\min}}%
\newcommand{\smax}{s_{\max}}

\newcommand{\T}{^{\top}}%
\newcommand{\prodsca}[2]{\langle #1,#2 \rangle}%
\newcommand{\norm}[1]{\|#1\|}%
\newcommand{\ind}[1]{\mathbf 1_{#1}}%

\newcommand{\sumim}{\sum_{i=1}^m}

\DeclareMathOperator*{\supp}{Supp}

\DeclareMathOperator{\limInf}{liminf}

\DeclareMathOperator*{\argmin}{argmin}

\newcommand{\1}{{\rm 1}\kern-0.24em{\rm I}}

\newtheorem{theorem}{Theorem}%
\newtheorem{lem}{Lemma}%
\newtheorem{assumption}{Assumption}%
\theoremstyle{remark}%
\newtheorem{rem}{Remark}

\endlocaldefs

\doi{10.1214/07-EJS077}
\pubyear{2007}
\volume{1}
\issue{0}
\firstpage{538}
\lastpage{573}

\begin{document}

\begin{frontmatter}
  \title{Optimal rates and adaptation in the single-index model using
    aggregation}%
  \runtitle{Optimal rates and adaptation in the single-index model}

  \begin{aug}
    \author{\fnms{St\'ephane} \snm{Ga\"iffas}\corref{}
      \ead[label=e1]{gaiffas@ccr.jussieu.fr}}
    \address{Universit\'e Paris 6  \\
      Laboratoire de Statistique Th\'eorique et Appliqu\'ee \\
      175 rue du Chevaleret \\
      75013 Paris \\
      \printead{e1}}%
  \end{aug}

  \begin{aug}
    \author{\fnms{ Guillaume} \snm{Lecu\'e}
      \ead[label=e2]{lecue@ccr.jussieu.fr}}
    \address{Universit\'e Paris 6  \\
      Laboratoire de Probabilit\'es et Mod\`eles Al\'eatoires \\
      175 rue du Chevaleret \\
      75013 Paris \\
      \printead{e2}} \runauthor{S. Ga\"iffas and G. Lecu\'e}
  \end{aug}

  \begin{abstract}
    We want to recover the regression function in the single-index
    model. Using an aggregation algorithm with local polynomial
    estimators, we answer in particular to the second part of
    Question~2 from \cite{stone82} on the optimal convergence
    rate. The procedure constructed here has strong adaptation
    properties: it adapts both to the smoothness of the link function
    and to the unknown index. Moreover, the procedure locally adapts
    to the distribution of the design. We propose new upper bounds for
    the local polynomial estimator (which are results of independent
    interest) that allows a fairly general design. The behavior of
    this algorithm is studied through numerical simulations. In
    particular, we show empirically that it improves strongly over
    empirical risk minimization.
  \end{abstract}



\received{\smonth{7} \syear{2007}}

\end{frontmatter}

\section{Introduction}

The single-index model is standard in statistical literature. It is
widely used in several fields, since it provides a simple trade-off
between purely nonparametric and purely parametric
approaches. Moreover, it is well-known that it allows to deal with the
so-called ``curse of dimensionality'' phenomenon. Within the minimax
theory, this phenomenon is explained by the fact that the minimax rate
linked to this model (which is multivariate, in the sense that the
number of explanatory variables is larger than $1$) is the same as in
the univariate model. Indeed, if $n$ is the sample size, the minimax
rate over an isotropic $s$-H\"older ball is $n^{-2 s / (2s + d)}$ for
mean integrated square error (MISE) in the $d$-dimensional regression
model without the single-index constraint, while in the single-index
model, this rate is conjectured to be $n^{-2s / (2s + 1)}$ by
\cite{stone82}. Hence, even for small values of $d$ (larger than $2$),
the dimension has a strong impact on the quality of estimation when no
prior assumption on the structure of the multivariate regression
function is made. In this sense, the single-index model provides a
simple way to reduce the dimension of the problem.

Let $(X,Y) \in \mathbb{R}^d \times \mathbb{R}$ be a random variable
satisfying
\begin{equation}
  \label{eq:model}
  Y = g(X) + \sigma(X) \varepsilon,
\end{equation}
where $\varepsilon$ is independent of $X$ with law $N(0,1)$ and where
$\sigma(\cdot)$ is such that $\sigma_0 < \sigma(X) \leq \sigma_1$
a.s. for some $\sigma_0 > 0$ and a known $\sigma_1 > 0$. We denote by
$P$ the probability distribution of $(X,Y)$ and by $P_X$ the margin
law in $X$ or \emph{design} law. In the single-index model, the
regression function as a particular structure. Indeed, we assume that
$g$ can be written has
\begin{equation}
  \label{eq:SI}
  g(x) = f(\vartheta^{\top} x)
\end{equation}
for all $x \in \mathbb R^d$, where $f : \mathbb R \rightarrow \mathbb
R$ is the \emph{link} function and where the direction $\vartheta \in
\mathbb R^d$, or \emph{index}. In order to make the
representation~\eqref{eq:SI} unique (identifiability), we assume the
following (see for instance the survey paper by
\cite{geenens_delecroix05}, or Chapter~2 in~\cite{horowitz98}):
\begin{itemize}
\item $f$ is not constant over the support of $\vartheta\T X$;
\item $X$ admits at least one continuously distributed coordinate
  (w.r.t. the Lebesgue measure);
\item the support of $X$ is not contained in any linear subspace of
  $\mathbb R^d$;
\item $\vartheta \in S_+^{d-1}$, where $S_+^{d-1}$ is the half-unit
  sphere defined by
  \begin{equation}
    \label{eq:half_unit_sphere}
    S_+^{d-1} = \big\{ v \in \mathbb R ^d \;|\; \norm{v}_2 = 1 \text{
      and } v_d \geq 0 \big\},
  \end{equation}
  where $\norm{\cdot}_2$ is the Euclidean norm over $\mathbb R^d$.
\end{itemize}
We assume that the available data
\begin{equation}
  \label{eq:whole_sample}
  D_n := [(X_i,Y_i); 1 \leq i \leq n]
\end{equation}
is a sample of $n$ i.i.d. copies of $(X,Y)$
satisfying~\eqref{eq:model} and~\eqref{eq:SI}. In this model, we can
focus on the estimation of the index $\vartheta$ based on $D_n$ when
the link function $f$ is unknown, or we can focus on the estimation of
the regression $g$ when both $f$ and $\vartheta$ are unknown. In this
paper, we consider the latter problem. It is assumed below that $f$
belongs to some family of H\"older balls, that is, we do not suppose
its smoothness to be known.

Statistical literature on this model is wide. Among many other
references, see \cite{horowitz98} for applications in econometrics, an
application in medical science can be found in \cite{hardle_xia06},
see also \cite{delecroix_etal03}, \cite{delecroix_etal06} and the
survey paper by \cite{geenens_delecroix05}. For the estimation of the
index, see for instance \cite{spok01}; for testing the parametric
versus the nonparametric single-index assumption, see
\cite{stute_zhu05}. See also a chapter in~\cite{kohler02} which is
devoted to dimension reduction techniques in the bounded regression
model. While the literature on single-index modelling is vast, several
problems remain open. For instance, the second part of Question~2 from
\cite{stone82} concerning the minimax rate over H\"older balls in
model~\eqref{eq:model},\eqref{eq:SI} is still open. The first part,
concerning additive modelling is handled in \cite{yang00}
and~\cite{yang_barron99}.

This paper provides new minimax results about the single-index model,
which provides an answer, in particual, to the latter
question. Indeed, we prove that in
model~\eqref{eq:model},\eqref{eq:SI}, we can achieve the rate $n^{-2 s
  / (2s + 1)}$ for a link function in a whole family of H\"older balls
with smoothness $s$, see Theorem~\ref{thm:UB}. The optimality of this
rate is proved in Theorem~\ref{thm:LB}. To prove the upper bound, we
use an estimator which adapts both to the index parameter and to the
smoothness of the link function. This result is stated under fairly
general assumptions on the design, which include any
``non-pathological'' law for $P_X$. Moreover, this estimator has a
nice ``design-adaptation'' property, since it does not depend within
its construction on $P_X$.

\section{Construction of the procedure}
\label{sec:estimator}

The procedure developed here for recovering the regression does not
use a plugin estimator by direct estimation of the index. Instead, it
\emph{adapts} to it, by aggregating several univariate estimators
based on projected samples
\begin{equation}
  \label{eq:proj_sample}
  D_m(v) := [ (v^{\top} X_i, Y_i), 1 \leq i \leq m],
\end{equation}
where $m < n$, for several $v$ in a regular lattice of $S_+^{d-1}$.
This ``adaptation to the direction'' uses a split of the sample. We
split the whole sample $D_n$ into a \emph{training sample}
\begin{equation*}
  D_m := [(X_i, Y_i); 1\leq i \leq m]
\end{equation*}
and a \emph{learning sample}
\begin{equation*}
  D_{(m)} := [(X_i, Y_i); m + 1 \leq i \leq n].
\end{equation*}
The choice of the split size can be quite general (see
Section~\ref{sec:results} for details). In the numerical study
(conducted in Section~\ref{sec:numerical} below), we consider simply
$m = 3n / 4$ (the learning sample size is a quarter of the whole
sample), which provides good results, but other splits can be
considered as well.

Using the training sample, we compute a family $\{ \bar g^{(\lambda)}
\,;\, \lambda \in \Lambda \}$ of linear (or \emph{weak}) estimators of
the regression $g$. Each of these estimators depend on a parameter
$\lambda = (v, s)$ which make them work based on the data ``as if''
the true underlying index were $v$ and ``as if'' the smoothness of the
link function were $s$ (in the H\"older sense, see
Section~\ref{sec:results}).

Then, using the learning sample, we compute a weight $w(\bar g) \in
[0,1]$ for each $\bar g \in \{ \bar g^{(\lambda)} \,;\, \lambda \in
\Lambda \}$, satisfying $\sum_{\lambda \in \Lambda} w(\bar
g^{(\lambda)}) = 1$. These weights give a level of significance to
each weak estimator. Finally, the adaptive, or \emph{aggregated}
estimator, is simply the convex combination of the weak estimators:
\begin{equation*}
  \hat g := \sum_{\lambda \in \Lambda} w(\bar g^{(\lambda)})
  \bar g^{(\lambda)}.
\end{equation*}
The family of weak estimators consists of univariate local polynomial
estimators (LPE), with a data-driven bandwidth that fits locally to
the amount of data. In the next section the parameter $\lambda = (v,
s)$ is fixed and known: we contruct a univariate LPE based on the
sample $D_m(v) = [(Z_i, Y_i); 1\leq i \leq m] = [(v\T X_i, Y_i); 1\leq
i \leq m]$.

\subsection{Weak estimators\textup: univariate LPE}

The LPE is standard in statistical literature, see for instance
\cite{fan_gijbels96, fan_gijbels95}, among many others.  We construct
an estimator $\bar f$ of $f$ based on i.i.d. copies $[(Z_i, Y_i) ; 1
\leq i \leq m]$ of a couple $(Z, Y) \in \mathbb R \times \mathbb R$
such that
\begin{equation}
  \label{eq:univariateModel}
  Y = f(Z) + \sigma(Z) \epsilon,
\end{equation}
where $\epsilon$ is standard Gaussian noise independent of $Z$,
$\sigma : \mathbb R \rightarrow [\sigma_0, \sigma_1] \subset
(0,+\infty)$ and $f \in H(s, L)$ where $H(s, L)$ is the set of
$s$-H\"olderian functions such that
\begin{equation*}
  | f^{(\lfloor s \rfloor)}(z_1) - f^{(\lfloor s \rfloor)}(z_2) | \leq L
  |z_1 - z_2|^{s - \lfloor s \rfloor}
\end{equation*}
for any $z_1, z_2 \in \mathbb R$, where $L > 0$ and $\lfloor s
\rfloor$ stands for the largest integer smaller than~$s$. This
H\"older assumption is standard in nonparametric literature.

\begin{figure}[htbp]
 \includegraphics[width=12cm]{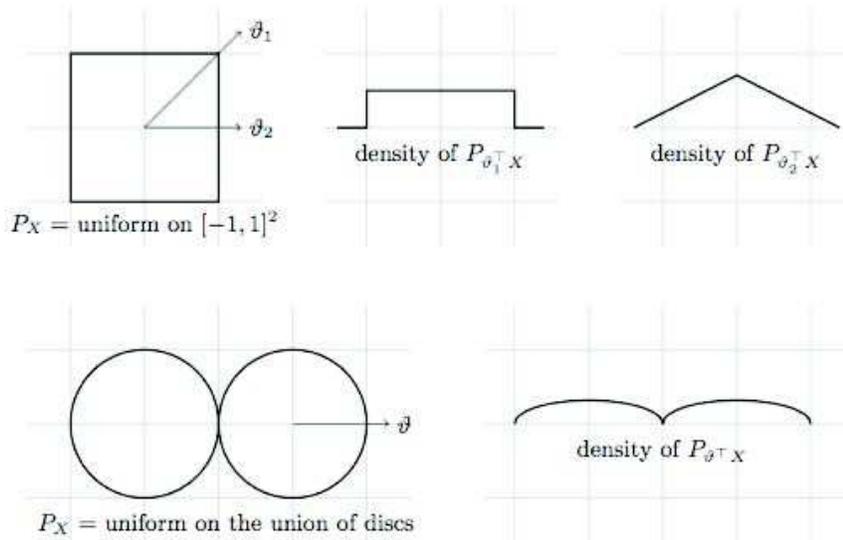}

  \caption{Simple design examples}
  \label{fig:design1}
\end{figure}

Let $r \in \mathbb N$ and $h > 0$ be fixed. If $z$ is fixed, we
consider the polynomial $\bar P_{(z, h)} \in \text{Pol}_r$ (the set of
real polynomials with degree at most $r$) which minimizes in $P$:
\begin{equation}
  \label{eq:defLPE}
  \sumim \big( Y_i - P(Z_i - z) \big)^2 \ind{Z_i \in I(z, h)},
\end{equation}
where $I(z, h) := [z - h, z + h]$ and we define the LPE at $z$ by
\begin{equation*}
  \bar f(z, h) := \bar P_{(z, h)}(0).
\end{equation*}
The polynomial $\bar P_{(z, h)}$ is well-defined and unique when the
symmetrical matrix $\bar {\mathbf Z}_m(z, h)$, with entries
\begin{equation}
  \label{eq:defMatrixZh}
  (\bar {\mathbf Z}_m(z, h))_{a, b} := \frac{1}{m \bar P_Z[ I(z, h) ]}
  \sumim \Big( \frac{Z_i - z}{h} \Big)^{a + b} \ind{Z_i \in I(z, h) }
\end{equation}
for $(a, b) \in \{ 0, \ldots, R\}^2$, is definite positive, where
$\bar P_Z$ is the empirical distribution of $(Z_i)_{1 \leq i \leq m}$,
given by
\begin{equation}
  \label{eq:empirical_design}
  \bar P_{Z} [A] := \frac{1}{m} \sum_{i=1}^m \ind{Z_i \in A}
\end{equation}
for any $A \subset \mathbb R$. When $\bar {\mathbf Z}_m(z, h)$ is
degenerate, we simply take $\bar f(z, h) := 0$. The tuning parameter
$h > 0$, which is called \emph{bandwidth}, localizes the least square
problem around the point $z$ in~\eqref{eq:defLPE}. Of course, the
choice of $h$ is of first importance in this estimation method (as
with any linear method). An important remark is then about the design
law. Indeed, the law of $Z = v\T X$ varies with $v$ strongly: even if
$P_X$ is very simple (for instance uniform over some subset of
$\mathbb R^d$ with positive Lebesgue measure), $P_{v\T X}$ can be
``far'' from the uniform law, namely with a density that can vanish at
the boundaries of its support, or inside the support, see the examples
in Figure~\ref{fig:design1}. This remark motivates the following
choice for the bandwidth.

If $f \in H(s, L)$ for known $s$ and $L$, a ``natural'' bandwidth,
which makes the balance between the bias and the variance of the LPE
is given by
\begin{equation}
  \label{eq:bandwidth}
  H_m(z) := \argmin_{h \in (0, 1)} \Big\{ L h^s \geq \frac{\sigma_1}{(m
    \bar P_{Z}[ I(z, h) ] )^{1/2}} \Big\}.
\end{equation}
This bandwidth choice stabilizes the LPE, since it fits point-by-point
to the local amount of data. We consider then
\begin{equation}
  \label{eq:LPE_bandwidth}
  \bar f(z) := \bar f(z, H_m(z)),
\end{equation}
for any $z \in \mathbb R$, which is in view of Theorem~\ref{thm:LPE}
(see Section~\ref{sec:results}) a minimax estimator over $H(s, L)$ in
model~\eqref{eq:univariateModel}.

\begin{rem}
  The reason why we consider local polynomials instead of some other
  method (like smoothing splines, for instance) is theoretical. It is
  linked with the fact that we need minimax weak estimators under the
  general design Assumption~\ref{assd}, so that the aggregated estimator is
  also minimax.
\end{rem}

\subsection{Adaptation by aggregation}

If $\lambda := (v, s)$ is fixed, we consider the LPE $\bar
f^{(\lambda)}$ given by~\eqref{eq:LPE_bandwidth}, and we take
\begin{equation}
  \label{eq:weak}
  \bar g^{(\lambda)}(x) := \tau_Q( \bar f^{(\lambda)}(\vartheta^{\top} x) ),
\end{equation}
for any $x \in \mathbb R^d$ as an estimator of $g$, where $\tau_Q(f)
:= \max(-Q, \min(Q, f))$ is the truncation operator by $Q > 0$. The
reason why we need to truncate the weak estimators is related to the
theoretical results concerning the aggregation procedure described
below, see Theorem~\ref{thm:oracle} in Section~\ref{sec:results}. In
order to adapt to the index $\vartheta$ and to the smoothness $s$ of
the link function, we aggregate the weak estimators from the family
$\{ \bar g^{(\lambda)} ; \lambda \in \Lambda \}$ with the following
algorithm: we take the convex combination
\begin{equation}
  \label{eq:aggregate}
  \hat g := \sum_{\lambda \in \Lambda} w(\bar g^{(\lambda)}) \bar
  g^{(\lambda)}
\end{equation}
where for a function $\bar g \in \{ \bar g^{(\lambda)} ; \lambda \in
\Lambda \}$, the weight is given by
\begin{equation}
  \label{eq:weights}
  w(\bar g) := \frac{\exp\big( -T R_{(m)}(\bar g) \big)}{\sum_{\lambda
      \in \Lambda} \exp\big(-T R_{(m)}(\bar g^{(\lambda)}) \big) },
\end{equation}
with a \emph{temperature} parameter $T > 0$ and
\begin{equation}
  \label{eq:least_training}
  R_{(m)}(\bar g) := \sum_{i = m + 1}^{n} (Y_i - \bar g(X_i))^2,
\end{equation}
which is the empirical sum of squares of $\bar g$ over the training
sample (up to a division by the sample size). This aggregation
algorithm (with Gibbs weights) can be found in \cite{leung_barron06}
in the regression framework, for projection-type weak
estimators. Cumulative versions of this algorithm can be found in
\cite{catbook:01},~\cite{juditsky_etal05}, \cite{juditsky_nazin05},
\cite{yang:00} and~\cite{yang04}.

We can understand the aggregation algorithm in the following way:
first, we compute the least squares of each weak estimators. This is
the most natural way of assessing the level of significance of some
estimator among the other ones. Then, we put a Gibbs law over the set
of weak estimators. The mass of each estimator relies on its least
squares (over the learning sample). Finally, the aggregate is simply
the mean expected estimator according to this law.

If $T$ is small, the weights~\eqref{eq:weights} are close to the
uniform law over the set of weak estimators, and of course, the
resulting aggregate is inaccurate. If $T$ is large, only one weight
will equal $1$, and the others equal to $0$: in this situation, the
aggregate is equal to the estimator obtained by empirical risk
minimization (ERM). This behavior can be also explained by
equation~\eqref{eq:oracle_minimization} in the proof of
Theorem~\ref{thm:oracle}. Indeed, the exponential
weights~\eqref{eq:weights} realize an optimal tradeoff between the ERM
procedure and the uniform weights procedure. Hence, $T$ is somehow a
regularization parameter of this tradeoff.

The ERM already gives good results, but if $T$ is chosen carefully, we
expect to obtain an estimator which outperforms the ERM. It has been
proved theoretically in \cite{lec9:07} that an aggregation procedure
outperforms the ERM in the regression framework. This fact is
confirmed by the numerical study conducted in
Section~\ref{sec:numerical}, where the choice of $T$ is done using a
simple leave-one-out cross-validation algorithm over the whole sample
for aggregates obtained with several $T$. Namely, we consider the
temperature
\begin{equation}
  \label{CVT}
  \hat T := \argmin_{T \in \mathcal T} \sum_{j=1}^n \sum_{i \neq j} \big(
  Y_i - \hat g_{-i}^{(T)} (X_i) \big)^2,
\end{equation}
where $\hat g_{-i}^{(T)}$ is the aggregated
estimator~(\ref{eq:aggregate}) with temperature $T$, based on the
sample $D_n^{-i} = [(X_j, Y_j) ; j \neq i]$, and where $\mathcal T$ is
some set of temperatures (in Section~\ref{sec:numerical}, we take $
\mathcal T = \{ 0.1, 0.2, \ldots, 4.9, 5 \}$).

The set of parameters $\Lambda$ is given by $ \Lambda := \bar S \times
G$, where $G$ is the grid with step $(\log n)^{-1}$ given by
\begin{equation}
  \label{eq:sgrid}
  G  := \big\{ \smin, \smin + (\log n)^{-1}, \smin + 2 (\log n)^{-1},
  \ldots, \smax \big\}.
\end{equation}
The tuning parameters $\smin$ and $\smax$ correspond to the minimum
and maximum ``allowed'' smoothness for the link function: for this
grid choice, the aggregated estimator converges with the optimal rate
for a link function in $H(s, L)$ for any $s \in [\smin, \smax]$ in
view of Theorem~\ref{thm:UB}.

\begin{figure}[htbp]
  \centering
  \includegraphics[width=6cm]{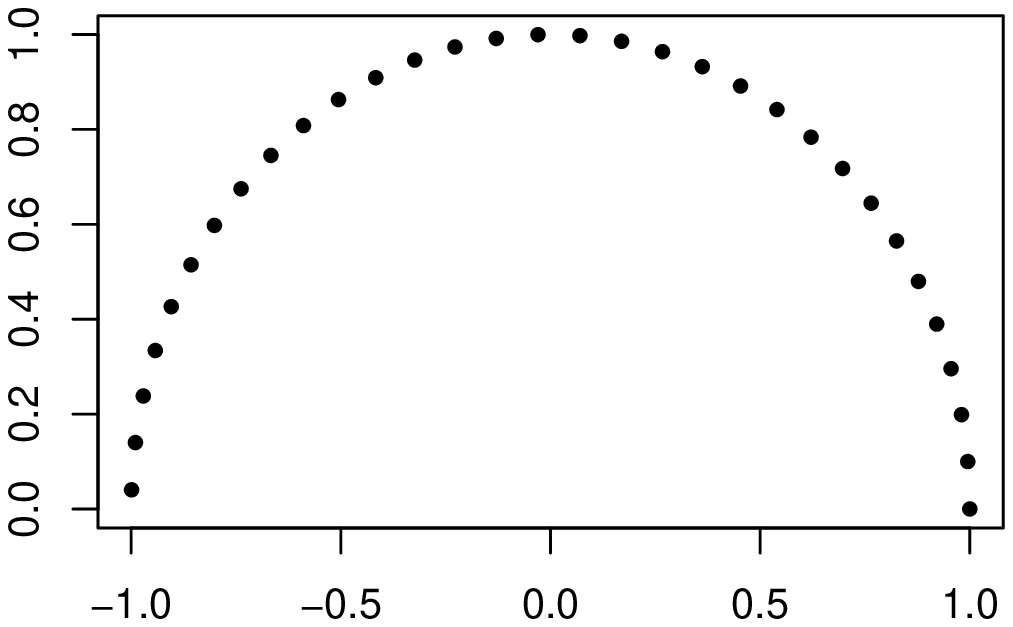}%
  \includegraphics[width=6cm]{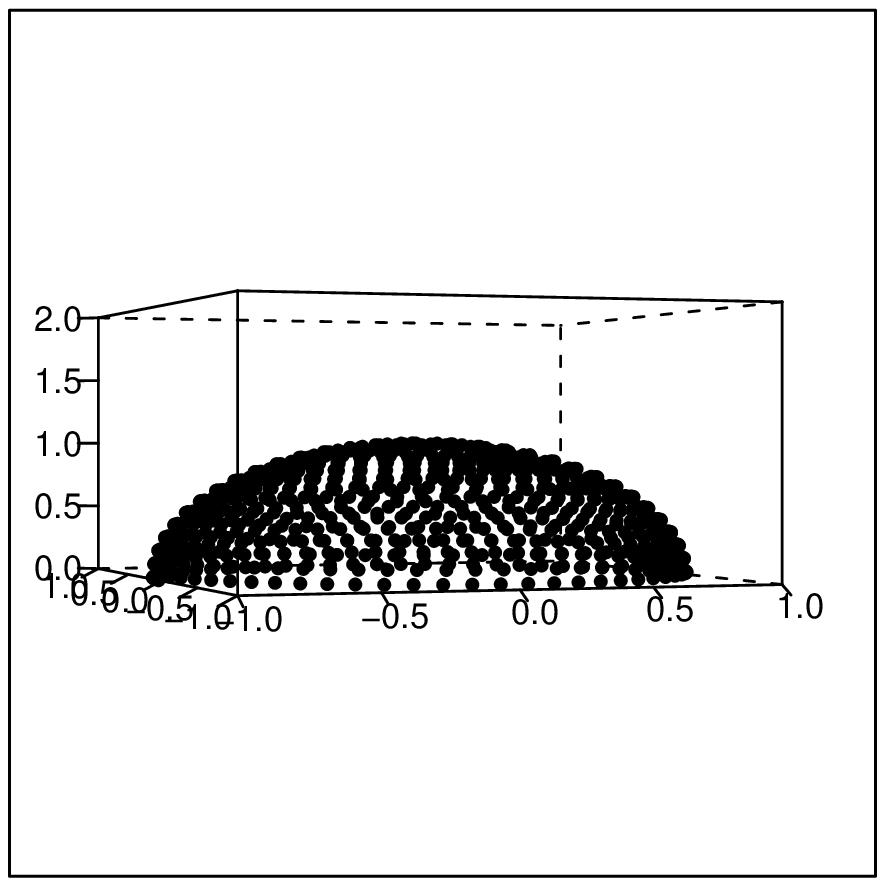}%
  \caption{Lattices $\bar S_\Delta^{d-1}$ for $\Delta = 0.1$ and
    $d=2, 3$}
  \label{fig:lattices}
\end{figure}

The set $\bar S = \bar S_{\Delta}^{d-1}$ is the regular lattice of the
half unit-sphere $S_+^{d-1}$ with discretization step~$\Delta$ which
is constructed as follows. Let us introduce $\Phi(\delta) :=
\cup_{\ell \geq 0} \{ \ell \delta \} \cap [0, \pi]$ and consider the
function $p : [0, \pi]^{d-1} \rightarrow S^{d-1}$ defined by
$p(\phi_1, \ldots, \phi_{d-1}) = (x_1, \ldots, x_d)$, where
\begin{equation*}
  \begin{cases}
    x_1 = \cos(\phi_1) \cos(\phi_2) \times \cdots \times
    \cos(\phi_{d-1}) \\
    x_2 = \sin(\phi_1) \cos(\phi_2) \times \cdots \times
    \cos(\phi_{d-1}) \\
    \phantom{x_2 = \sin(\phi_1) \cos} \vdots \\
    x_\ell = \sin(\phi_{\ell-1}) \cos(\phi_\ell) \times
    \cdots \times \cos(\phi_{d-1}) \\
    \phantom{x_2 = \sin(\phi_1) \cos} \vdots \\
    x_{d-1} = \sin(\phi_{d-2}) \cos(\phi_{d-1}) \\
    x_d = \sin(\phi_{d-1}).
  \end{cases}
\end{equation*}
Then, the regular lattice $\bar S_\Delta^{d-1}$ is constructed using
Algorithm~\ref{alg:lattice}. In Figure~\ref{fig:lattices} we show
$\bar S_\Delta^{d-1}$ for $\Delta=0.1$ and $d=2, 3$. The step is taken
as
\begin{equation}
  \label{eq:lattice_step}
  \Delta = (n \log n)^{-1 / (2 \smin)},
\end{equation}
which relies on the minimal allowed smoothness of the link
function. For instance, if we want the estimator to be adaptive for
link functions at least Lipschitz, we take $\Delta = (n \log
n)^{-1/2}$.

\begin{algorithm}[H]
  \SetLine%
  \KwIn{$d$ (dimension parameter) and $\Delta$ (discretization step)}%
  \KwOut{$\bar S_\Delta^{d-1}$ (regular discretization of $S^{d-1}$)}
  $\bar S_\Delta^{d-1} = \emptyset$ \\
  $\Phi_{d-1} = \Phi(\arccos(1 - \Delta^2 / 2))$ \\
  \ForEach{$\phi_{d-1} \in \Phi_{d-1}$}{
    $\Phi_{d-2} = \Phi(\Delta / \arccos(\phi_{d-1}))$ \\
    \ForEach{$\phi_{d-2} \in \Phi_{d-2}$}{
      $\Phi_{d-3} = \Phi(\Delta / \arccos(\phi_{d-2}))$ \\
      $\vdots$ \\
      \ForEach{$\phi_{2} \in \Phi_{2}$}{
        $\Phi_{1} = \Phi(\Delta / \arccos(\phi_{2}))$ \\
        \ForEach{$\phi_{1} \in \Phi_{1}$}{ add the point of
          coordinates $h(\phi_1, \ldots,
          \phi_{d-1})$ in $\bar S_\Delta^{d-1}$ \\
        } } } }
  \caption{Construction of the regular lattice $\bar
    S_\Delta^{d-1}$.}
  \label{alg:lattice}
\end{algorithm}

\subsection{Reduction of the complexity of the algorithm}
\label{sec:reduccomplex}

The adaptive procedure described previously requires the computation
of the LPE for each parameter $\lambda \in \tilde \Lambda := \Lambda
\times \mathcal L$ (actually, we do also a grid $\mathcal L$ over the
radius parameter $L$ in the simulations). Hence, there are $| \bar
S_{\Delta}^{d-1} | \times | G | \times |\mathcal L| $ LPE to
compute. Namely, this is $ (\pi / \Delta)^{d-1} \times |G| \times
|\mathcal L|$, which equals, if $|G| = |\mathcal L| = 4$ and $\Delta =
(n \log n)^{-1/2}$ (see Section~\ref{sec:numerical}) to $1079$ when $d
= 2$ and to $72722$ when $d = 3$, which is much too large. Hence, the
complexity of this procedure must be reduced: we propose a recursive
algorithm which improves strongly the complexity of the
estimator. Actually, the coefficients $w(\bar g^{(\lambda)})$ are very
close to zero (see Figures~\ref{fig:weightsd2} and~\ref{fig:weightsd3}
in Section~\ref{sec:numerical}) when $\lambda = (v, s)$ is such that
$v$ is ``far'' from the true index $\vartheta$. Hence, these
coefficients should not be computed at all, since the corresponding
weak estimators do not contribute to the aggregated
estimator~\eqref{eq:aggregate}. Thus, instead of using a lattice of
the whole half unit-sphere for detecting the index, we only build a
part of it, which corresponds to the coefficients which are the most
significative. This is done with an iterative algorithm, see
Alogorithm~\ref{alg:preselection}, which makes a preselection of weak
estimators to aggregate ($B^{d}(v, \delta)$ stands for the ball in
$(\mathbb R^d$, $\norm{\cdot}_2)$ centered at $v$ with radius $\delta$
and $R_{(m)}(\bar g)$ is given by~\eqref{eq:least_training}).

\begin{algorithm}[H]
  \SetLine%
  \KwIn{$(X_i, Y_i)$ (Data), $G$ (smoothness grid)}%
  \KwOut{$\hat S$ (a section of $\bar S_\Delta^{d-1}$)}%
  Put $\Delta = (n \log n)^{-1/2}$ and $\Delta_0 = (2 d
  n)^{-1/(2(d-1))}$ \\
  Compute the lattice $\hat S = \bar
  S_{\Delta_0}^{d-1}$ and put $\hat \Lambda := \hat S \times G$ \\
  \While{$\Delta_0 > \Delta$ }{ find the point $\hat v$ such that
    $(\hat v, \hat s) = \hat \lambda = \argmin_{\lambda \in \hat
      \Lambda} R_{(m)}(\bar g^{(\lambda)})$ \\
    put $\Delta_0 = \Delta_0 / 2$ \\
    put $\hat S = \bar S_{\Delta_0}^{d-1} \cap B^{d}(\hat v, 2
    \Delta_0)$ and $\hat \Lambda := \hat S \times G$ \;
  }
  \caption{Preselection of the coefficients}
  \label{alg:preselection}
\end{algorithm}
When the algorithm exits, $\hat S$ is a section of the lattice $\bar
S_{\Delta}^{d-1}$ centered at $\hat v$ with radius $2^{d-1} \Delta$,
which contains (with a high probability) the points $v \in \bar
S_\Delta^{d-1}$ corresponding to the largest coefficients $w(\bar
g^{(\lambda)})$ where $\lambda = (v, s, L) \in \bar S_\Delta^{d-1}
\times G \times \mathcal L$. The aggegate is then computed for a set
of parameters $\hat \Lambda = \hat S \times G \times \mathcal L$
using~\eqref{eq:aggregate} with weights \eqref{eq:weights}. The
parameter $\Delta_0$ is chosen so that the surface of $B^{d}(v,
\Delta_0)$ is $C_d (2 d n)^{-1/2}$: $n$ is not a power of
$d$. Moreover, the number of iterations is $O(\log n)$, thus the
complexity is much smaller than the full aggregation algorithm. This
procedure gives nice empirical results, see
Section~\ref{sec:numerical}. We show the iterative construction of
$\hat S$ in Figure~\ref{fig:iterations}.

\newlength{\figlen}%
\setlength{\figlen}{6cm}
\begin{figure}[t]
  \includegraphics[width=\figlen]{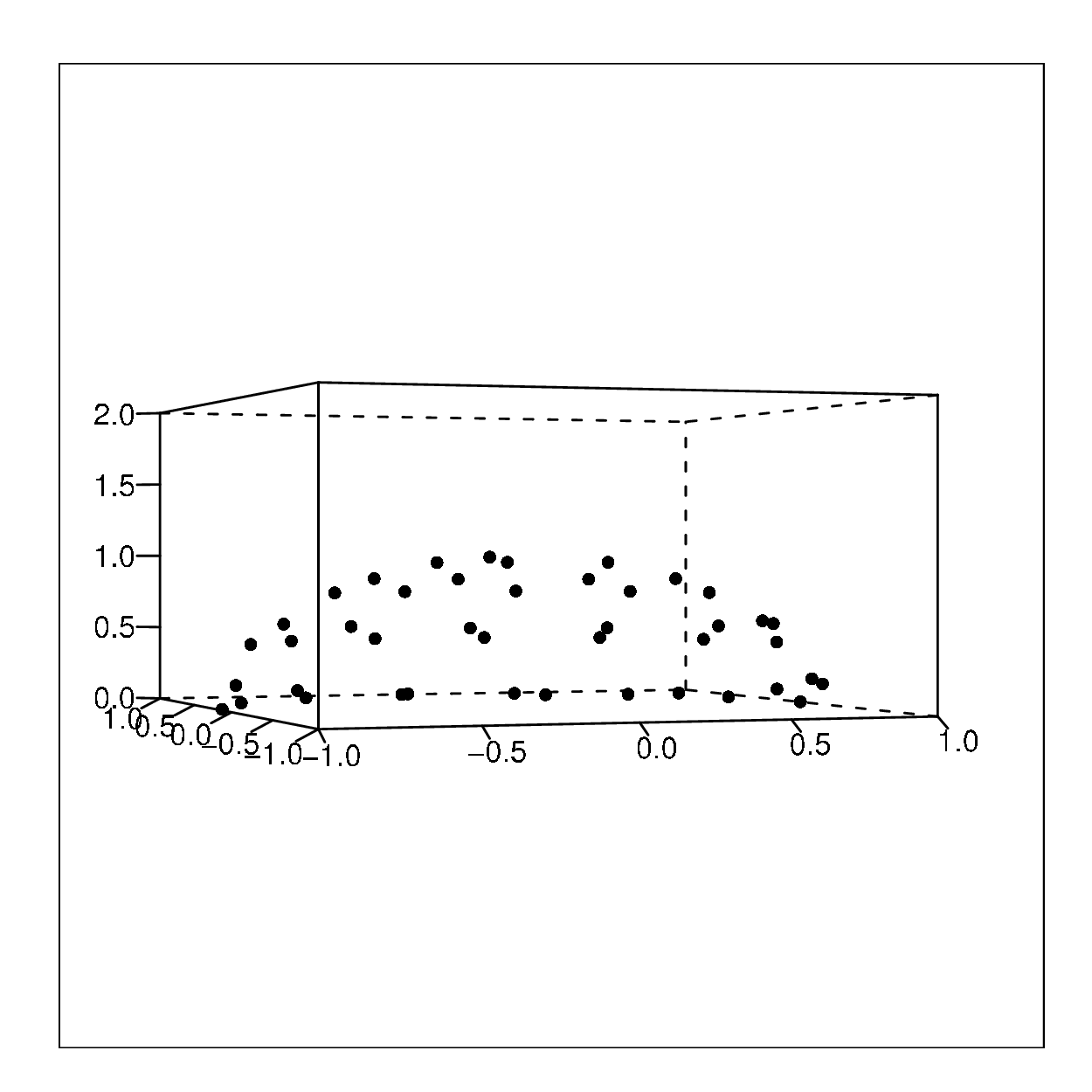}
  \includegraphics[width=\figlen]{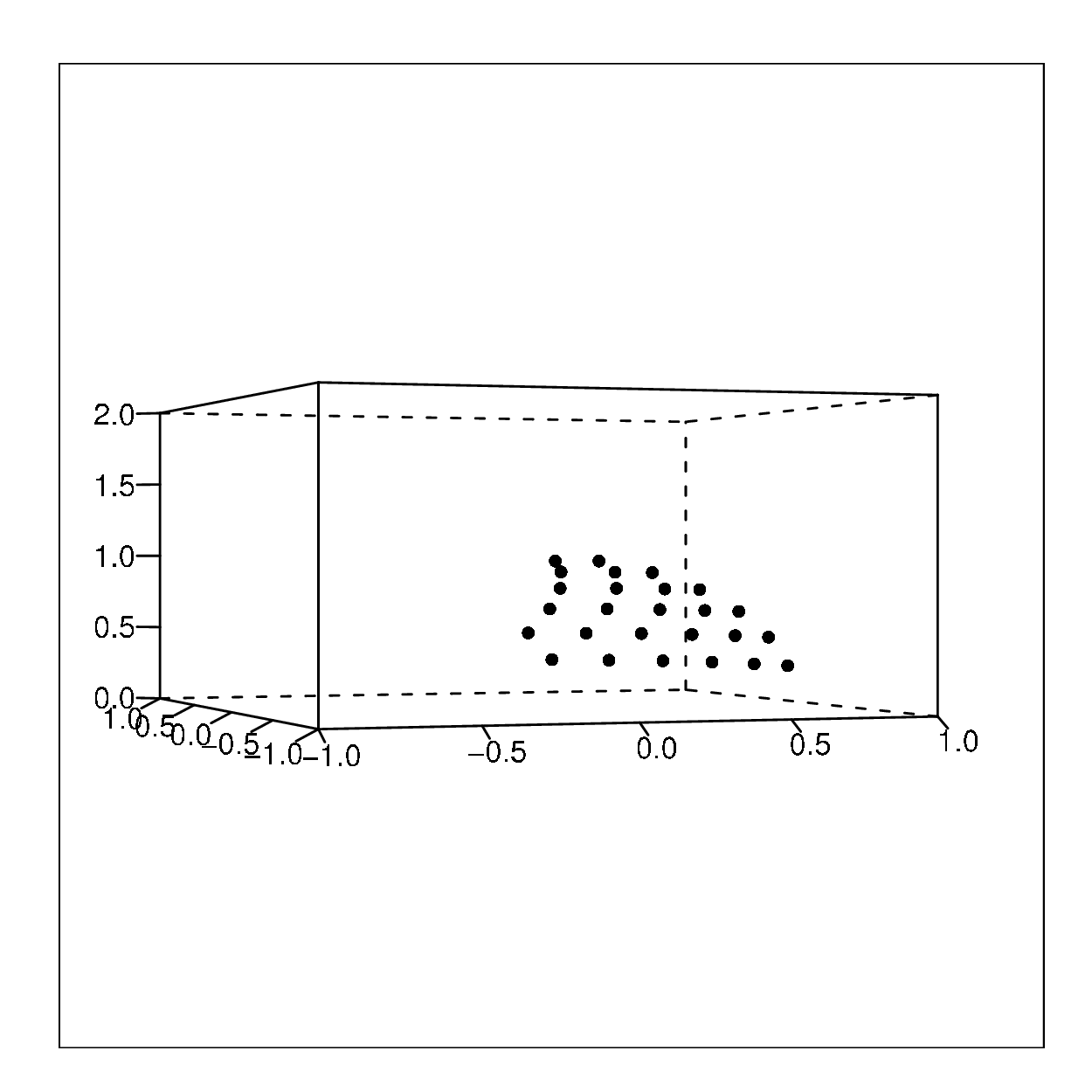} \\
  \includegraphics[width=\figlen]{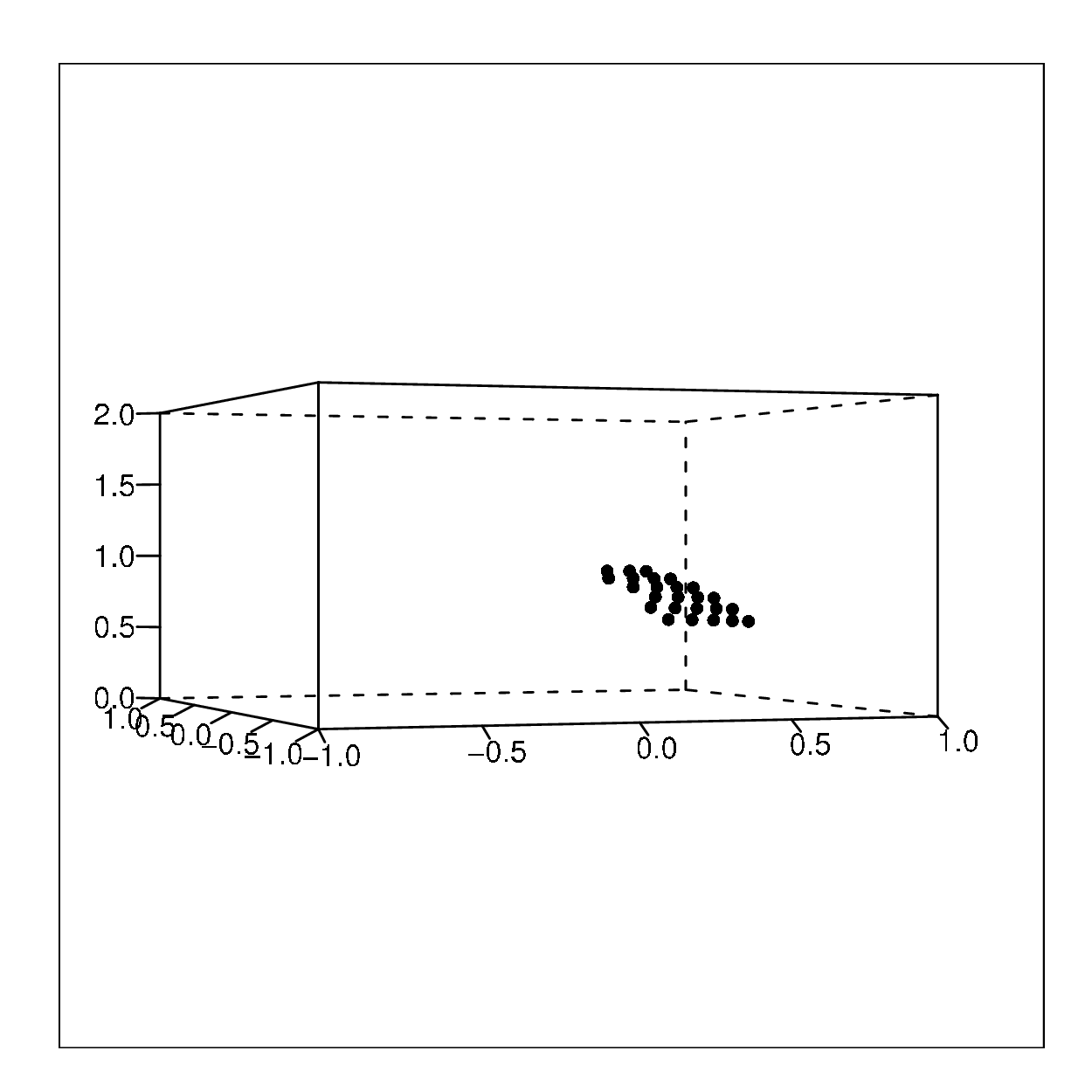}
  \includegraphics[width=\figlen]{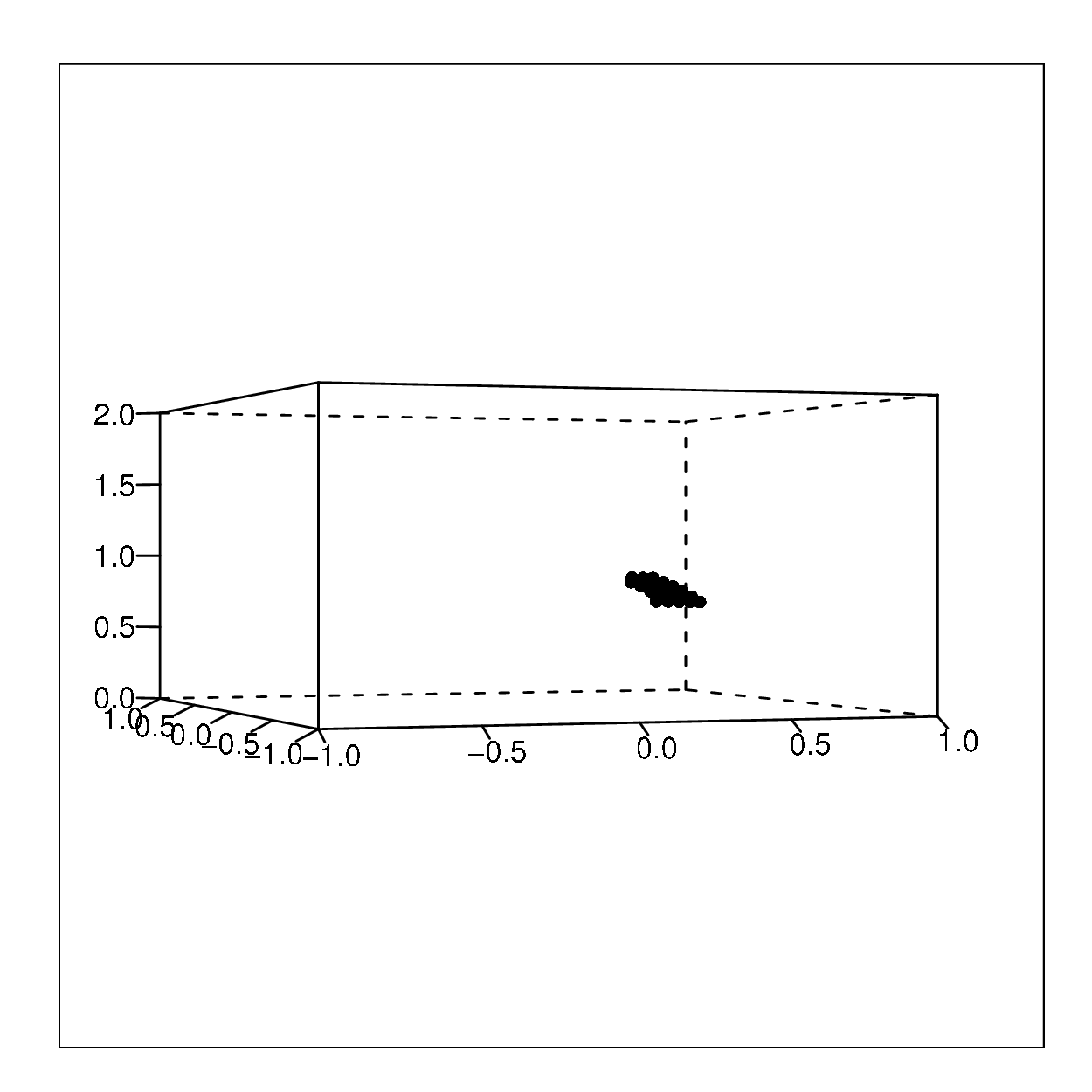}
  \caption{Iterative construction of the set $\hat S$ of
    preselectioned weak estimators indexes. Weak estimators are
    aggregated only for $v \in \hat S$ \textup(bottom right\textup),
    which is concentrated around the true index.}
  \label{fig:iterations}
\end{figure}

\section{Main results}
\label{sec:results}

The error of estimation is measured with the $L^2(P_X)$-norm, defined
by
\begin{equation*}
  \norm{g}_{L^2(P_X)} := \Big( \int_{\mathbb{R}^d} g(x)^2 P_X(dx)
  \Big)^{1/2},
\end{equation*}
where we recall that $P_X$ is the design law. We consider the set
$H^Q(s, L) := H(s, L) \cap \{ f : \mathbb R \rightarrow \mathbb R
\;|\; \norm{f}_{\infty} := \sup_x |f(x)| \leq Q\}$. Since we want the
adaptive procedure to work whatever $\vartheta \in S_+^{d-1}$ is, we
need to work with as general assumptions on the law of $\vartheta\T X$
as possible. The following assumption generalizes the usual
assumptions on random designs (when $P_X$ has a density with respect
to the Lebesgue measure) that can be met in literature.  Namely, we do
not assume that the density of $P_{v\T X}$ is bounded away from zero.
Indeed, even with a very simple $P_X$, this assumption holds for
specific $v$ only (see Figure~\ref{fig:design1}). We say that a real
random variable $Z$ satisfies Assumption~\ref{assd}~if:

\def\theassumption{(D)}
\begin{assumption}
\label{assd}
  There is a density $\mu$ of $P_Z$ with respect to the Lebesgue
  measure which is continuous. Moreover, we assume that
  \begin{itemize}
  \item $\mu$ is compactly supported\textup;
  \item There is a finite number of $z$ in the support of $\mu$ such
    that $\mu(z) = 0$\textup;
  \item For any such $z$, there is an interval $I_z = [z - a_z, z +
    b_z]$ such that $\mu$ is decreasing over $[z - a_z, z]$ and
    increasing over $[z, z + b_z]$\textup;
  \item There is $\beta \geq 0$ and $\gamma > 0$ such that
  \begin{equation*}
    P_Z[I] \geq \gamma |I|^{\beta + 1}
  \end{equation*}
  for any $I$, where $|I|$ stands for the length of $I$.
\end{itemize}
\end{assumption}

This assumption includes any design with continuous density with
respect to the Lebesgue measure that can vanish at several points, but
not faster than some power function.

\subsection{Upper and lower bounds}

The next Theorem provides an upper bound for the adaptive estimator
constructed in Section~\ref{sec:estimator}. This upper bound holds for
quite general tuning parameters. The temperature $T > 0$ can be
arbitrary (but not in practice of course). The training sample size is
given by
\begin{equation}
  \label{eq:learning_size}
  m = [n ( 1 - \ell_n)],
\end{equation}
where $[x]$ is the integral part of $x$, and where $\ell_n$ is a
positive sequence such that for all $n$, $(\log n)^{-\alpha} \leq
\ell_n < 1$ with $\alpha > 0$. Note that in methods involving data
splitting, the optimal choice of the split size is open. The degree
$r$ of the LPE and the grid choice $G$ must be such that $\smax \leq r
+ 1$.

The upper bound below shows that the estimator converges with the
optimal rate for a link function in a whole family of H\"older
classes, and for any index. In what follows, $E^n$ stands for the
expectation with respect to the joint law $P^n$ of the whole sample
$D_n$.

\begin{theorem}
  \label{thm:UB}
  Let $\hat g$ be the aggregated estimator given
  by~\eqref{eq:aggregate} with the weights~\eqref{eq:weights}. If for
  all $v \in S_+^{d-1}$, $v\T X$ satisfies Assumption~\ref{assd}, we have
  \begin{equation*}
    \sup_{\vartheta \in S_+^{d-1}} \sup_{f \in H^Q(s, L)} E^n
    \norm{\hat g - g}_{L^2(P_X)}^2 \leq C n^{-2s / (2s + 1)},
  \end{equation*}
  for any $s \in [\smin, \smax]$ when $n$ is large enough, where we
  recall that $g(\cdot) = f(\vartheta\T \cdot)$. The constant $C > 0$
  depends on $\sigma_1, L, \smin, \smax$ and $P_X$ only.
\end{theorem}

Note that $\hat g$ does not depend within its construction on the
index $\vartheta$, nor the smoothness $s$ of the link function $f$,
nor the design law $P_X$. The assumption that $v\T X$ satisfies
Assumption~\ref{assd} for any $v \in S_+^{d-1}$ holds, for instance, for the
multivariate designs from Figure~\ref{fig:design1}. More generally,
this property holds for any uniform law over a support that does not
have very ``spiky'' boundary. Note that this assumption is more
general than the one considered in~\cite{audibert_tsybakov07}.

In Theorem~\ref{thm:LB} below, we prove in our setting (when
Assumption~\ref{assd} holds on the design) that $n^{-2s / (2s + 1)}$ is a
lower bound for a link function in $H(s, L)$ in the single-index
model.

\begin{theorem}
  \label{thm:LB}
  Let $s, L, Q > 0$ and $\vartheta \in S_+^{d-1}$ be such that
  $\vartheta\T X$ satisfies Assumption~\ref{assd}. We have
  \begin{equation*}
    \inf_{\tilde g} \sup_{f \in H^Q(s, L)} E^n \norm{\tilde g -
      g}_{L^2(P_X)}^2 \geq C' n^{-2s / (2s + 1)},
  \end{equation*}
  where the infimum is taken among all estimators based on data
  from~\eqref{eq:model},\eqref{eq:SI}, and where $C' > 0$ is a
  constant depending on $\sigma_1, s, L$ and $P_{\vartheta\T X}$ only.
\end{theorem}

Theorem~\ref{thm:UB} and Theorem~\ref{thm:LB} together entail that
$n^{-2s / (2s + 1)}$ is the minimax rate for the estimation of $g$ in
model~\eqref{eq:model} under the constraint~\eqref{eq:SI} when the
link function belongs to an $s$-H\"older class. It answers in
particular to Question~2 from \cite{stone82}.

\subsection{A new result for the LPE}

In this section, we give upper bounds for the LPE in the univariate
regression model~\eqref{eq:univariateModel}. Despite the fact that the
literature about LPE is wide, the Theorem below is new. It provides a
minimax optimal upper bound for the $L^2(P_Z)$-integrated risk of the
LPE over H\"older balls under Assumption~\ref{assd}, which is a general
assumption for random designs having a density with respect to the
Lebesgue
measure.

In this section, the smoothness $s$ is supposed known and fixed, and
we assume that the degree $r$ of the local polynomials satisfies $r +
1 \geq s$. First, we give an upper bound for the pointwise risk
conditionally on the design. Then, we derive from it an upper bound
for the $L^2(P_Z)$-integrated risk, using standard tools from
empirical process theory (see Appendix). Here, $E^m$ stands for the
expectation with respect to the joint law $P^m$ of the observations
$[(Z_i, Y_i); 1 \leq i \leq m]$ from
model~\eqref{eq:univariateModel}. Let us define the matrix
\begin{equation*}
  \bar {\mathbf Z}_m(z) := \bar {\mathbf Z}_m(z, H_m(z))
\end{equation*}
where $\bar {\mathbf Z}_m(z, h)$ is given by~\eqref{eq:defMatrixZh}
and $H_m(z)$ is given by~\eqref{eq:bandwidth}. Let us denote by
$\lambda(M)$ the smallest eigenvalue of a matrix $M$ and introduce
$Z_1^m := (Z_1, \ldots, Z_m)$.

\begin{theorem}
  \label{thm:LPE}
  For any $z \in \supp P_Z$, let $\bar f(z)$ be given
  by~\eqref{eq:LPE_bandwidth}. We have on the event $\{ \lambda( \bar
  {\mathbf Z}_m(z)) > 0 \}$\textup:
  \begin{equation}
    \label{eq:LPE_point}
    \sup_{f \in H(s, L)} E^m \big[ (\bar f(z) - f(z))^2 | Z_1^m \big]
    \leq 2 \lambda( \bar {\mathbf Z}_m(z))^{-2} L^2 H_m(z)^{2s}.
  \end{equation}
  Moreover, if $Z$ satisfies Assumption~\ref{assd}, we have
  \begin{equation}
    \label{eq:LPE_L2}
    \sup_{f \in H^Q(s, L)} E^m \big[ \norm{\tau_Q(\bar f) - f}_{L^2(P_Z)}^2
    \big] \leq C_2 m^{-2s / (2s + 1)}
  \end{equation}
  for $m$ large enough, where we recall that $\tau_Q$ is the
  truncation operator by $Q > 0$ and where $C_2 > 0$ is a constant
  depending on $s$, $Q$, and $P_Z$ only.
\end{theorem}

\begin{rem}
  While inequality~\eqref{eq:LPE_point} in Theorem~\ref{thm:LPE} is
  stated over $\{ \lambda( \bar {\mathbf Z}_m(z)) > 0 \}$, which
  entails the existence and the unicity of a solution to the linear
  system~\eqref{eq:defLPE} (this inequality is stated conditionally on
  the design), we only need Assumption~\ref{assd} for
  inequality~\eqref{eq:LPE_L2} to hold.
\end{rem}

\subsection{Oracle inequality}

In this section, we provide an oracle inequality for the aggregation
algorithm~\eqref{eq:aggregate} with weights \eqref{eq:weights}. This
result, which is of independent interest, is stated for a general
finite set $\{ \bar g^{(\lambda)} ; \lambda \in \Lambda \}$ of
deterministic functions such that $\norm{\bar g^{(\lambda)}}_{\infty}
\leq Q$ for all $\lambda \in \Lambda$. These functions are for
instance weak estimators computed with the training sample (or
\emph{frozen} sample), which is independent of the learning
sample. Let $D := [(X_i, Y_i); 1 \leq i \leq |D|]$ (where $|D|$ stands
for the cardinality of $D$) be an i.i.d. sample of $(X, Y)$ from the
multivariate regression model~\eqref{eq:model}, where no particular
structure like~\eqref{eq:SI} is assumed.

The aim of aggregation schemes is to mimic (up to an additive
residual) the oracle in $\{ \bar g^{(\lambda)} ; \lambda \in \Lambda
\}$.  This aggregation framework has been considered, among others, by
\cite{b:04}, \cite{catbook:01}, \cite{jn:00}, \cite{leung_barron06},
\cite{n:00}, \cite{tsy:03} and \cite{yang:00}.

\begin{theorem}
  \label{thm:oracle}
   The aggregation procedure $\hat g$ based on the
  learning sample $D$ defined by~\eqref{eq:aggregate}
  and~\eqref{eq:weights} satisfies
  \begin{equation*}
    E^D \norm{\hat g - g}_{L^2(P_X)}^2 \leq (1 + a) \min_{\lambda \in
      \Lambda} \norm{\bar g^{(\lambda)} - g}_{L^2(P_X)}^2 + \frac{C
      \log |\Lambda| (\log |D|)^{1/2}}{|D|}
  \end{equation*}
  for any $a > 0$, where $|\Lambda|$ denotes the cardinality of
  $\Lambda$, where $E^D$ stands for the expectation with respect to
  the joint law of $D$, and where $C := 3[8 Q^2 (1+a)^2/a + 4 (6 Q^2+2
  \sigma_1 2\sqrt{2})(1 + a)/3]+ 2 + 1 / T$.
\end{theorem}

This theorem is a model-selection type oracle inequality for the
aggregation procedure given by~\eqref{eq:aggregate}
and~\eqref{eq:weights}. Sharper oracle inequalities for more general
models can be found in~\cite{juditsky_etal05}, where the algorithm
used therein requires an extra cumulative sum.

\begin{rem}
  Inspection of the proof of Theorem~\ref{thm:oracle} shows that the
  ERM (which is the estimator minimizing the empirical risk
  $R_{(m)}(g) := \sum_{i = m + 1}^{n} (Y_i - g(X_i))^2$ over all $g$
  in $\{ \bar g^{(\lambda)} ; \lambda \in \Lambda \}$) satisfies the
  same oracle inequality. Nevertheless, it has been proved in
  \cite{lec9:07} that the ERM is theoretically suboptimal in this
  framework, when we want to mimic the oracle without the extra factor
  $1 + a$ in front of the biais term $\min_{\lambda \in \Lambda}
  \norm{\bar g^{(\lambda)} - g}_{L^2(P_X)}^2$. The simulation study of
  Section~\ref{sec:numerical} (especially
  Figures~\ref{fig:misehardsined2}, \ref{fig:misehardsined3},
  \ref{fig:misehardsined4}) confirms this suboptimality.
\end{rem}

\section{Numerical illustrations}
\label{sec:numerical}

\begin{figure}
  \includegraphics[width = 6cm]{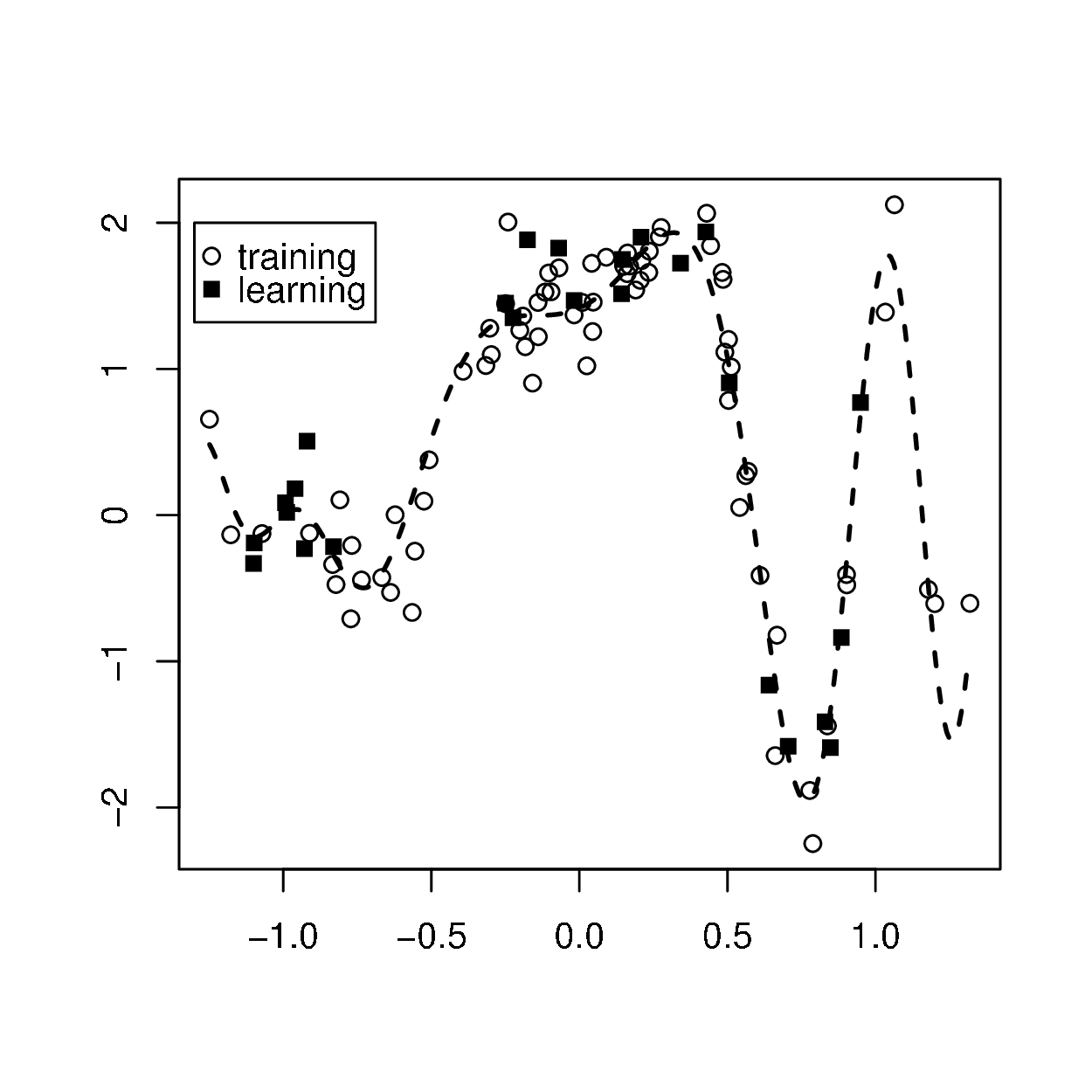}
  \includegraphics[width = 6cm]{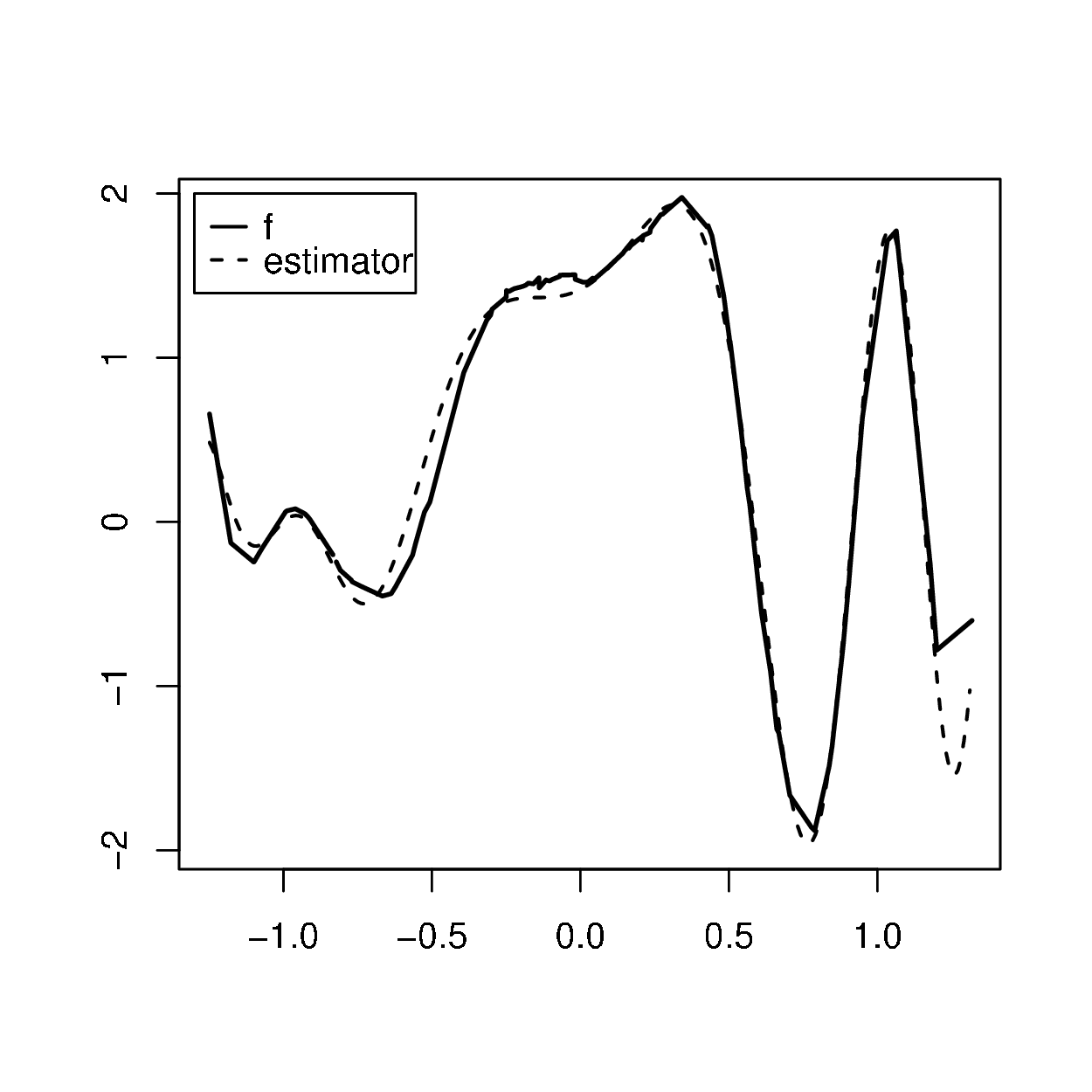}
  \includegraphics[width = 6cm]{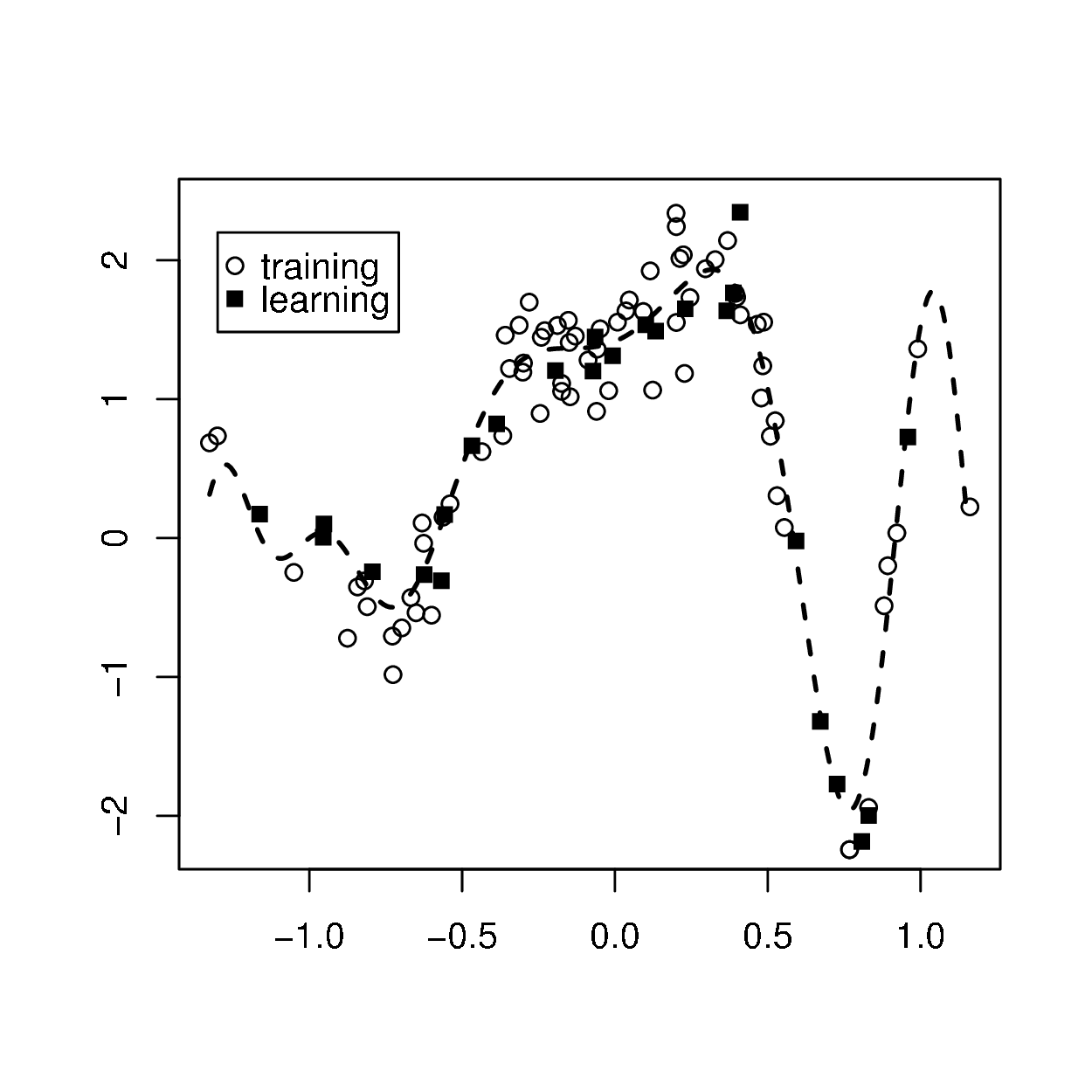}
  \includegraphics[width = 6cm]{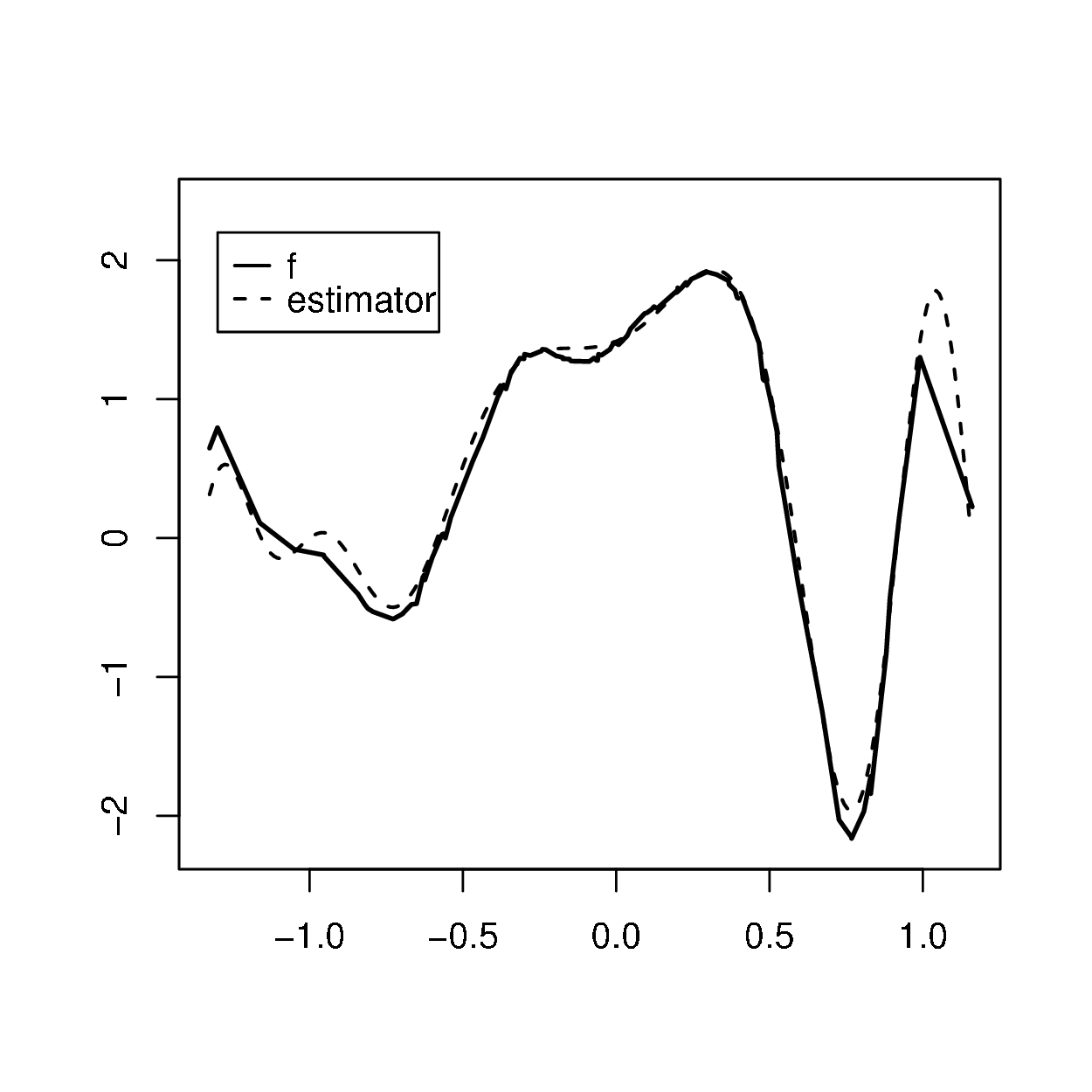}
  \includegraphics[width = 6cm]{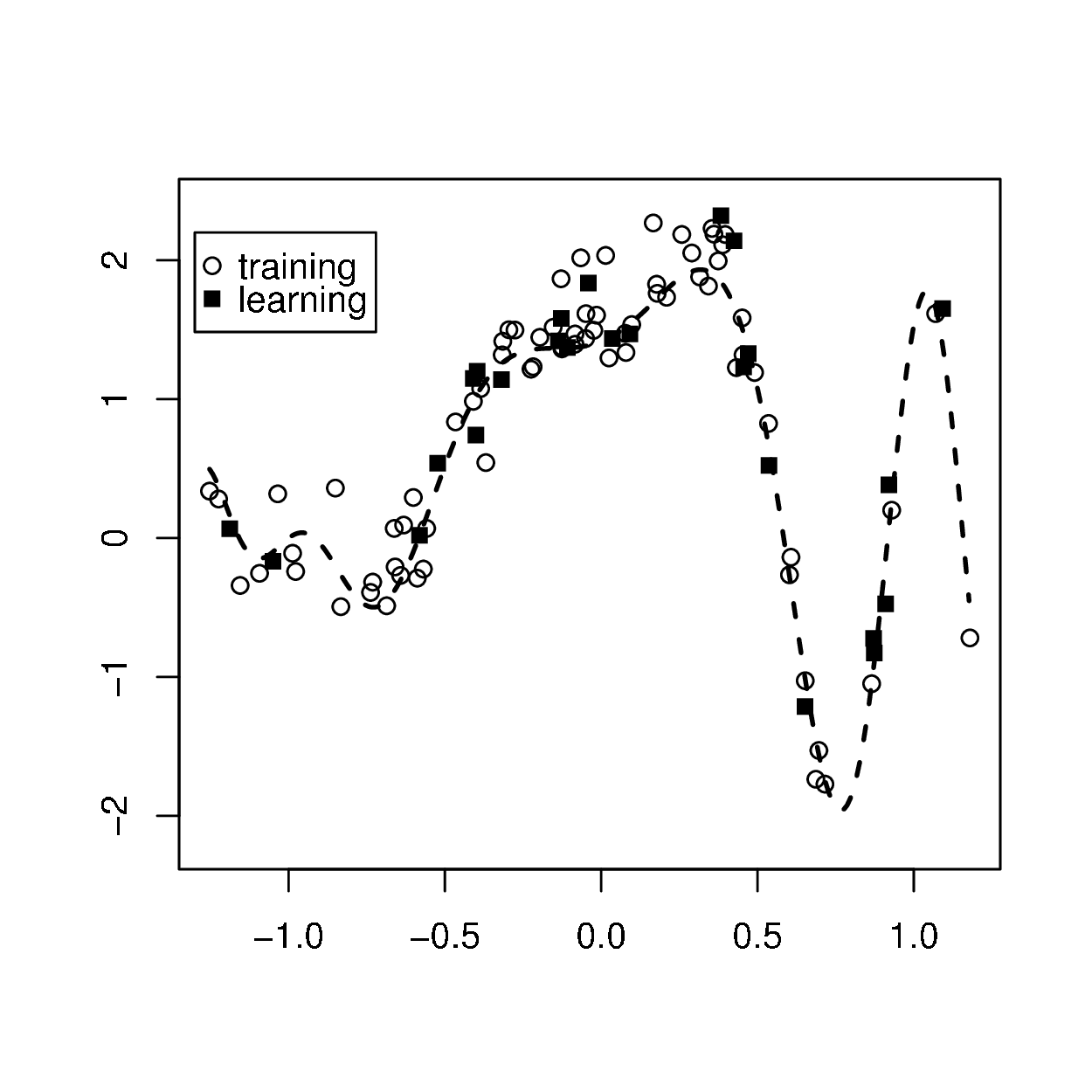}
  \includegraphics[width = 6cm]{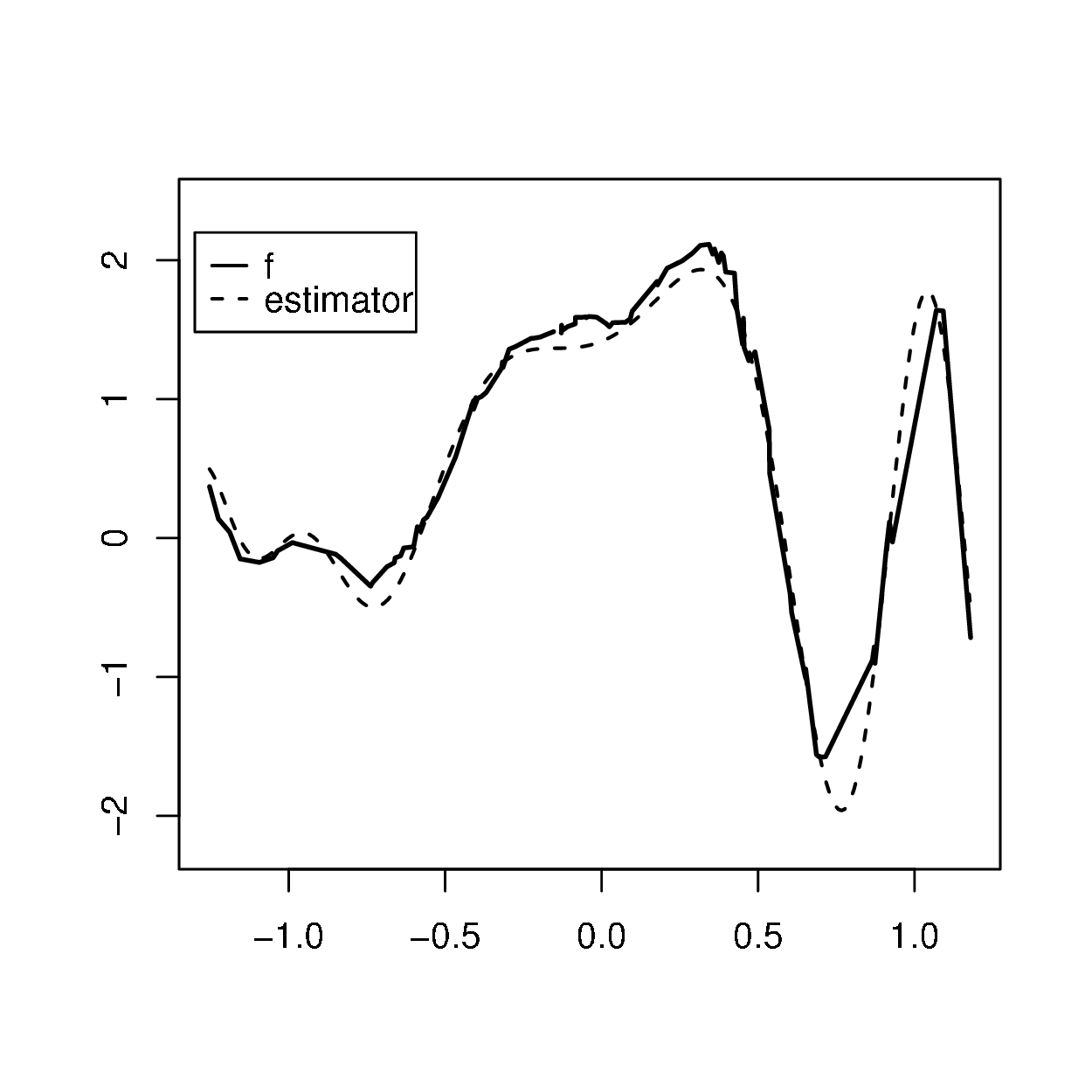}
  \caption{Simulated datasets and aggregated estimators with
    cross-validated temperature for $f = \mathtt{hardsine}$, $n=100$,
    and indexes $\vartheta = (1 / \sqrt{2}, 1 / \sqrt{2})$, $\vartheta
    = (2 / \sqrt{14}, 1/\sqrt{14}, 3/\sqrt{14})$,
    $\vartheta=(1/\sqrt{21}, -2/\sqrt{21}, 0, 4/\sqrt{21} )$ from top
    to bottom.}
  \label{fig:n100hardsine}
\end{figure}

\begin{figure}[htbp]
  \includegraphics[width = 6cm]{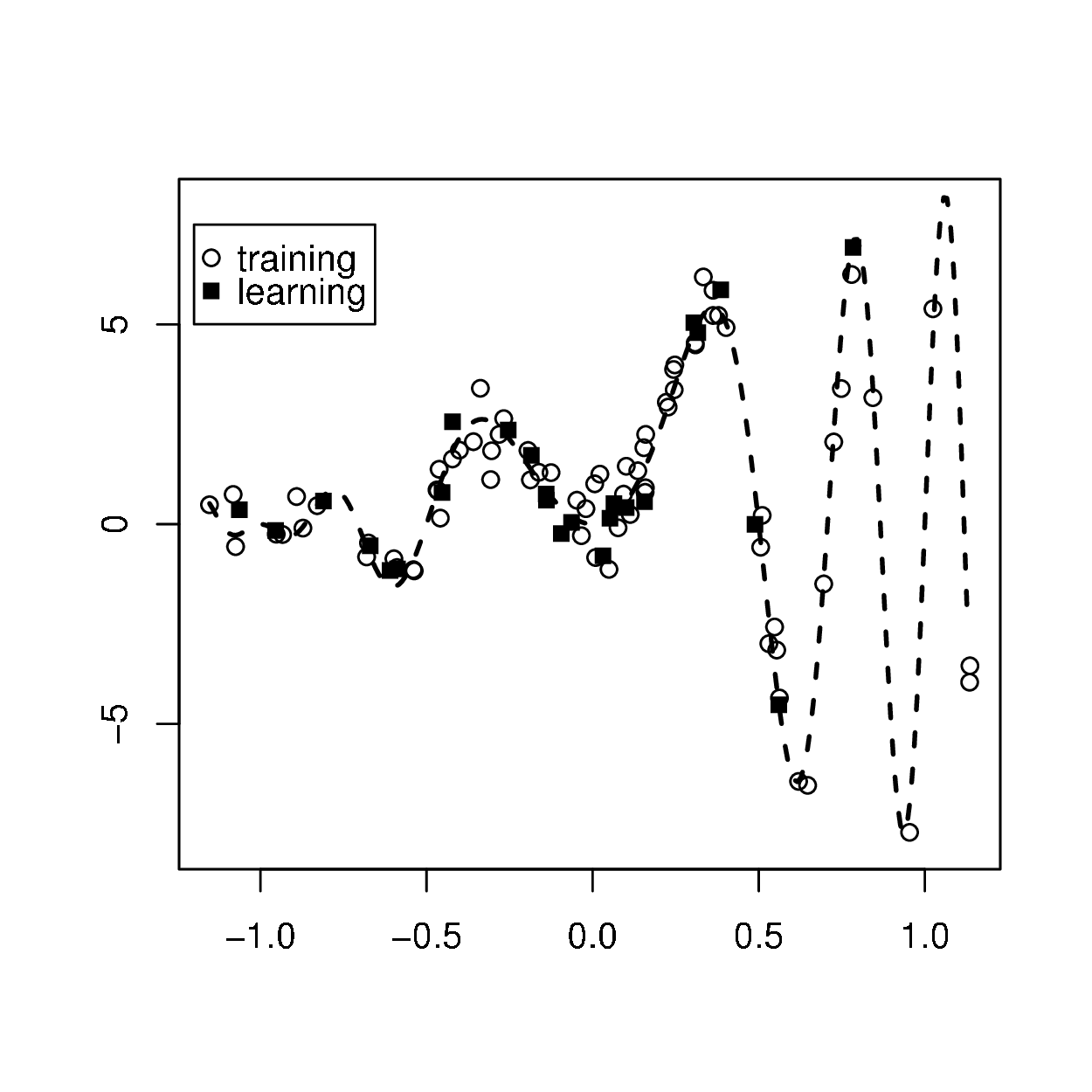}
  \includegraphics[width = 6cm]{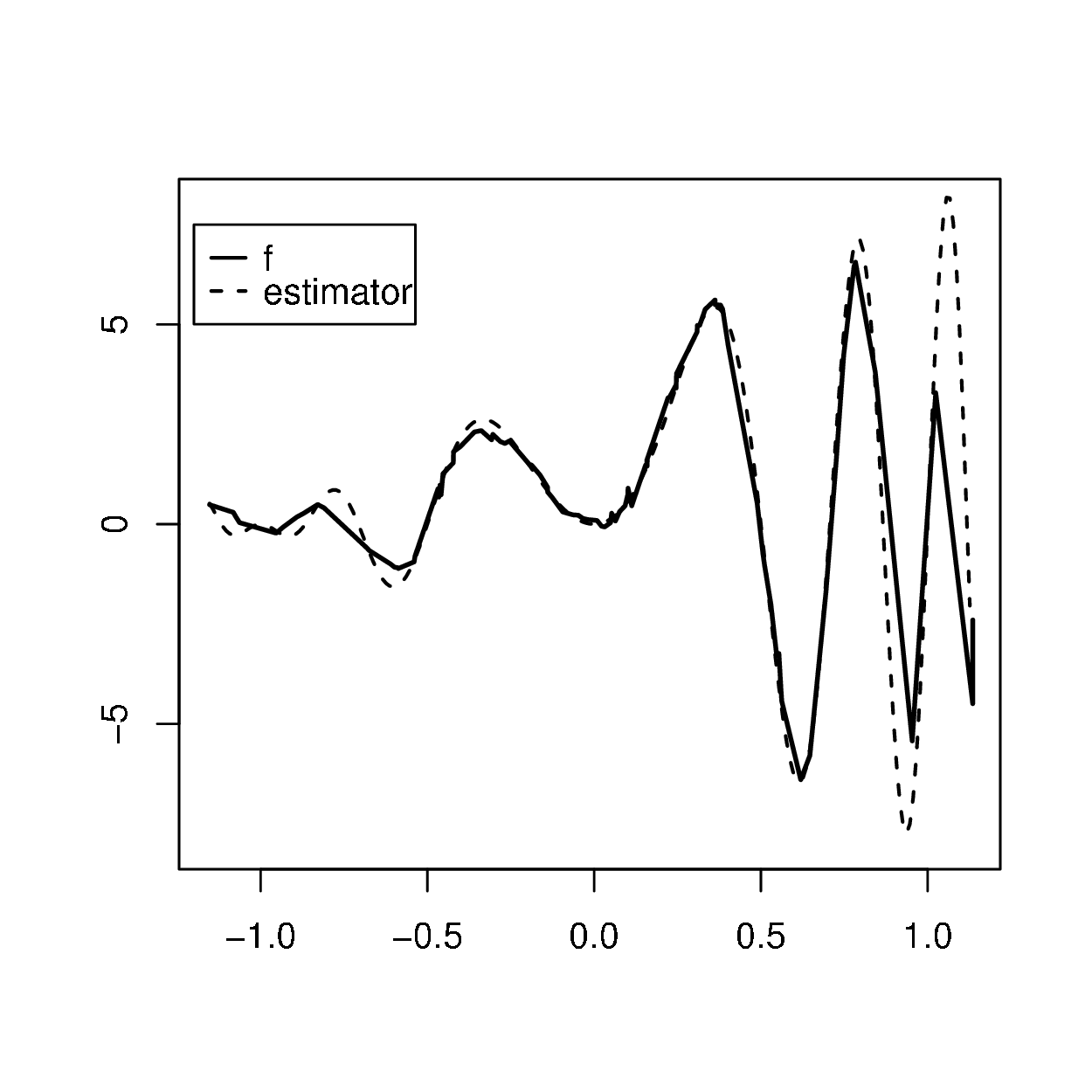}
  \includegraphics[width = 6cm]{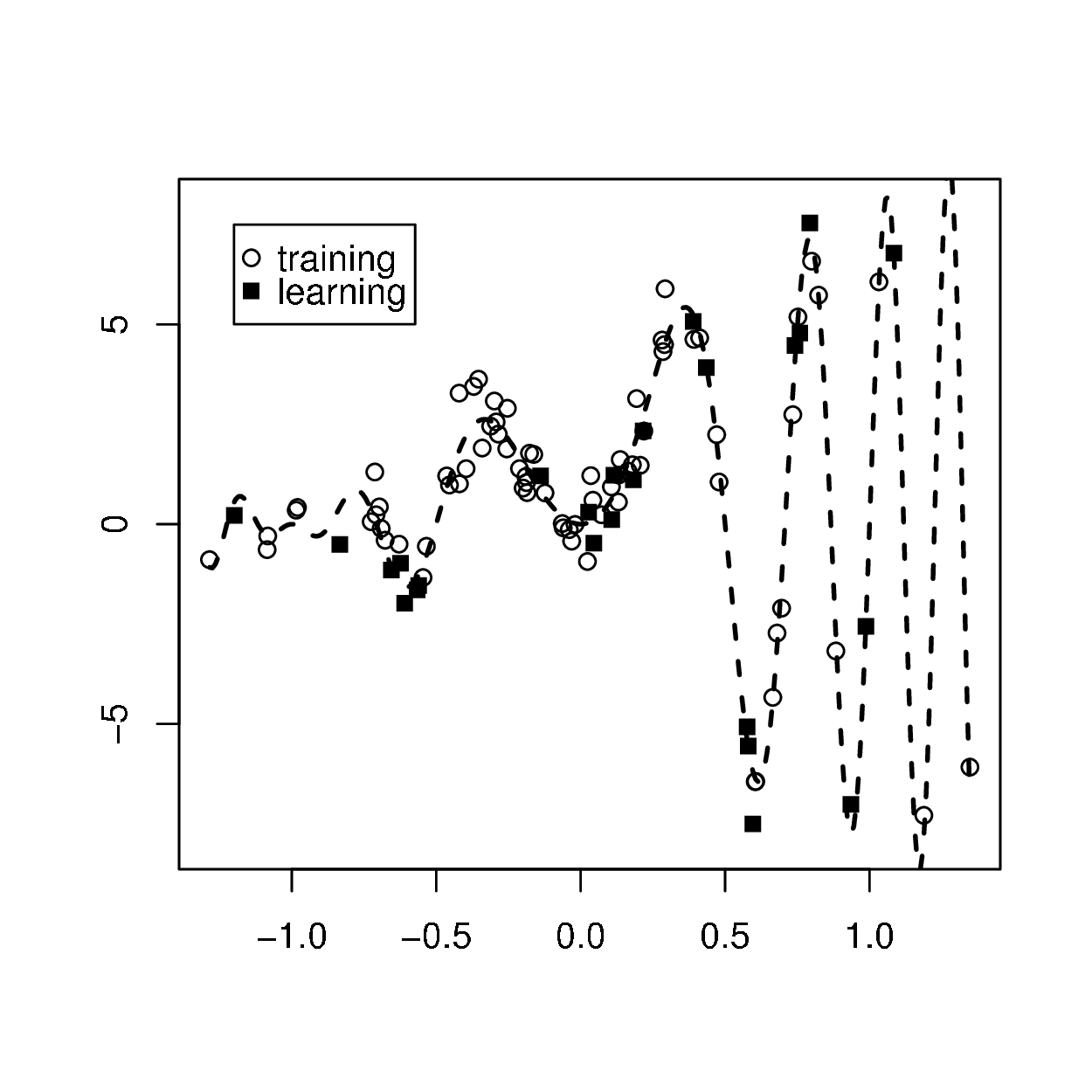}
  \includegraphics[width = 6cm]{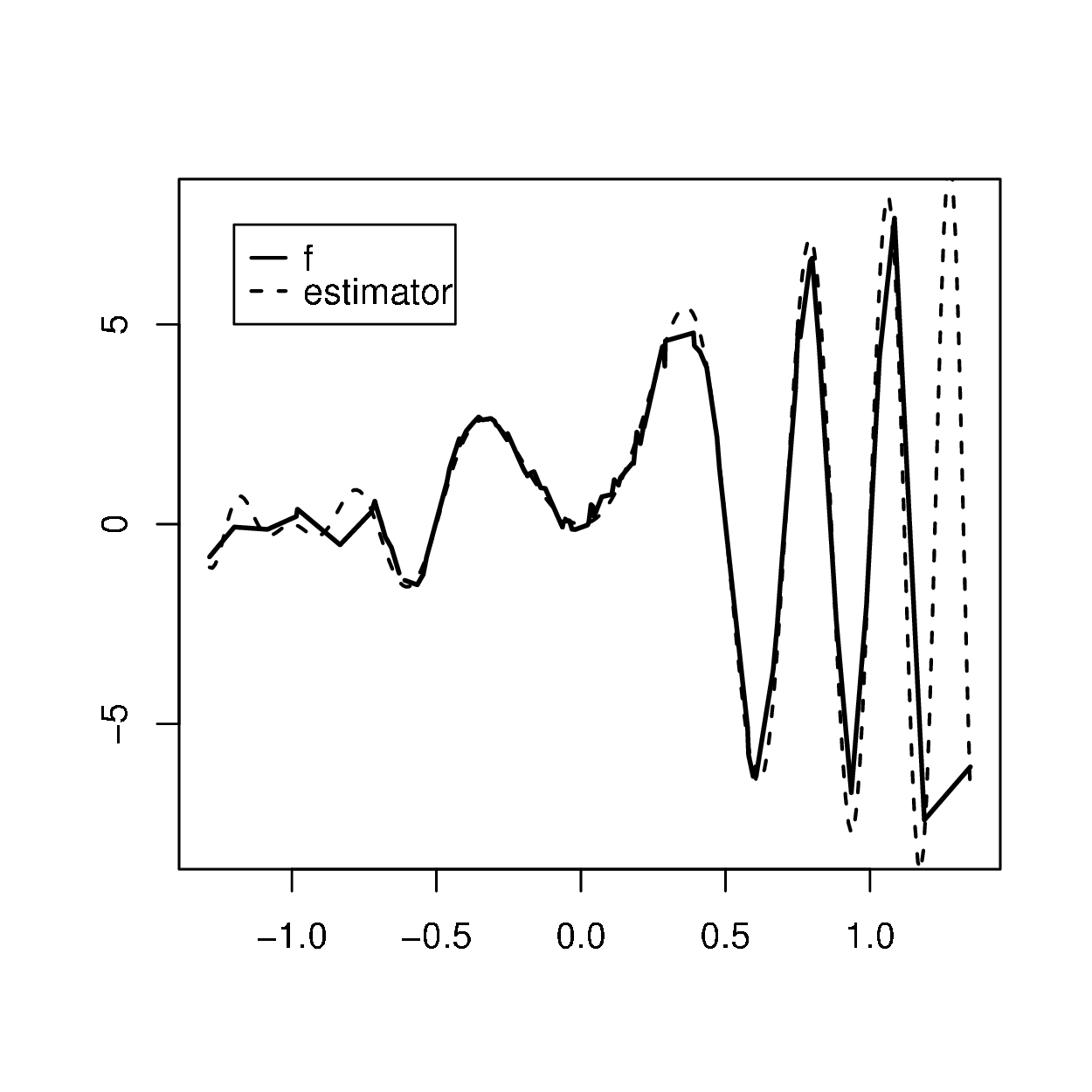}
  \includegraphics[width = 6cm]{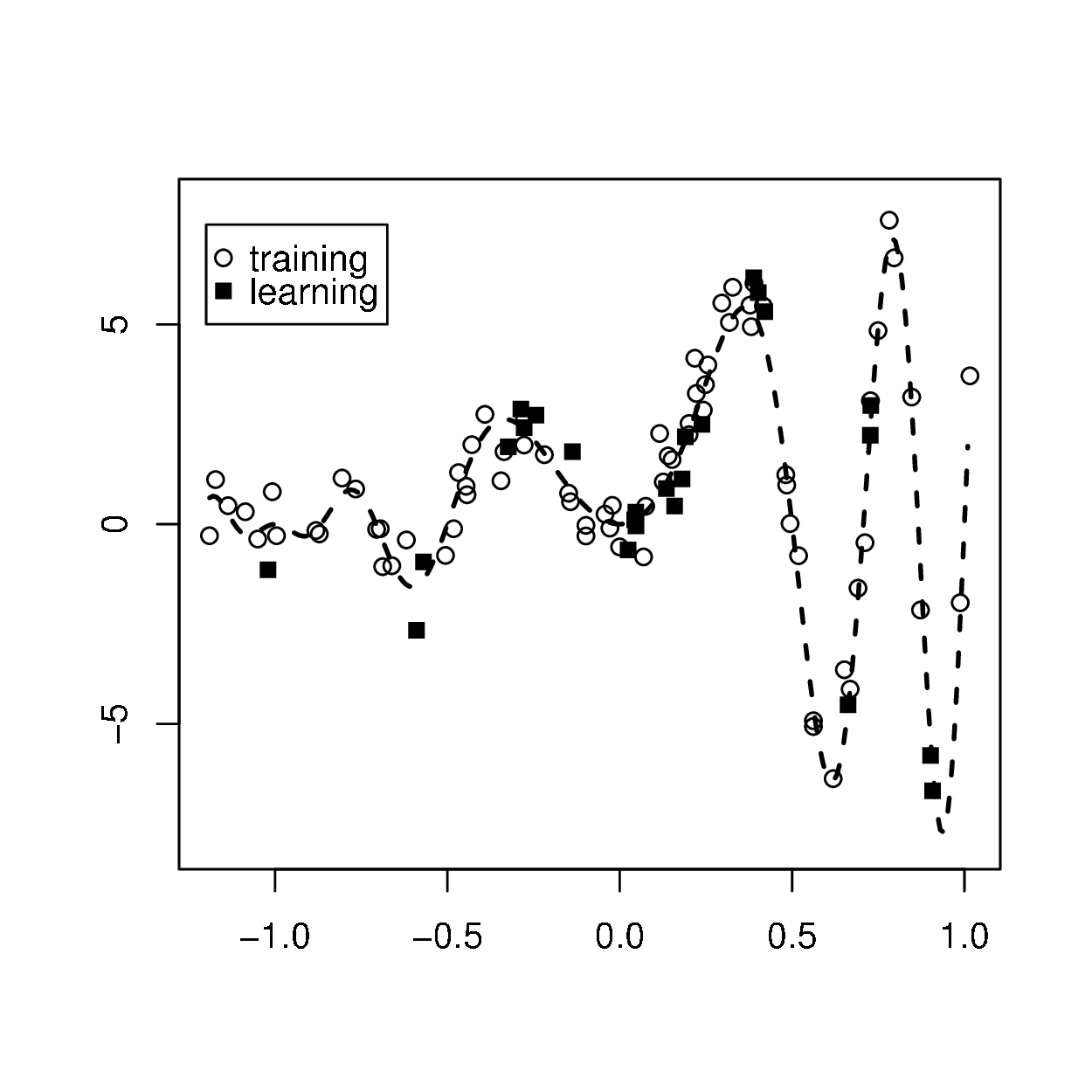}
  \includegraphics[width = 6cm]{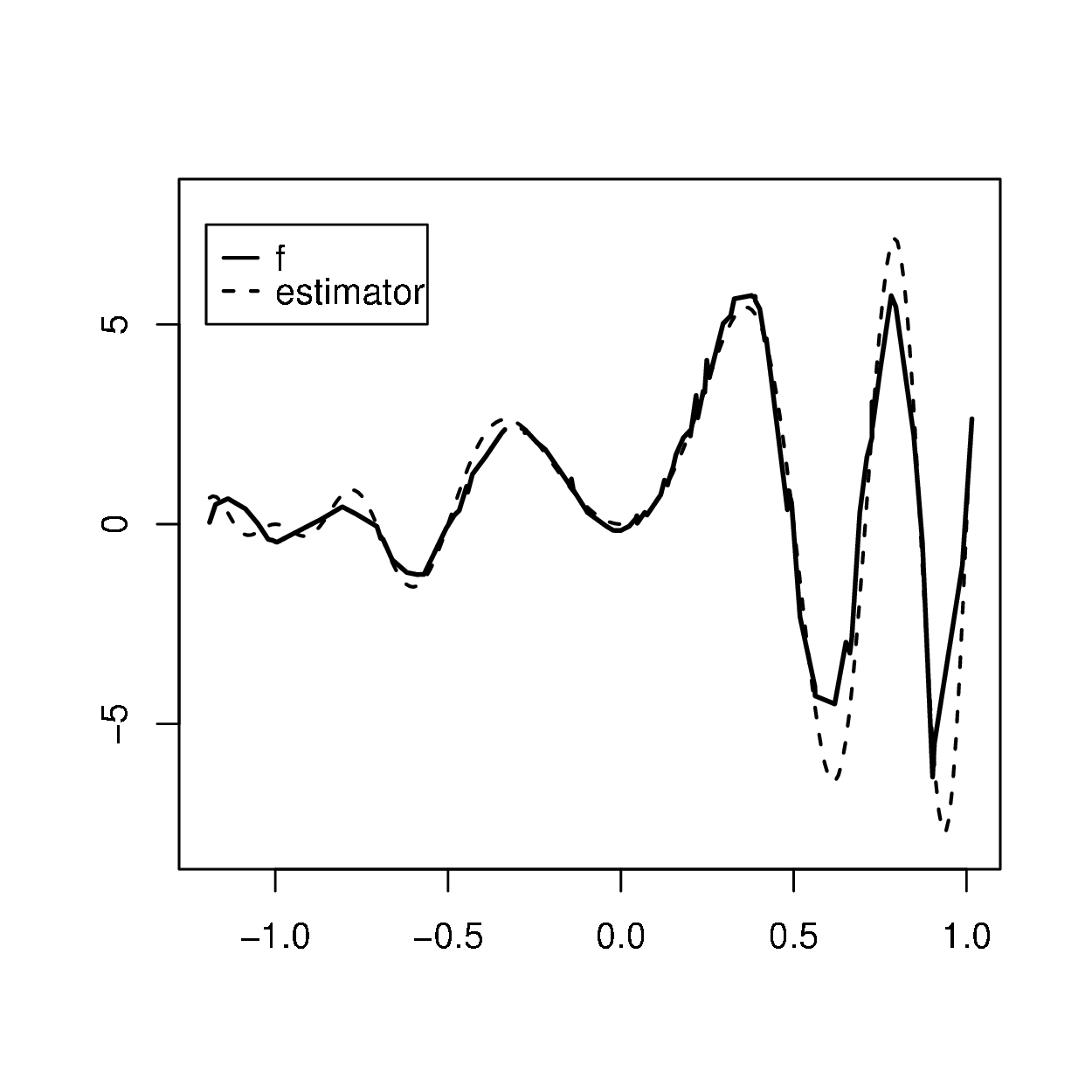}
  \caption{Simulated datasets and aggregated estimators with
    cross-validated temperature for $f = \mathtt{oscsine}$, $n=100$,
    and indexes $\vartheta = (1 / \sqrt{2}, 1 / \sqrt{2})$, $\vartheta
    = (2 / \sqrt{14}, 1/\sqrt{14}, 3/\sqrt{14})$,
    $\vartheta=(1/\sqrt{21}, -2/\sqrt{21}, 0, 4/\sqrt{21} )$ from top
    to bottom.}
  \label{fig:n100oscsine}
\end{figure}

We implemented the procedure described in Section~\ref{sec:estimator}
using the \texttt{R} software (see
\texttt{http://www.r-project.org/}). In order to increase computation
speed, we implemented the computation of local polynomials and the
bandwidth selection~\eqref{eq:bandwidth} in \texttt{C} language. The
simulated samples satisfy~\eqref{eq:model},\eqref{eq:SI}, where the
noise is centered Gaussian with homoscedastic variance
\begin{equation*}
  \sigma = \big[ \sum_{1 \leq i \leq n} f(\vartheta\T X_i)^2 / (n
  \times \mathtt{rsnr}) \big]^{1/2},
\end{equation*}
where $\mathtt{rsnr} = 5$. This choice of $\sigma$ makes the
root-signal-to-noise ratio, which is a commonly used assessment of the
complexity of estimation, equals to~$5$. We consider the following
link functions (see the dashed lines in Figures~\ref{fig:n100hardsine}
and~\ref{fig:n100oscsine}):
\begin{itemize}
\item $\mathtt{oscsine}(x) = 4 (x+1) \sin( 4 \pi x^2 )$,
\item $\mathtt{hardsine}(x) = 2 \sin(1 + x) \sin( 2 \pi x^2 + 1)$.
\end{itemize}
The simulations are done with a uniform design on $[-1, 1]^d$, with
dimensions $d \in \{ 2, 3, 4 \}$ and we consider several indexes
$\vartheta$ that make $P_{\vartheta\T X}$ not uniform.

In all the computations below, the parameters for the procedure are
$\Lambda = \hat S \times G \times \mathcal L$ where $\hat S$ is
computed using the algorithm described in
Section~\ref{sec:reduccomplex} and where $G = \{ 1, 2, 3, 4 \}$ and
$\mathcal L = \{ 0.1, 0.5, 1, 1.5 \}$. The degree of the local
polynomials is $r = 5$. The learning sample has size $[n / 4]$, and is
chosen randomly in the whole sample. We do not use a jackknife
procedure (that is, the average of estimators obtained with several
learning subsamples), since the results are stable enough (at least
when $n \geq 100$) when we consider only one learning sample.

In Tables~\ref{tab:misehardsined2}, \ref{tab:misehardsined3},
\ref{tab:misehardsined4} and Figures \ref{fig:misehardsined2},
\ref{fig:misehardsined3}, \ref{fig:misehardsined4}, we show the mean
MISE for 100 replications and its standard deviation for several Gibbs
%
\begin{table*}[h]
  \caption{MISE against the Gibbs temperature \textup($f= \mathtt{hardsine}$,
    $d=2$, $\vartheta=(1/\sqrt{2}, 1/\sqrt{2})$.\textup)}
  \begin{tabular}{lccccccccc}
    \hline
    Temperature & 0.1 &  0.5 &  0.7 & 1.0 & 1.5 & 2.0 & ERM &
    \texttt{aggCVT} \\
    \hline
    n = 100 & 0.026 & 0.017 & 0.015 & \textbf{0.014} & \textbf{0.014}
    & 0.015 & 0.034 & 0.015 \\
    & (.009) & (.006) & (.006) & (.005) & (.005) & (.006) &
    (.018) & (.005) \\
    n = 200 & 0.015 & 0.009 & \textbf{0.008} & \textbf{0.008} & 0.009 &
    0.011 & 0.027 & 0.009 \\
    & (.004) & (.002) & (.003) & (.003) & (.005) & (.007) & (.014) &
    (.004) \\
    n = 400 & 0.006 & 0.005 & \textbf{0.004} & 0.005 & 0.006 & 0.007 &
    0.016 & 0.005  \\
    & (.001) & (.001) & (.001) & (.001) & (.002) & (.002) & (.003) &
    (.002) \\
    \hline
  \end{tabular}
  \label{tab:misehardsined2}
\end{table*}
%
\begin{table*}[h]
  \caption{MISE against the Gibbs temperature \textup($f = \mathtt{hardsine}$,
    $d=3$, $\vartheta=(2 / \sqrt{14}, 1/\sqrt{14}, 3/\sqrt{14})$\textup).}
  \begin{tabular}{lcccccccc}
    \hline
    Temperature & 0.1 &  0.5 &  0.7 & 1.0 & 1.5 & 2.0 & ERM &
    \texttt{aggCVT} \\
    \hline
    n = 100 & 0.029 & 0.021 & 0.019 & 0.018 & \textbf{0.017} & 0.018 &
    0.037 & 0.020 \\
    & (.011) & (.008) & (.008) & (.007) & (.008) & (.009) & (.022) &
    (.008) \\
    n = 200 & 0.016 & 0.010 & 0.010 & \textbf{0.009} & \textbf{0.009}
    & 0.010 & 0.026 & 0.010 \\
    & (.005) & (.003) & (.003) & (.002) & (.002) & (.003) & (0.008) &
    (.003) \\
    n = 400 & 0.007 & 0.006 & \textbf{0.005} & \textbf{0.005} & 0.006
    & 0.007 & 0.017 & 0.006 \\
    & (.002) & (.001) & (.001) & (.001) & (.001) & (.002) & (.003) &
    (.001) \\
    \hline
  \end{tabular}
  \label{tab:misehardsined3}
\end{table*}
%
%
\begin{table*}[h]
  \caption{MISE against the Gibbs temperature \textup($f= \mathtt{hardsine}$,
    $d=4$, $\vartheta=(1/\sqrt{21}, -2/\sqrt{21}, 0, 4/\sqrt{21})$\textup)}
  \begin{tabular}{lcccccccc}
    \hline
    Temperature & 0.1 &  0.5 &  0.7 & 1.0 & 1.5 & 2.0 & ERM &
    \texttt{aggCVT} \\
    \hline
    n = 100 & 0.038 & 0.027 & 0.021 & 0.019 & \textbf{0.017} &
    \textbf{0.017} & 0.038 & 0.020 \\
    &  (.016) & (.010) & (.009) & (.008) & (.007) & (.007) & (.025) &
    (.010) \\
    n = 200 & 0.019 & 0.013 & \textbf{0.012} & \textbf{0.012} & 0.013
    & 0.014 & 0.031  & 0.013 \\
    & (.014) & (.009) & (.010) & (.011) & (.012) & (.012) & (.016)
    & (.010) \\
    n = 400 &  0.009 & 0.006 & \textbf{0.005} & \textbf{0.005} & 0.006
    & 0.007 & 0.017 & 0.006 \\
    & (.002) & (.001) & (.001) & (.001) & (.001) & (.002) & (.004) &
    (.001) \\
    \hline
  \end{tabular}
  \label{tab:misehardsined4}
\end{table*}
temperatures, several sample sizes and indexes. These results
empirically prove that the aggregated estimator outperforms the ERM
(which is computed as the aggregated estimator with a large
temperature $T = 30$) since in each case, the aggregated estimator
with cross-validated temperature (\texttt{aggCVT}, given
by~\eqref{CVT}, with $\mathcal T = \{ 0.1, 0.2, \ldots, 4.9, 5 \}$),
has a MISE much smaller than the MISE of the ERM. Moreover,
\texttt{aggCVT} is more stable than the ERM in view of the standard
deviations (in brackets). Note also that as expected, the dimension
parameter has no impact on the accuracy of estimation: the MISEs are
barely the same when $d=2, 3, 4$.

\setlength{\figlen}{4.6cm}

\begin{figure}[t]
  \includegraphics[width = \figlen]{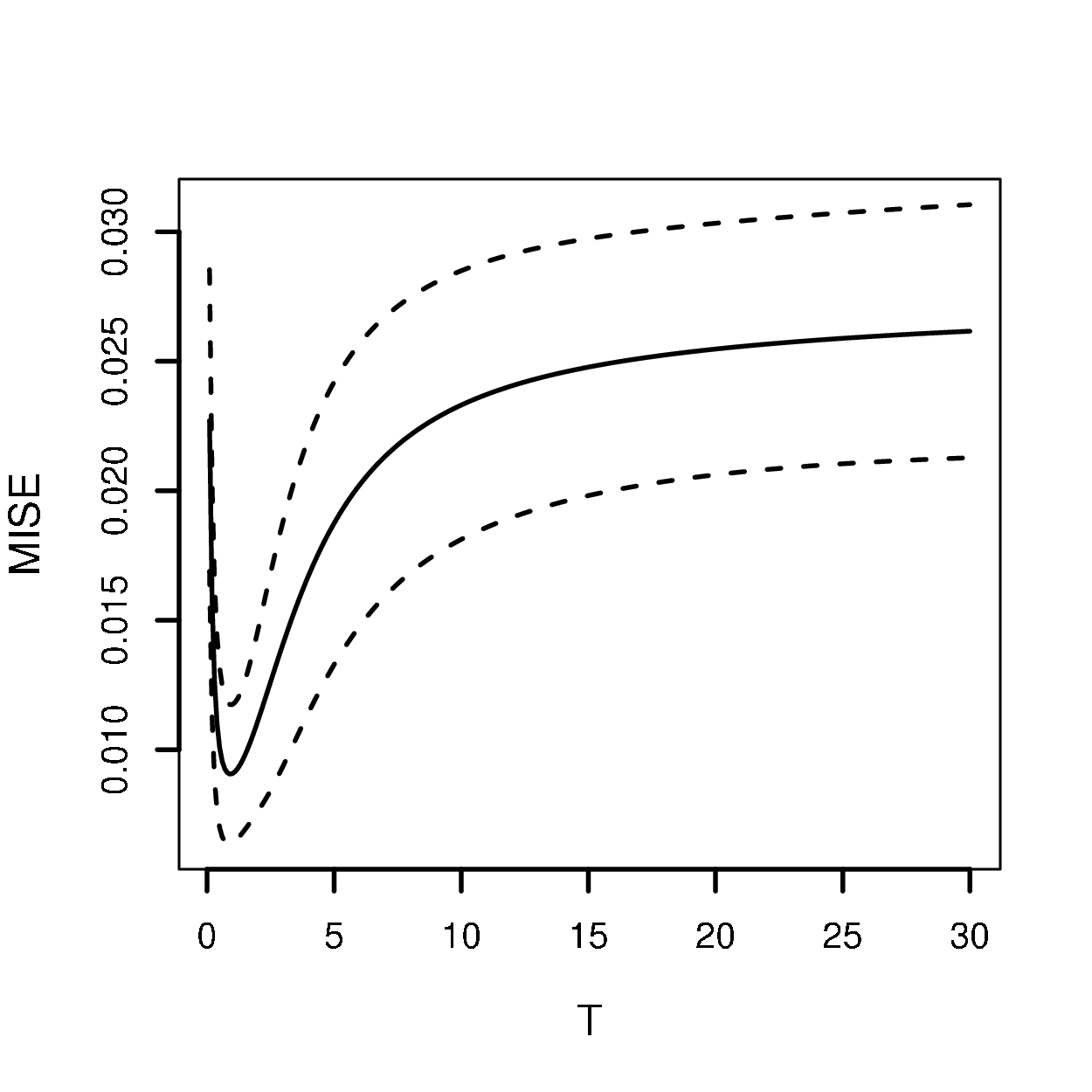}
  \includegraphics[width = \figlen]{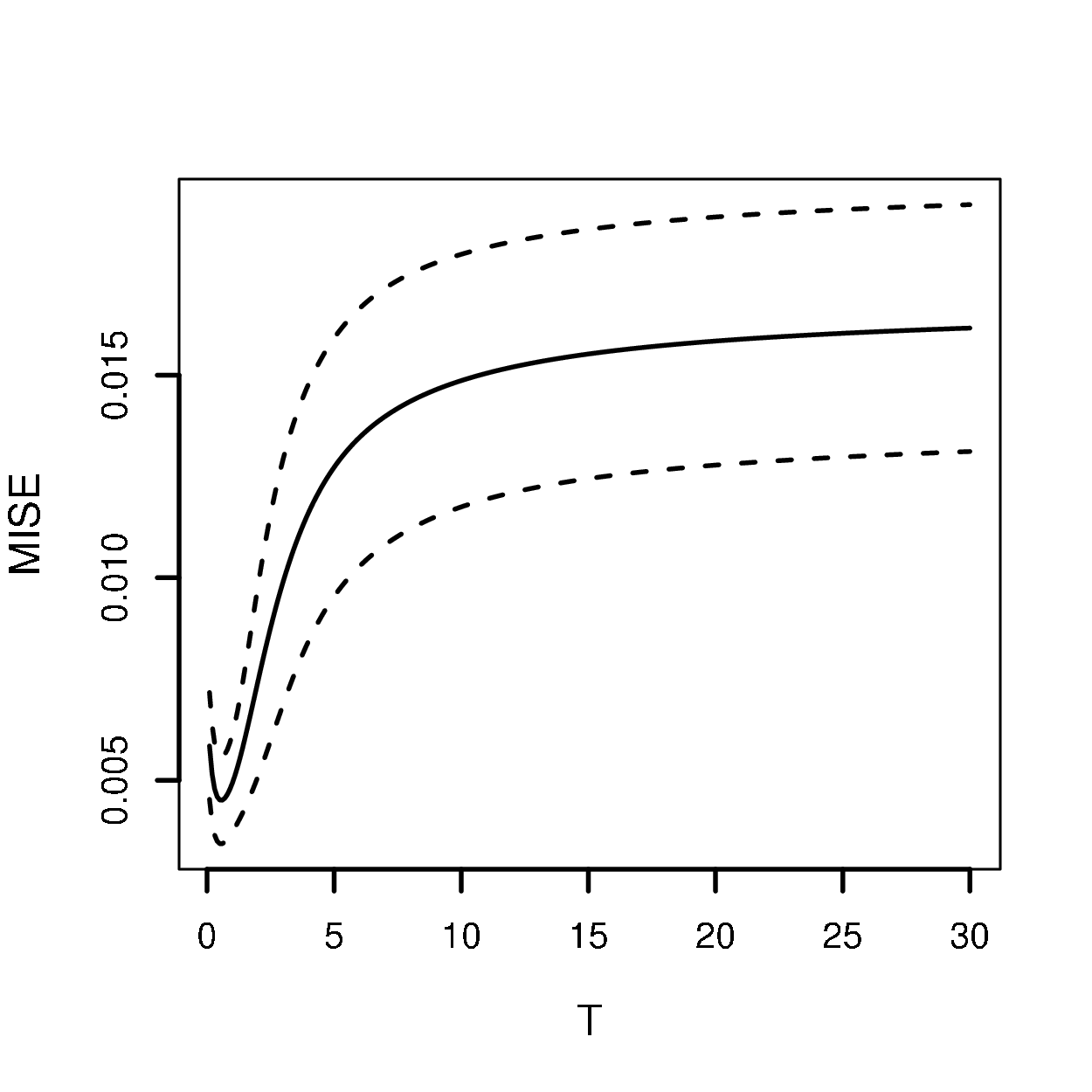}
  \caption{MISE against the Gibbs temperature for $f =
    \mathtt{hardsine}$, $\vartheta = (1 / \sqrt{2}, 1 / \sqrt{2})$, $n
    = 200, 400$ \textup(solid line = mean of the MISE for 100
    replications, dashed line = mean MISE $\pm$ standard
    deviation.\textup)}
  \label{fig:misehardsined2}
\end{figure}
\begin{figure}
  \includegraphics[width = \figlen]{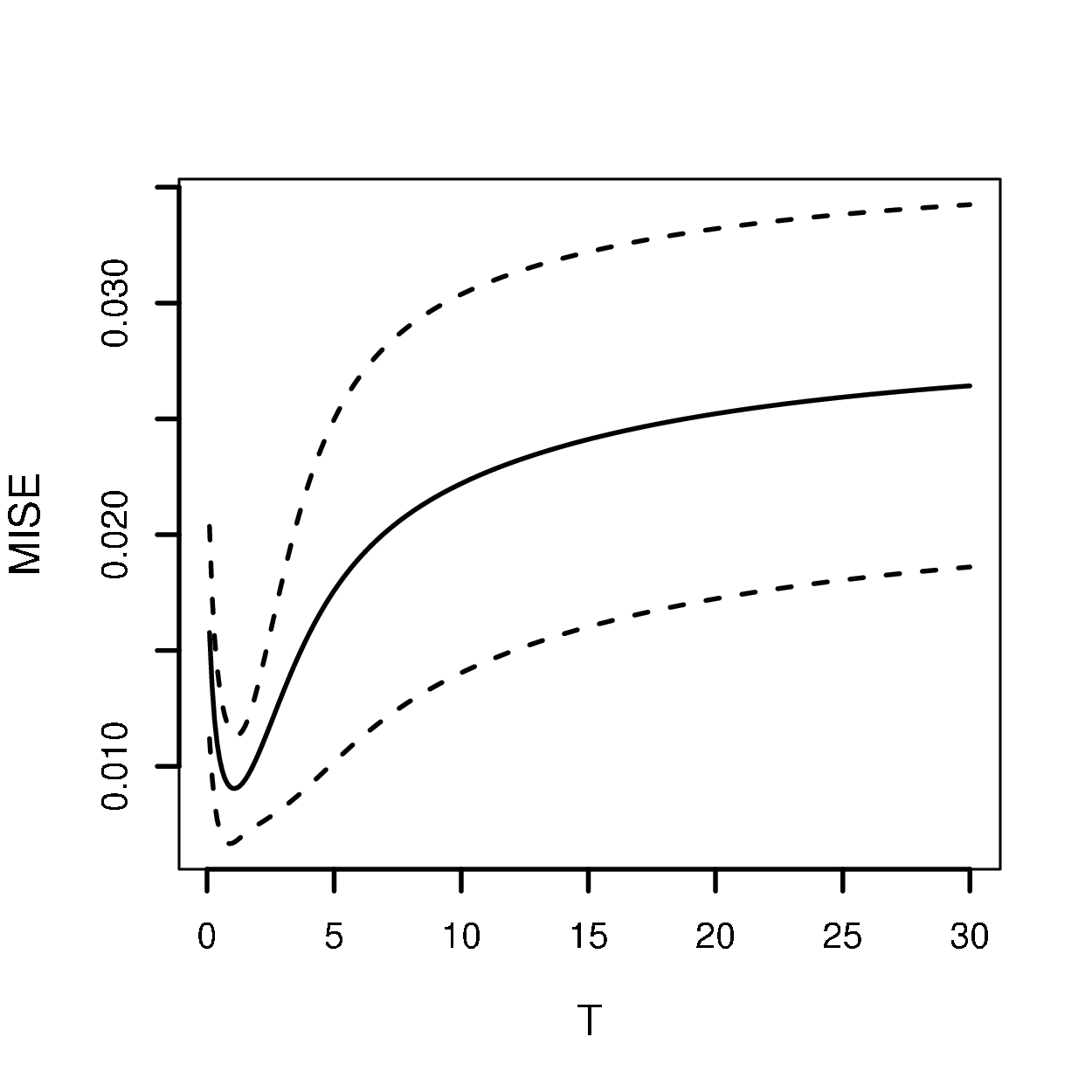}
  \includegraphics[width = \figlen]{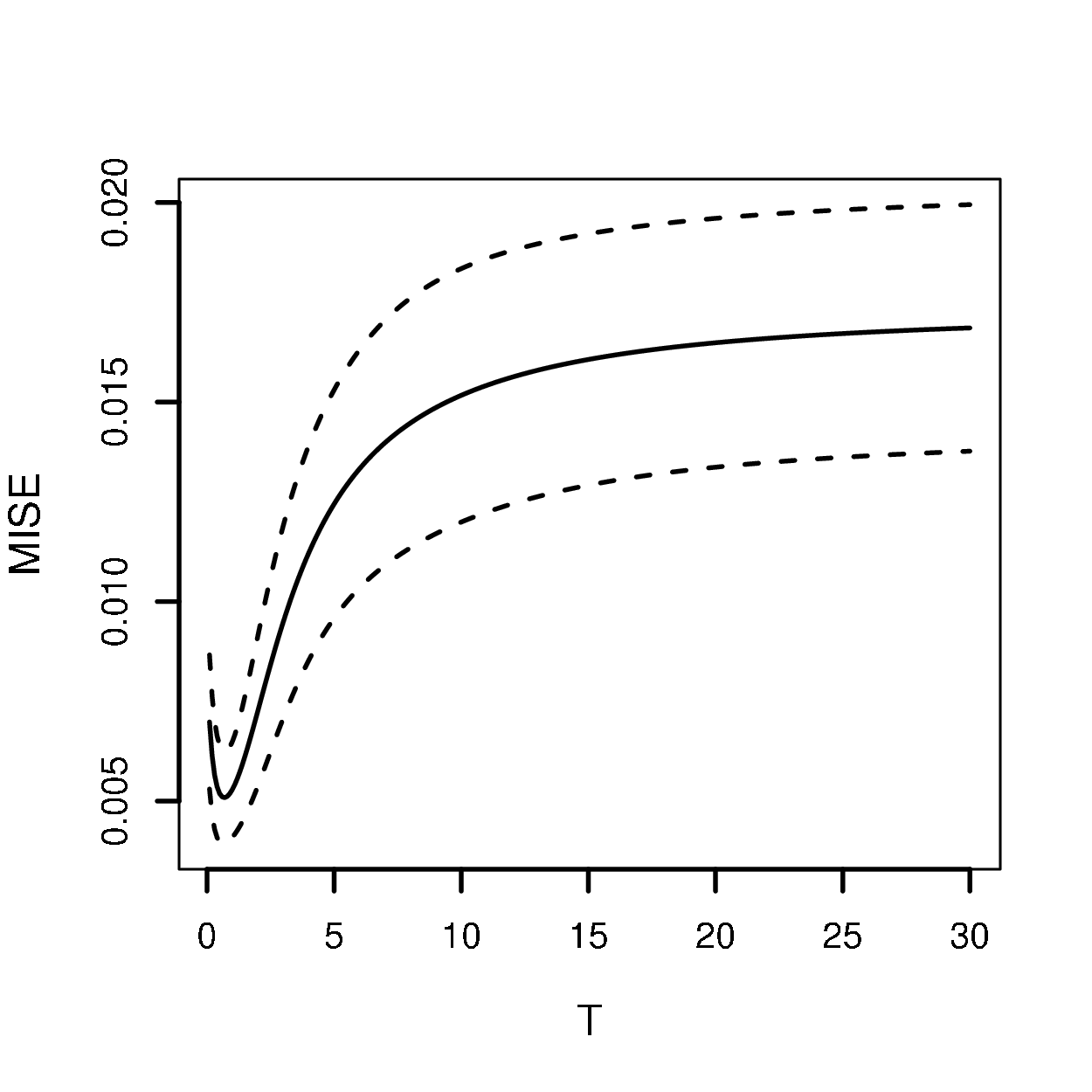}
  \caption{MISE against the Gibbs temperature for $f =
    \mathtt{hardsine}$, $\vartheta = (2 / \sqrt{14}, 1/\sqrt{14},
    3/\sqrt{14})$, $n = 200, 400$ \textup(solid line = mean of the
    MISE for 100 replications, dashed line = mean MISE $\pm$ standard
    deviation.\textup)}
  \label{fig:misehardsined3}
\end{figure}

\begin{figure}[h!]
  \vspace*{-12pt}
  \includegraphics[width = \figlen]{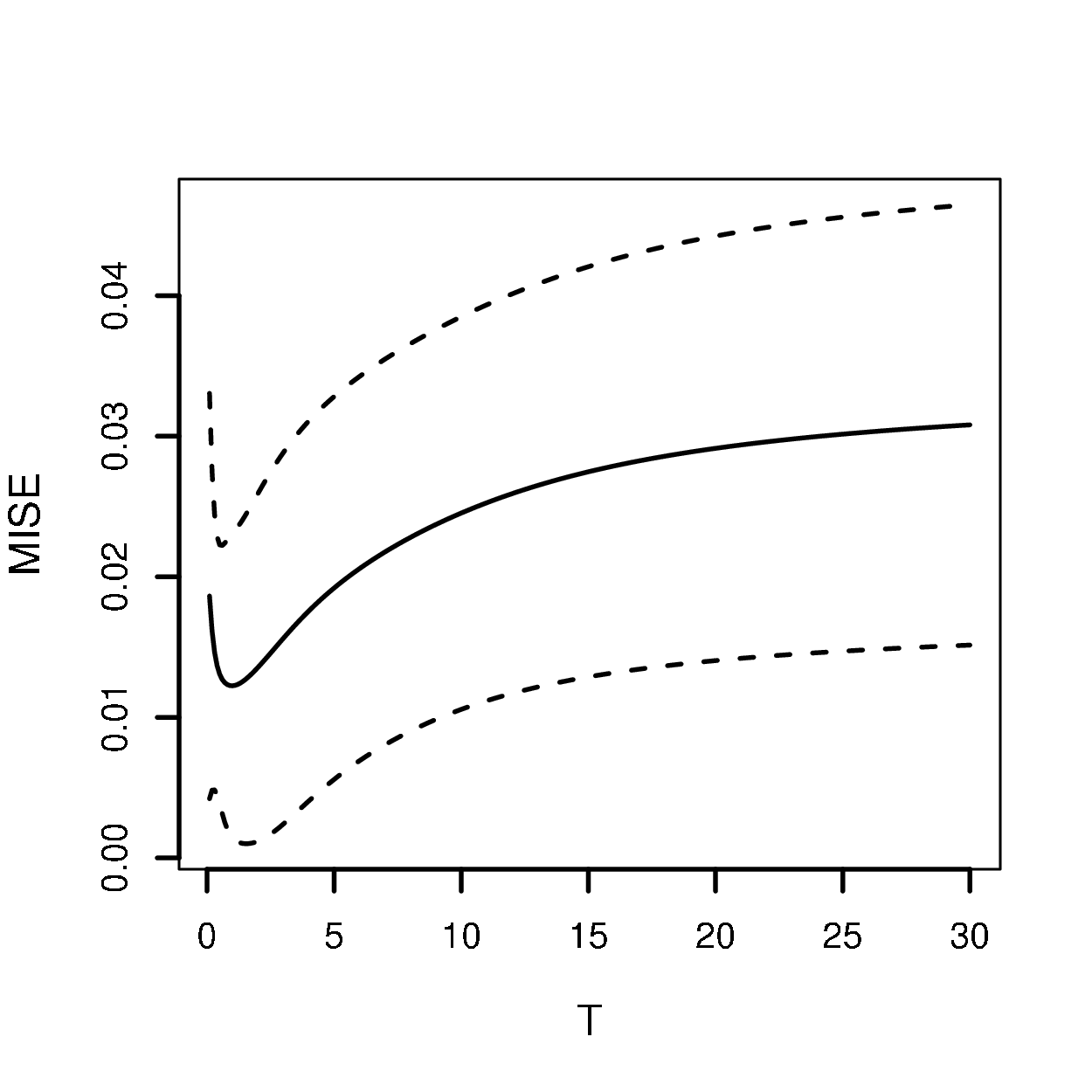}
  \includegraphics[width = \figlen]{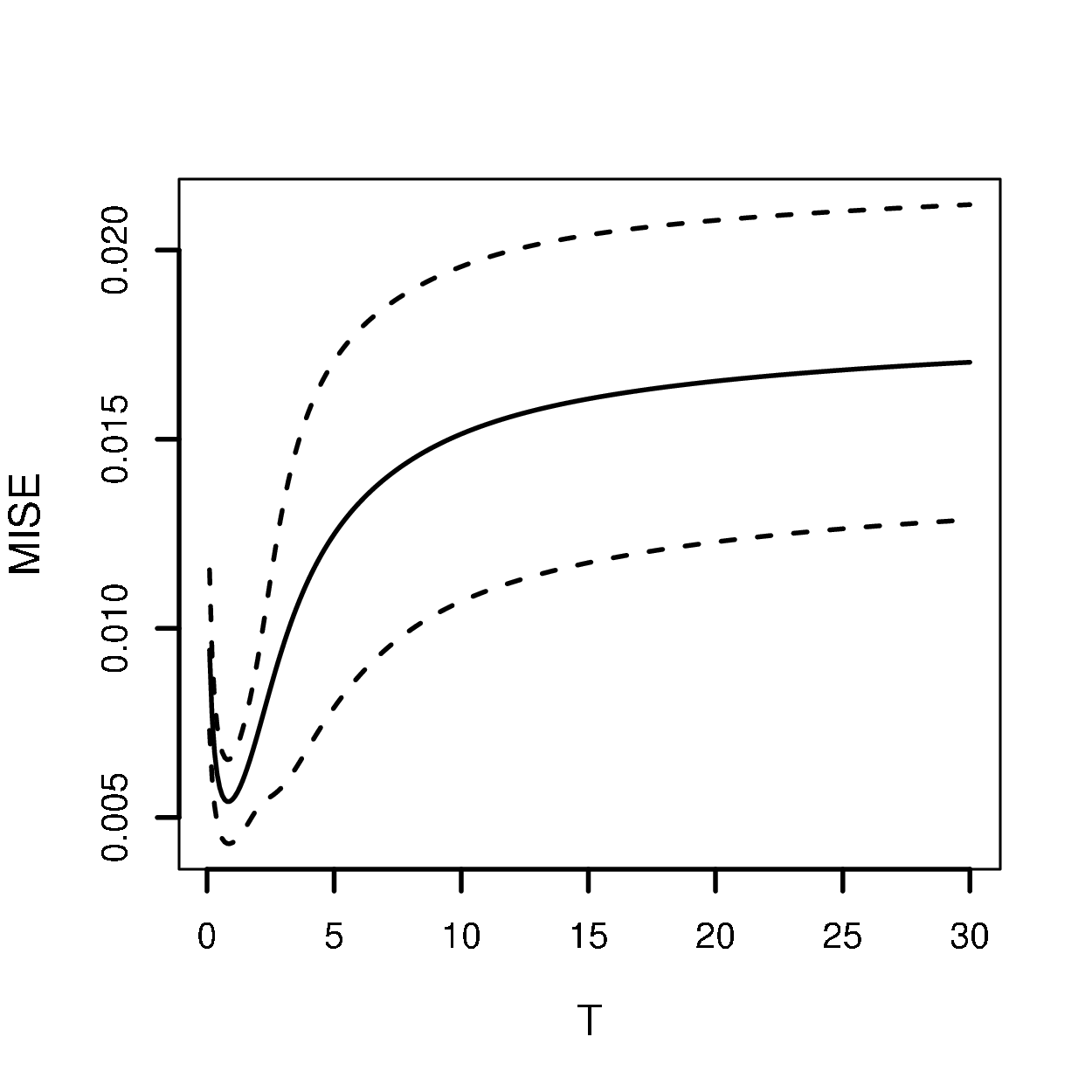}
  \vspace*{-6pt}
  \caption{MISE against the Gibbs temperature for $f =
    \mathtt{hardsine}$, $\vartheta = (2 / \sqrt{14}, 1/\sqrt{14},
    3/\sqrt{14})$, $n = 200, 400$.}
  \label{fig:misehardsined4}
  \vspace*{-12pt}
\end{figure}


\setlength{\figlen}{6cm}

\begin{figure}[h!]
  \includegraphics[width = \figlen]{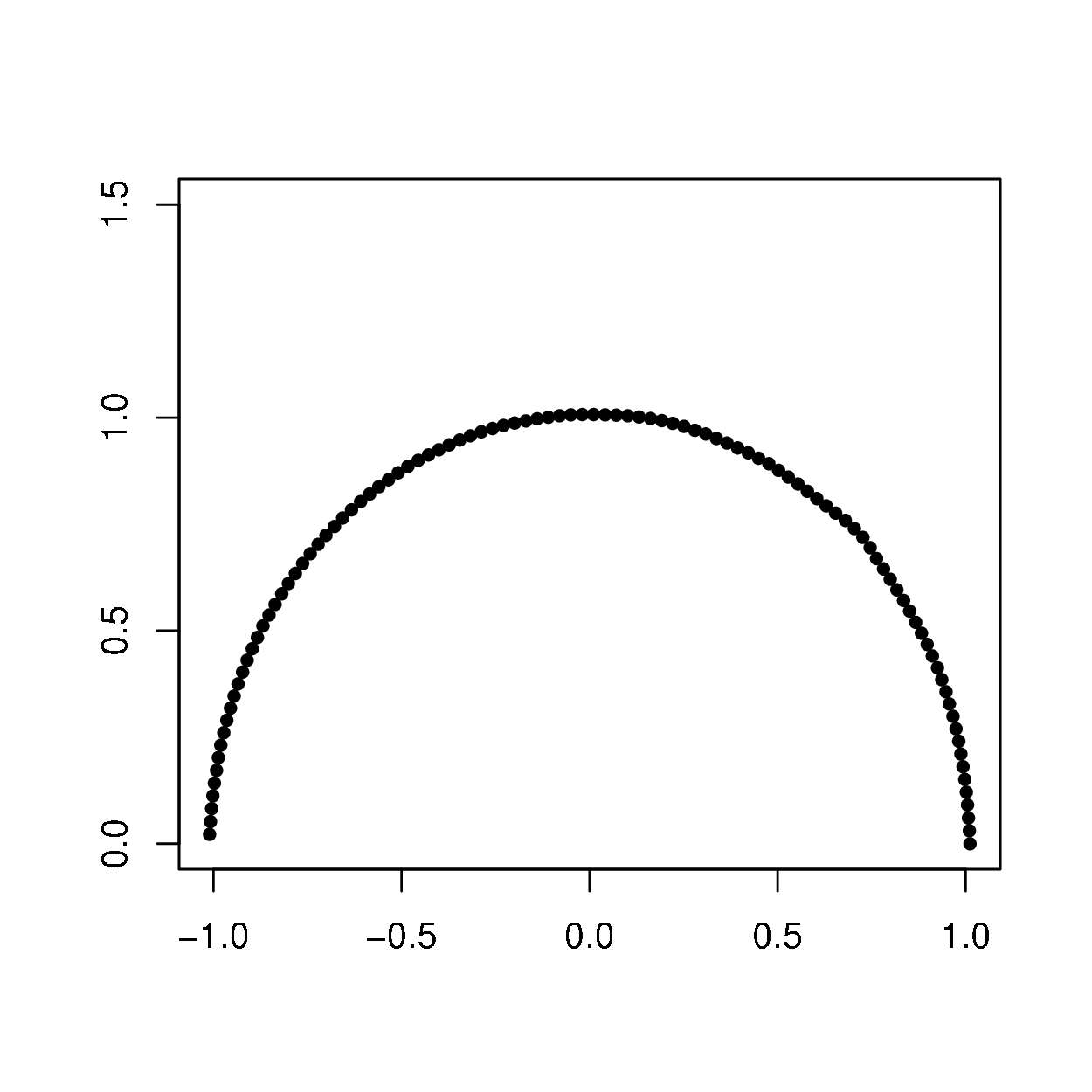}
  \includegraphics[width = \figlen]{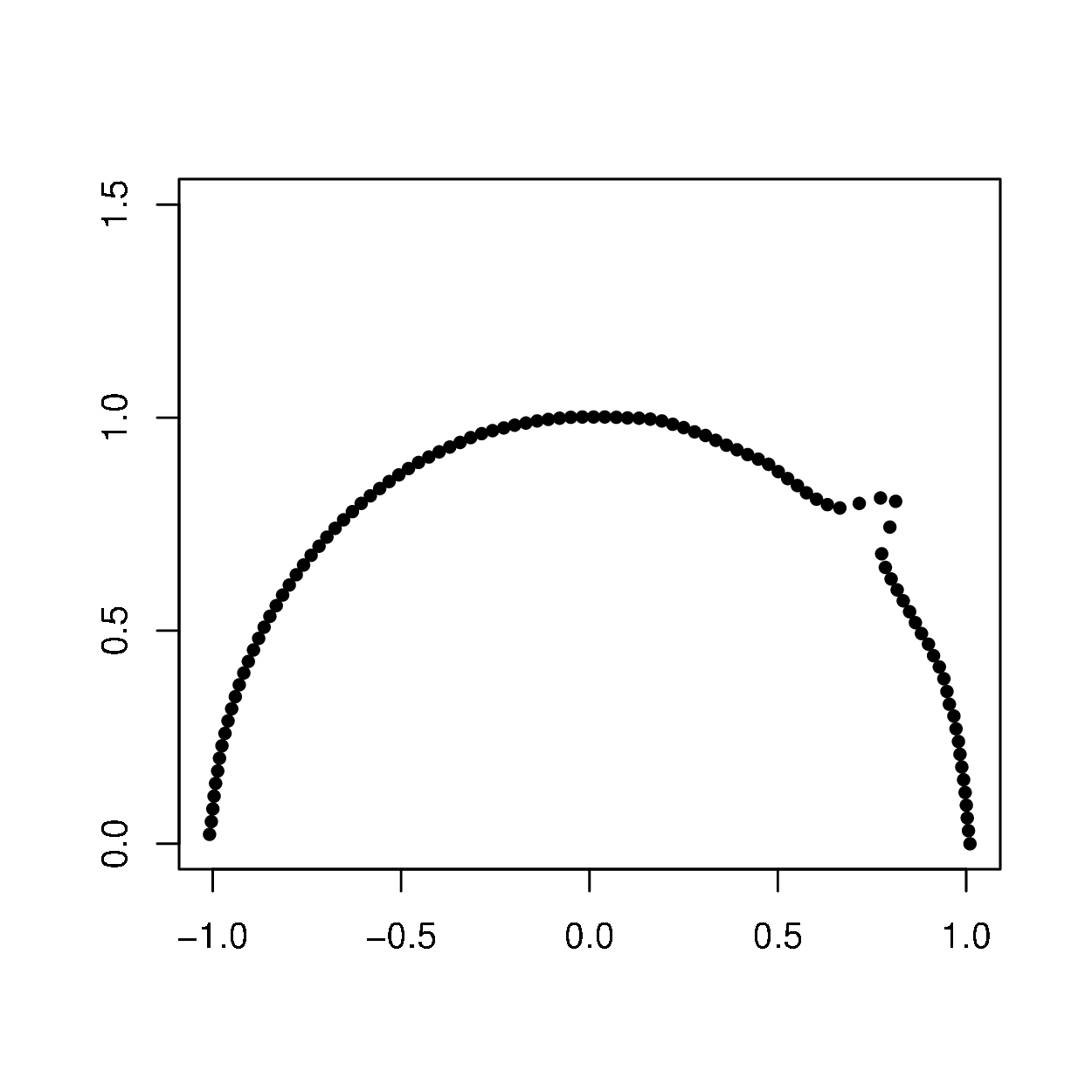} \\[-12pt]
  \includegraphics[width = \figlen]{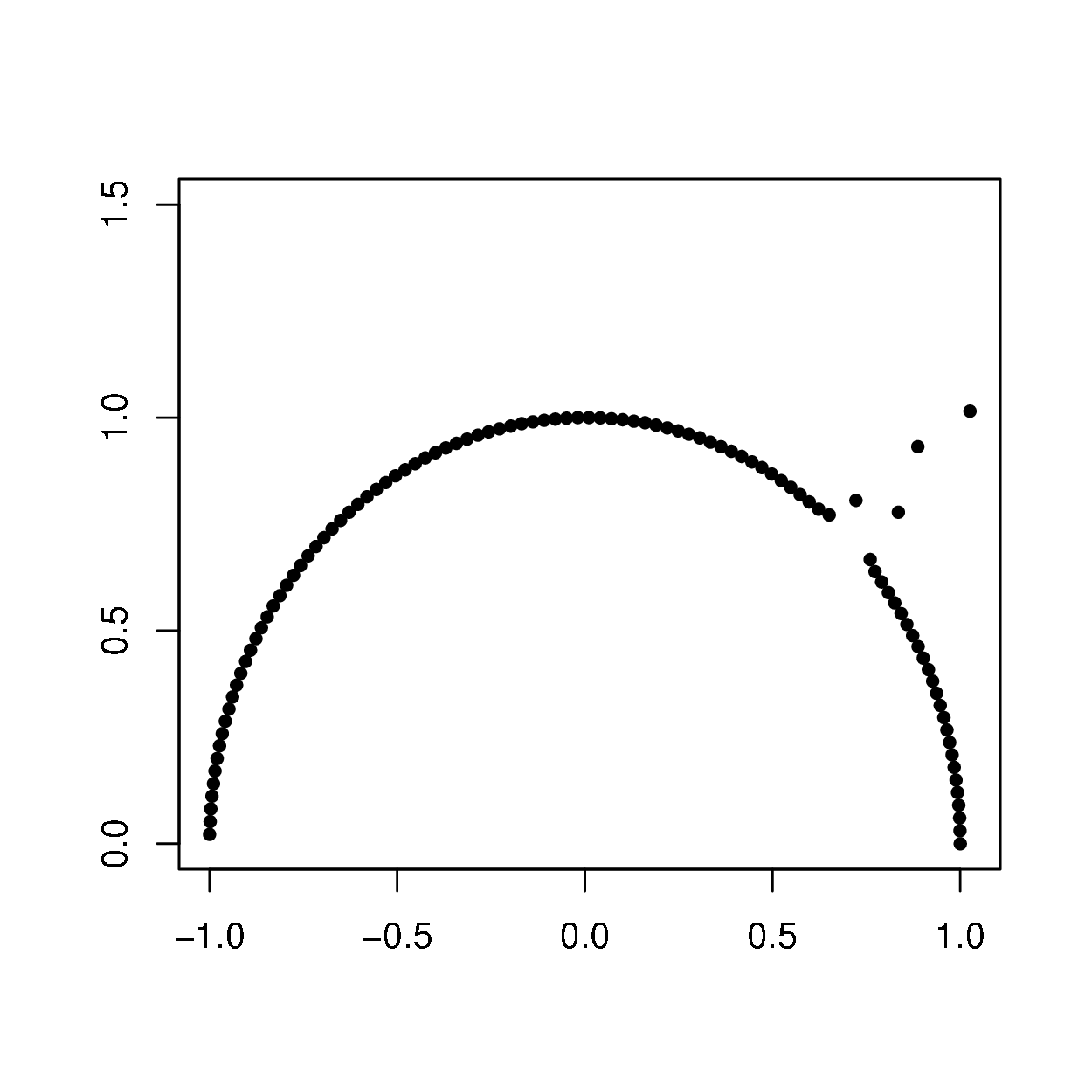}
  \includegraphics[width = \figlen]{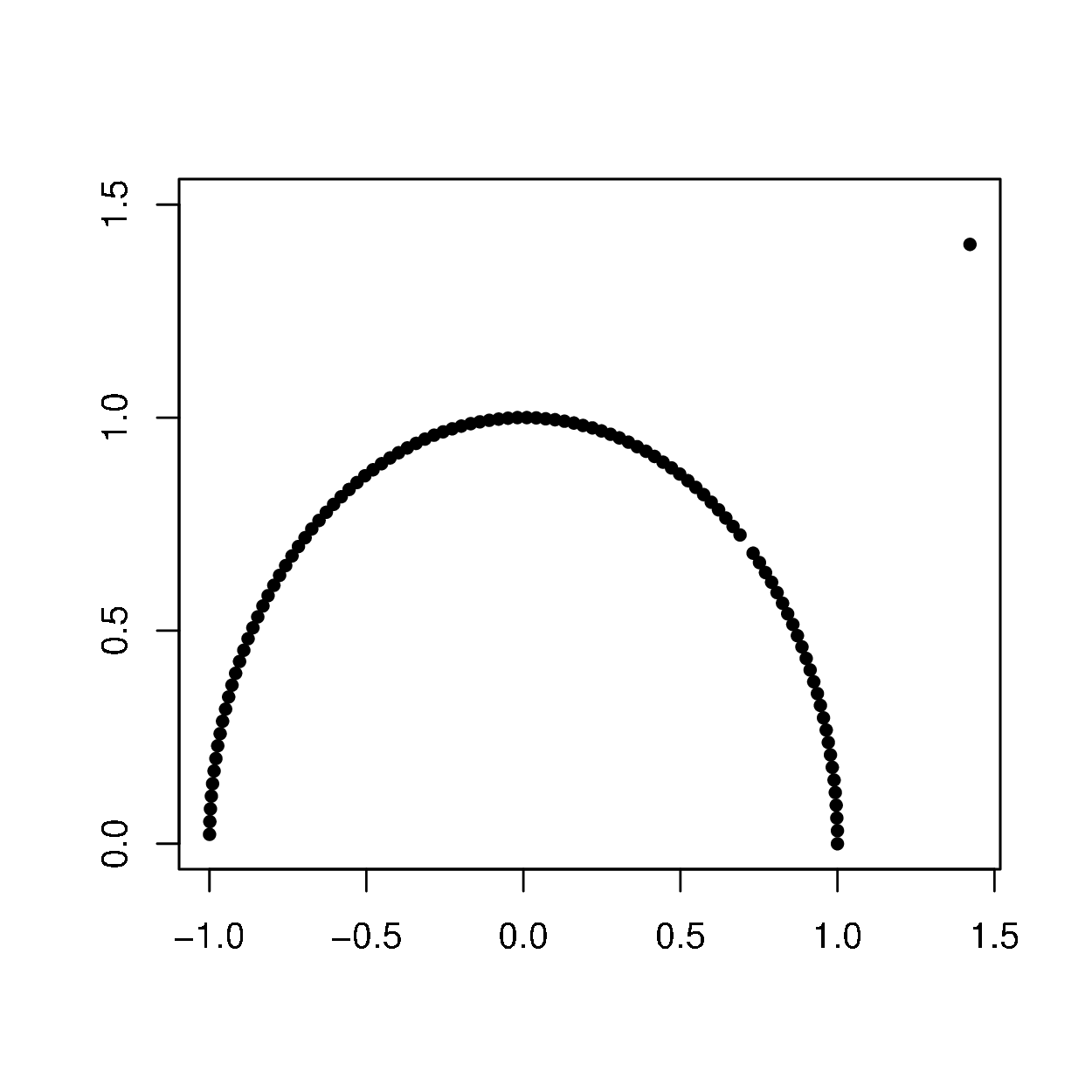}
  \vspace*{-12pt}
  \caption{Weights associated to each point
    \textup(see~\protect\eqref{eq:weight_points}\textup) of the lattice $\bar
    S_\Delta^1$ for $\Delta = 0.03$, $\vartheta = (1 / \sqrt{2}, 1 /
    \sqrt{2})$ and $T = 0.05, 0.2, 0.5, 10$ \textup(from top to bottom
    and left to right.\textup)}
  \label{fig:weightsd2}
\end{figure}

\begin{figure}[t]
  \includegraphics[width = \figlen]{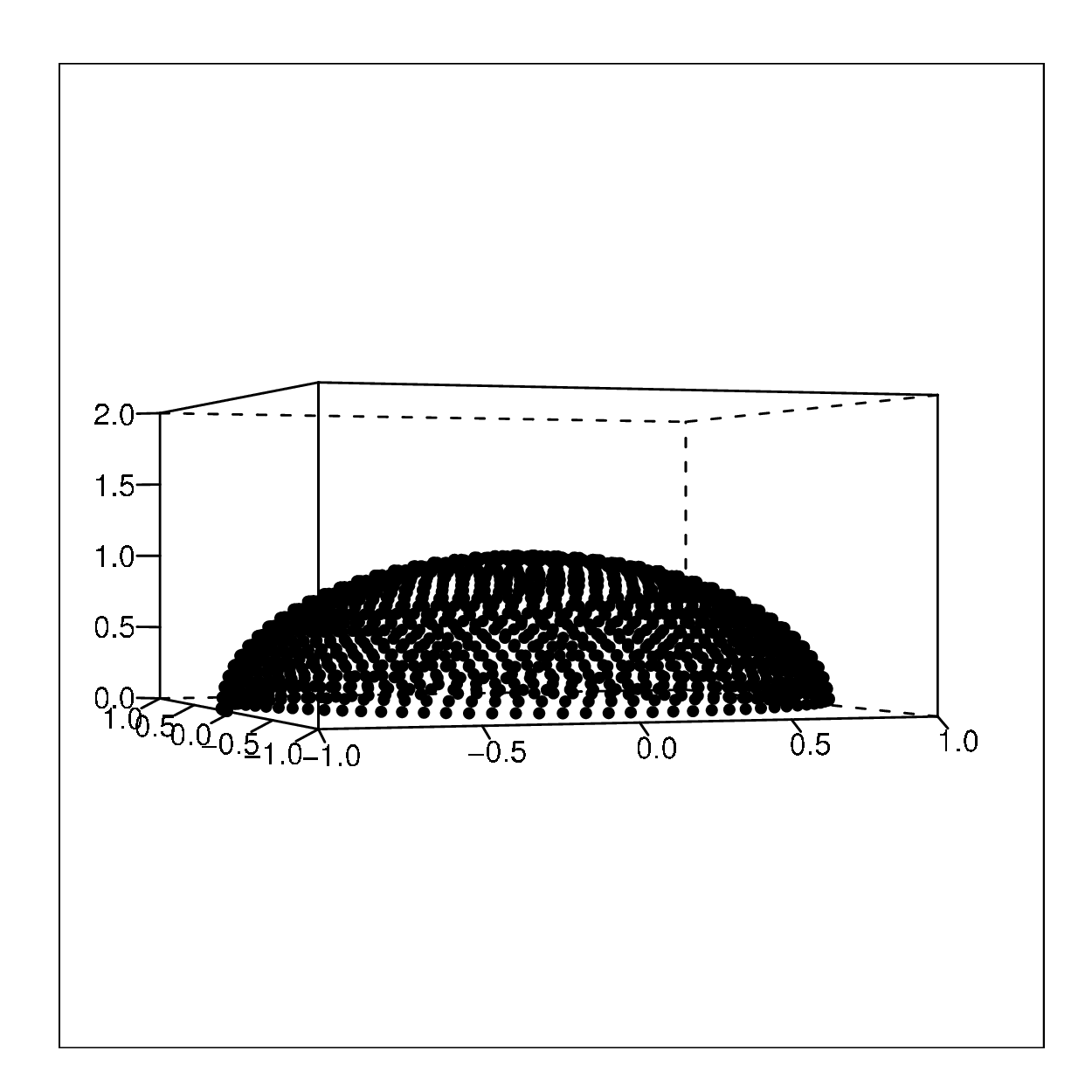}
  \includegraphics[width = \figlen]{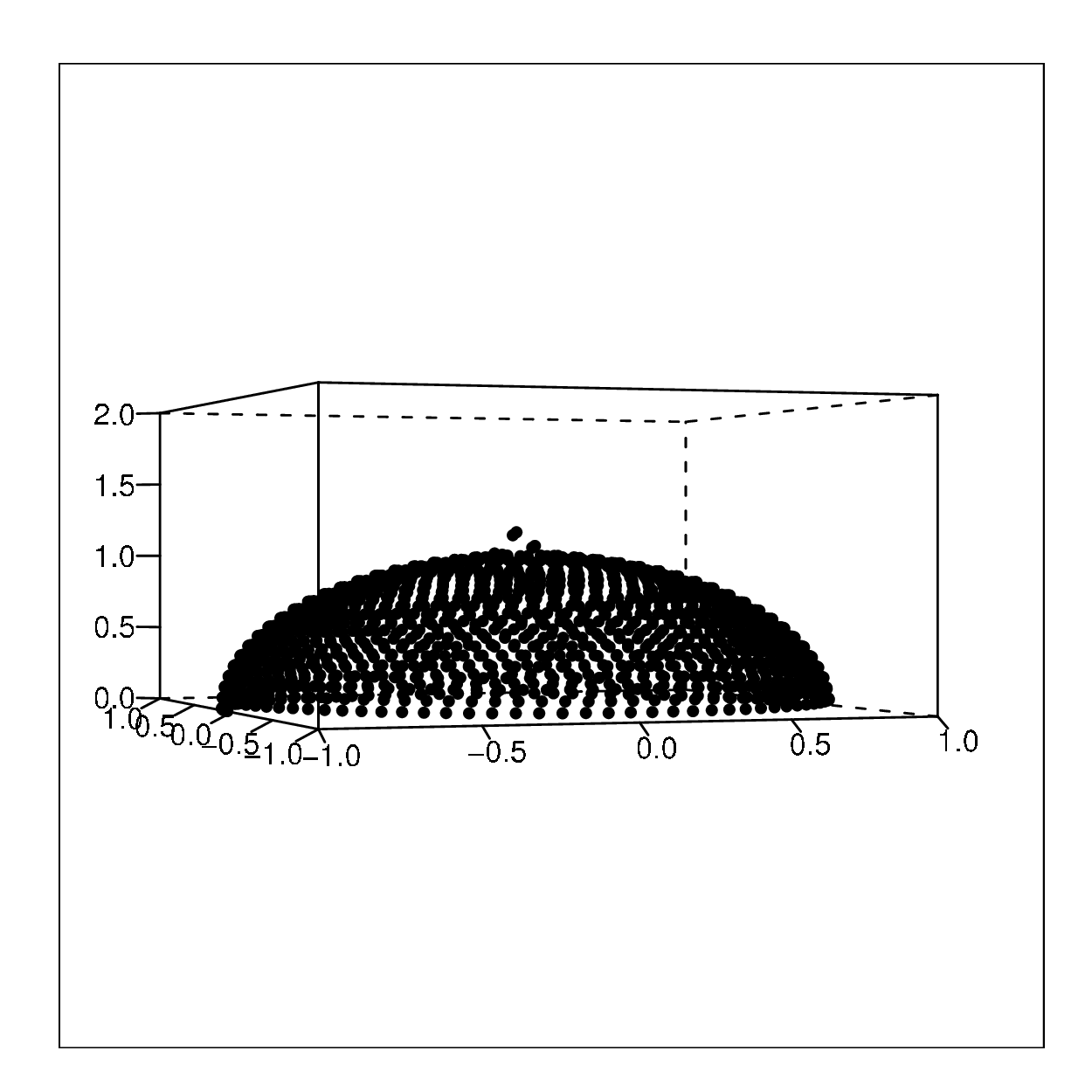} \\
  \includegraphics[width = \figlen]{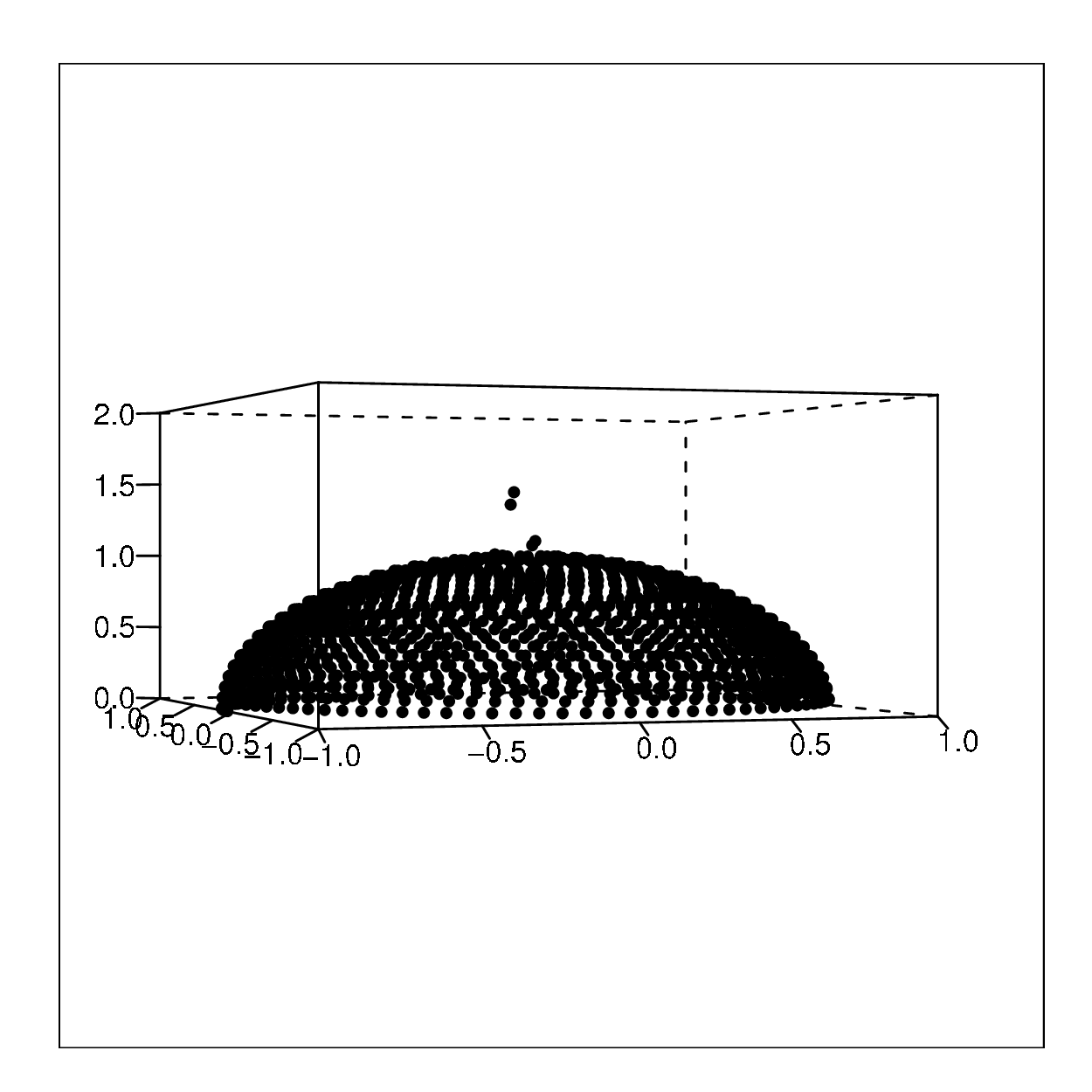}
  \includegraphics[width = \figlen]{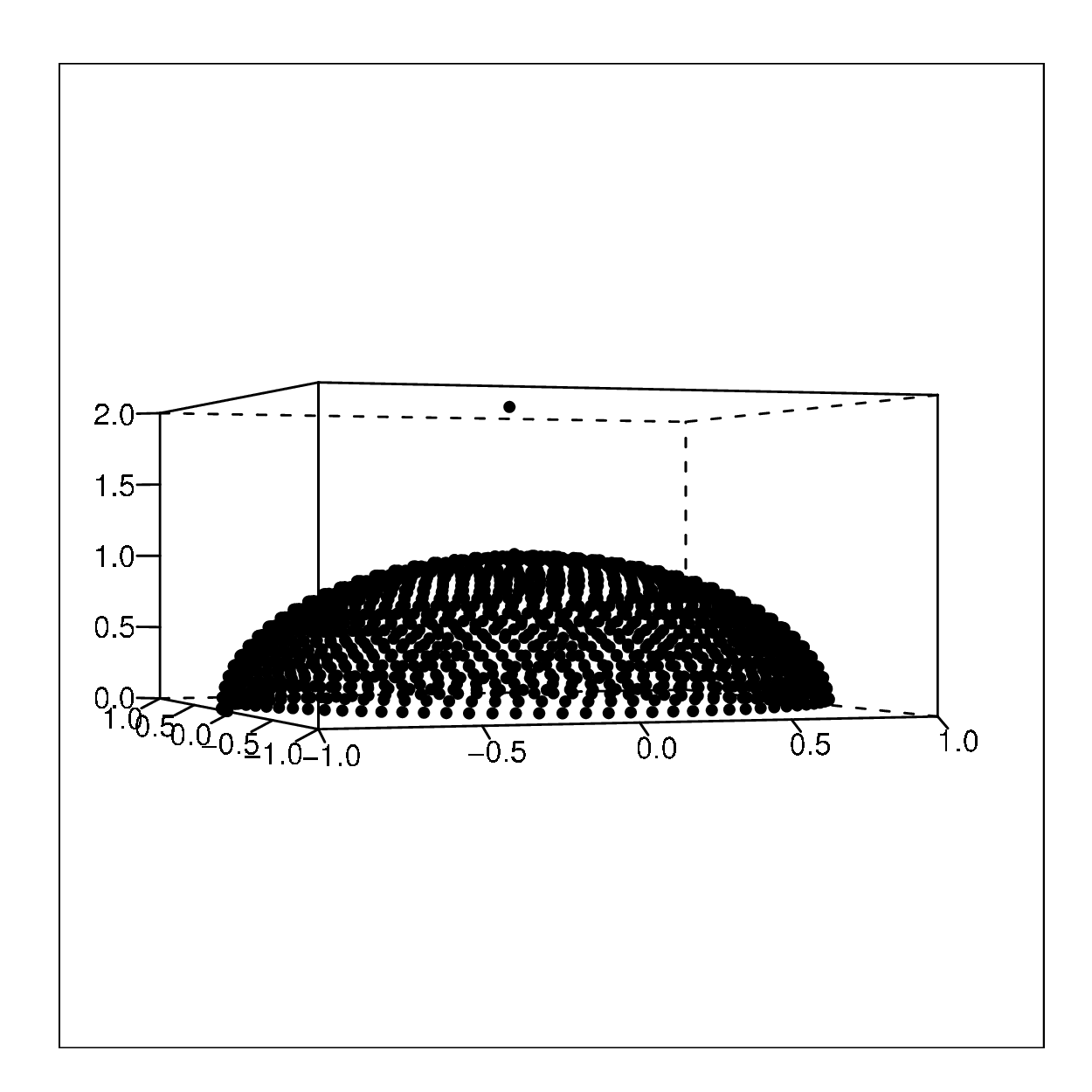}
  \caption{Weights associated to each
    points~\textup(see~\protect\eqref{eq:weight_points}\textup) of the lattice
    $\bar S_\Delta^2$ for $\Delta = 0.07$, $\vartheta = (0, 0, 1)$,
    and $T = 0.05, 0.3, 0.5, 10$ \textup(from top to bottom and left
    to right\textup).}
  \label{fig:weightsd3}
\end{figure}

\begin{figure}[htbp]
  \includegraphics[width = 6cm]{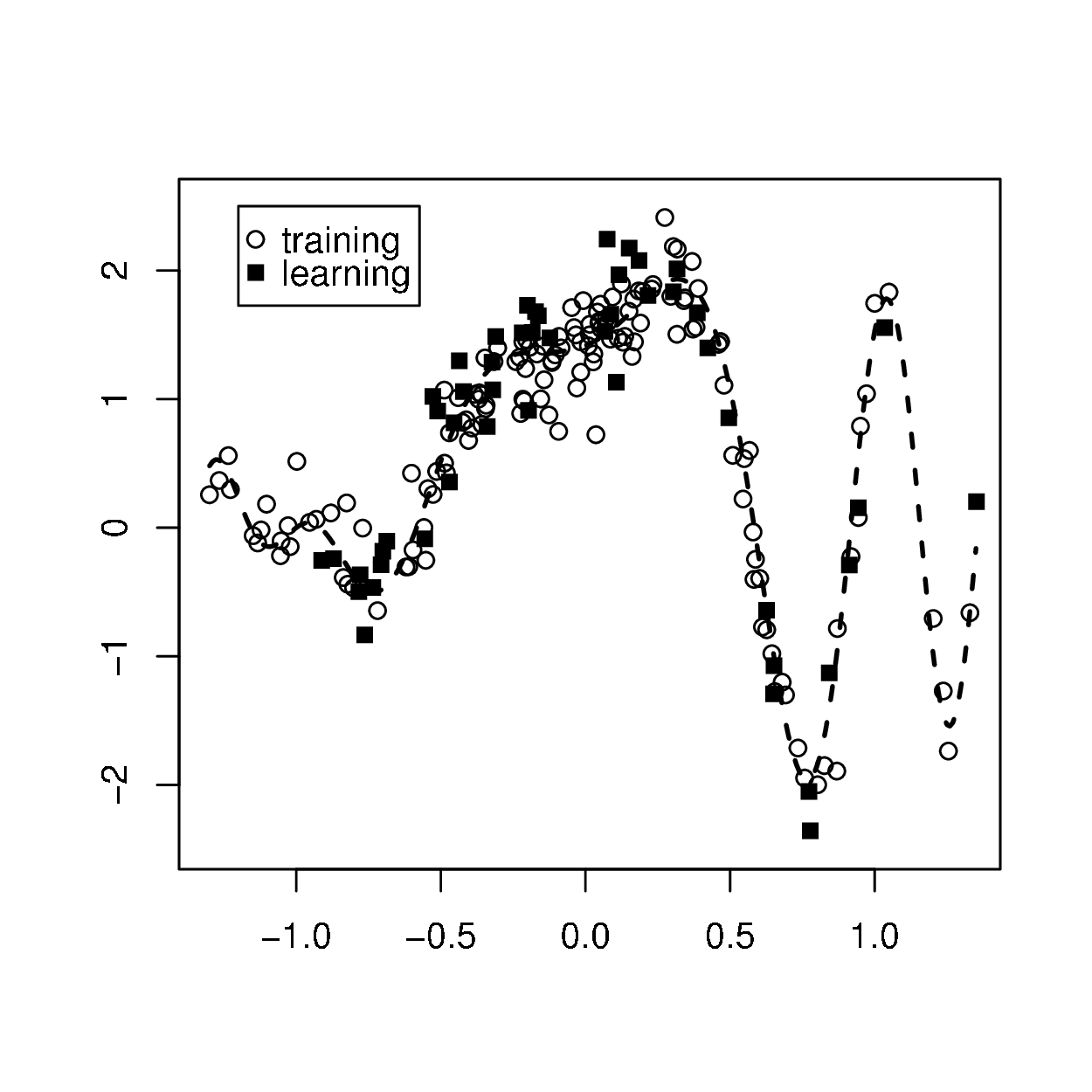}
  \includegraphics[width = 6cm]{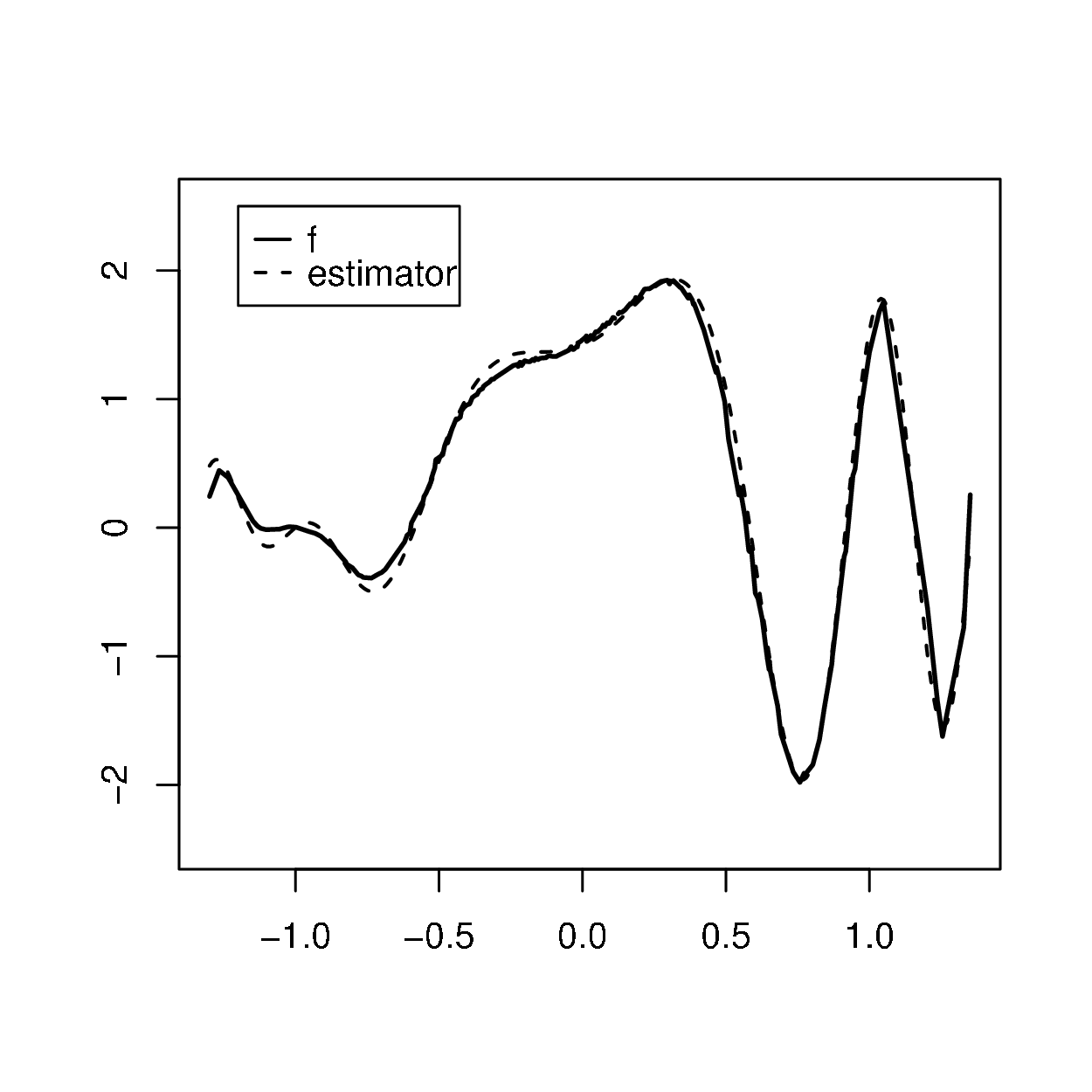}
  \includegraphics[width = 6cm]{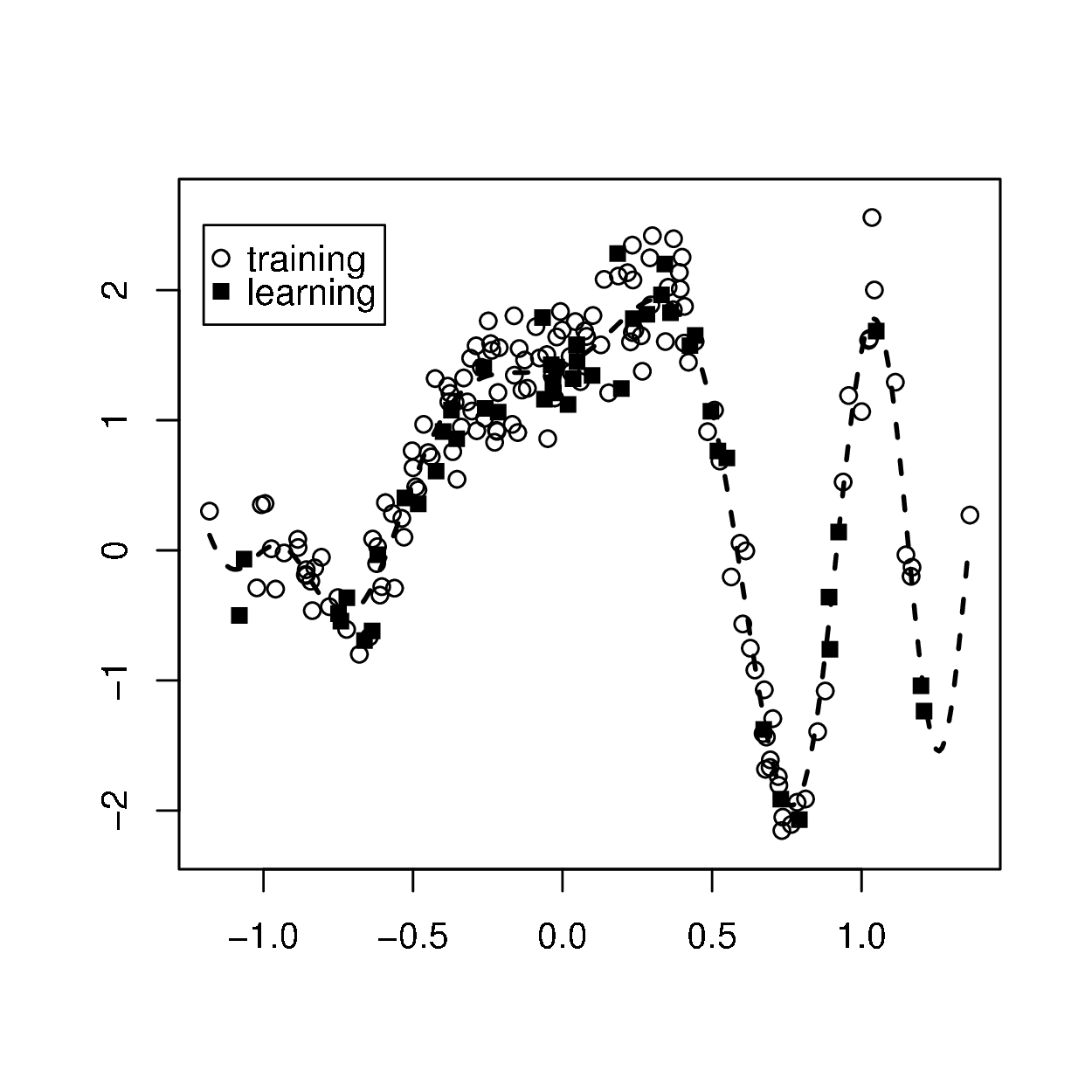}
  \includegraphics[width = 6cm]{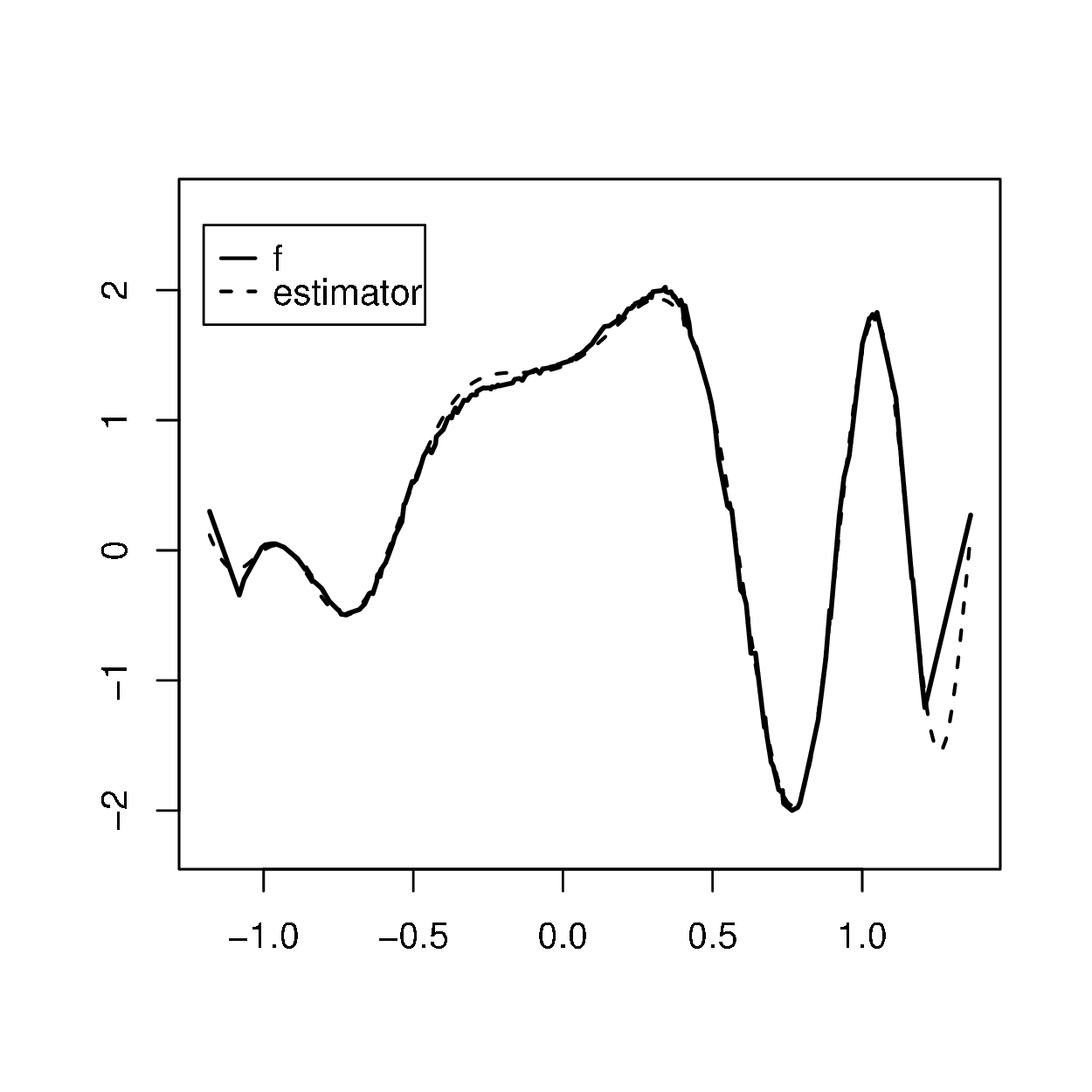}
  \includegraphics[width = 6cm]{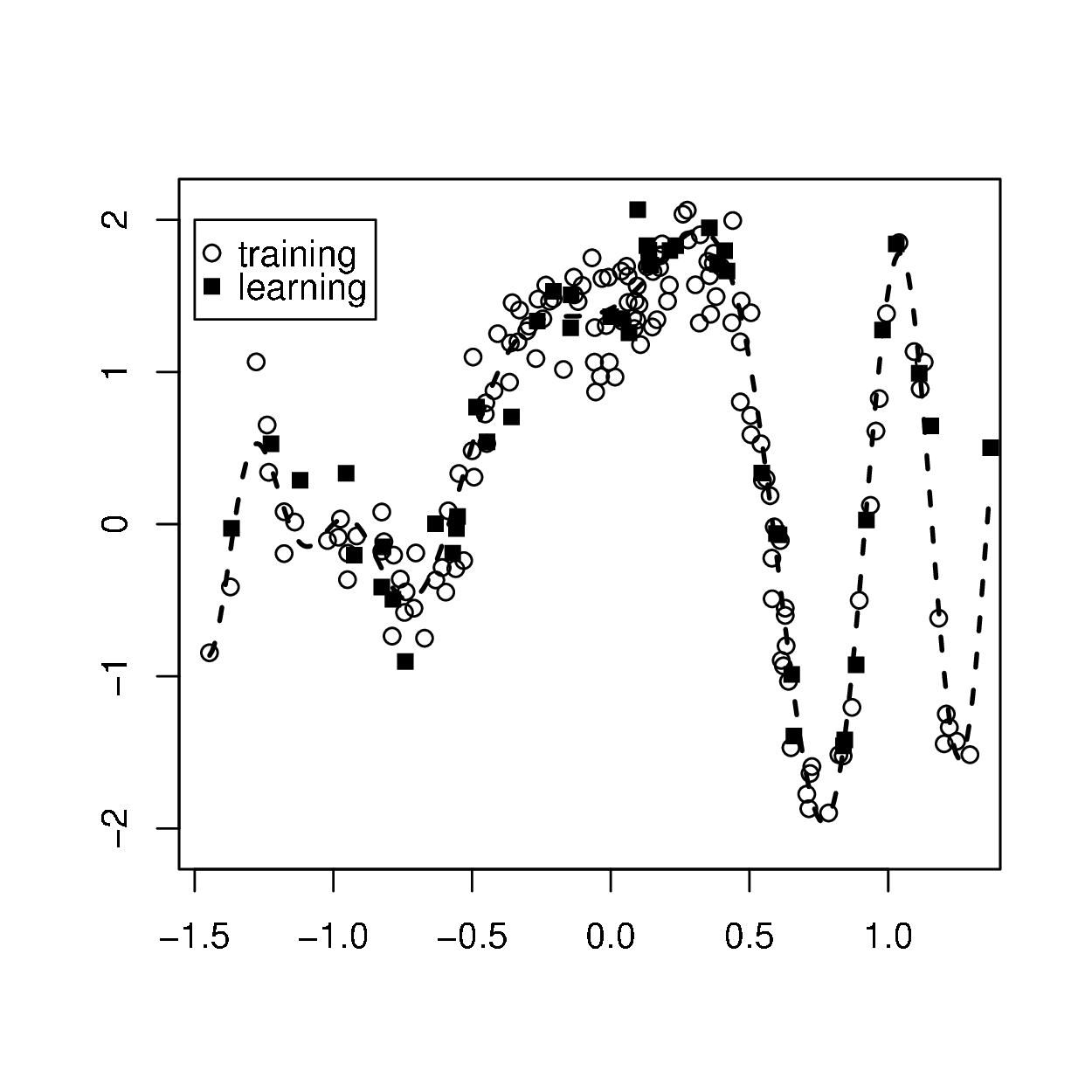}
  \includegraphics[width = 6cm]{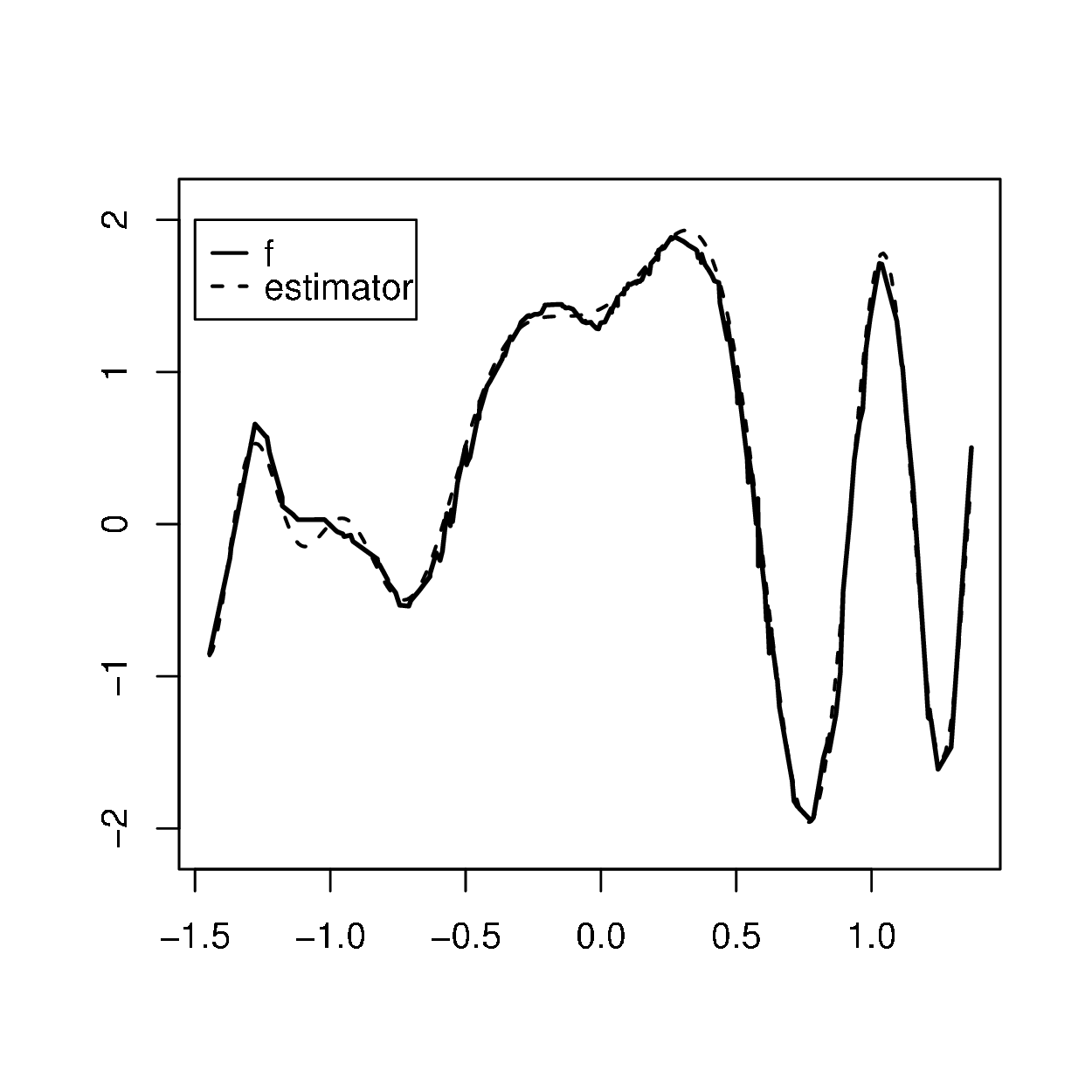}
  \caption{Simulated datasets and aggregated estimators with
    cross-validated temperature for $f = \mathtt{hardsine}$, $n=200$,
    and indexes $\vartheta = (1 / \sqrt{2}, 1 / \sqrt{2})$, $\vartheta
    = (2 / \sqrt{14}, 1/\sqrt{14}, 3/\sqrt{14})$,
    $\vartheta=(1/\sqrt{21}, -2/\sqrt{21}, 0, 4/\sqrt{21} )$ from top
    to bottom.}
  \label{fig:n200hardsine}
\end{figure}

\begin{figure}[htbp]
  \includegraphics[width = 6cm]{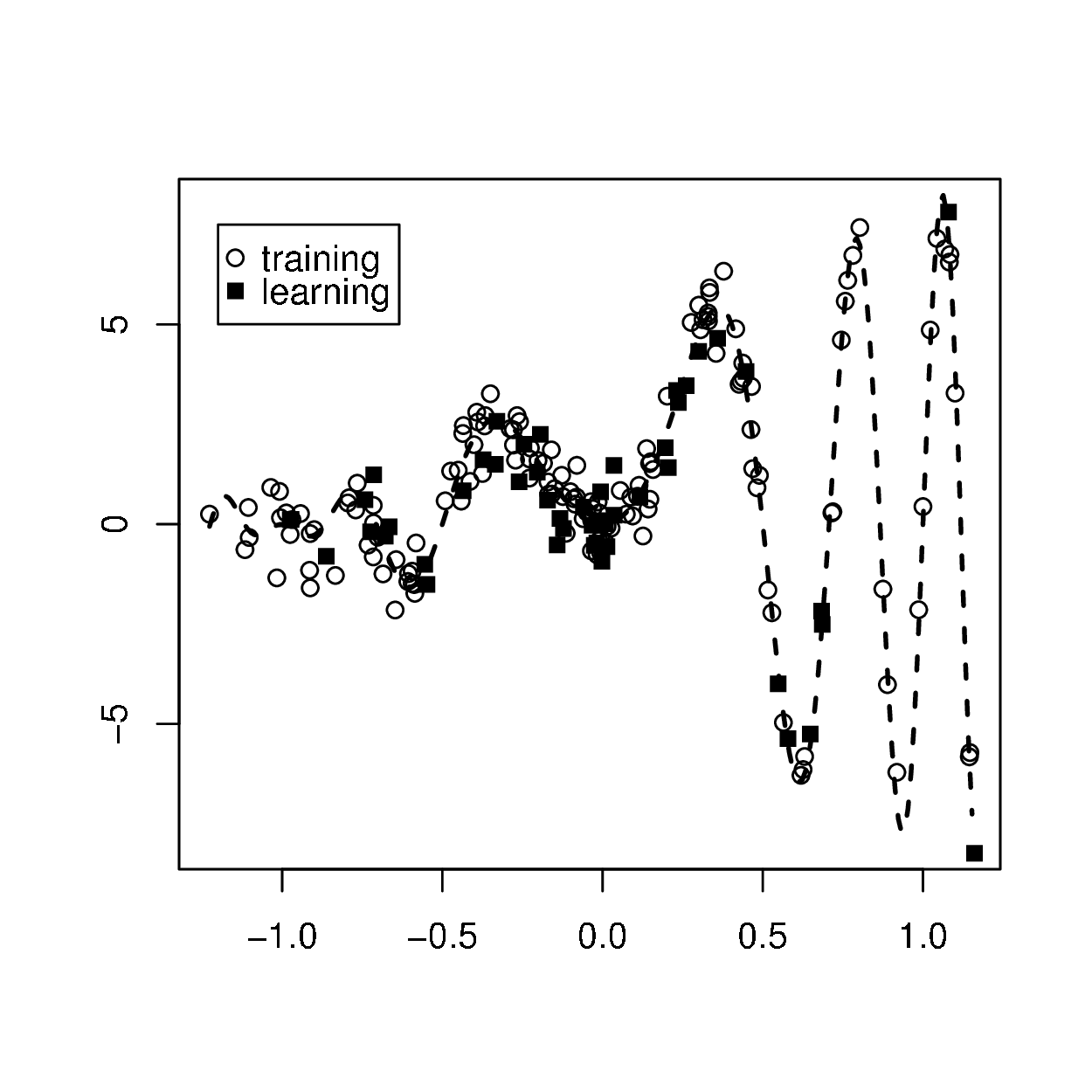}
  \includegraphics[width = 6cm]{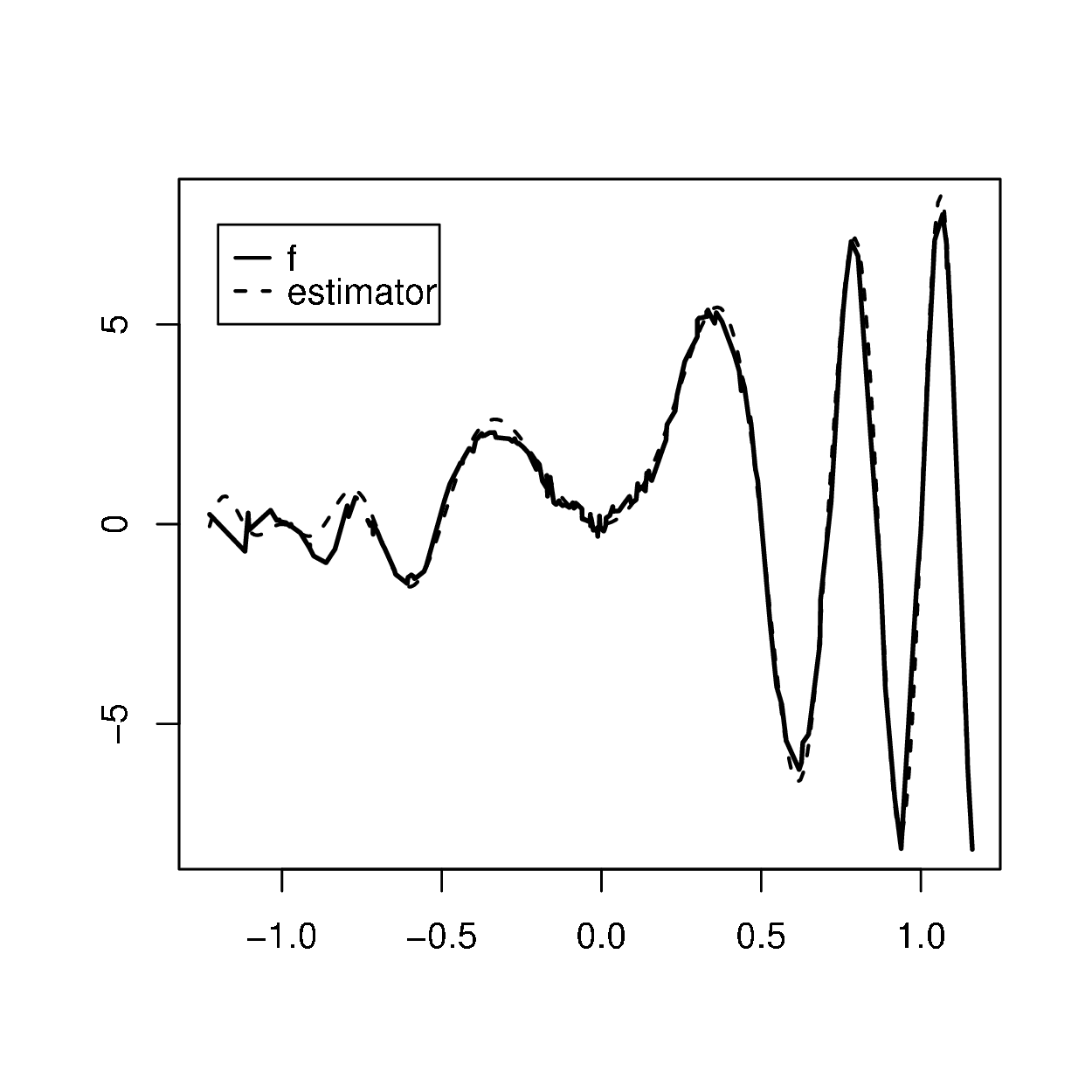}
  \includegraphics[width = 6cm]{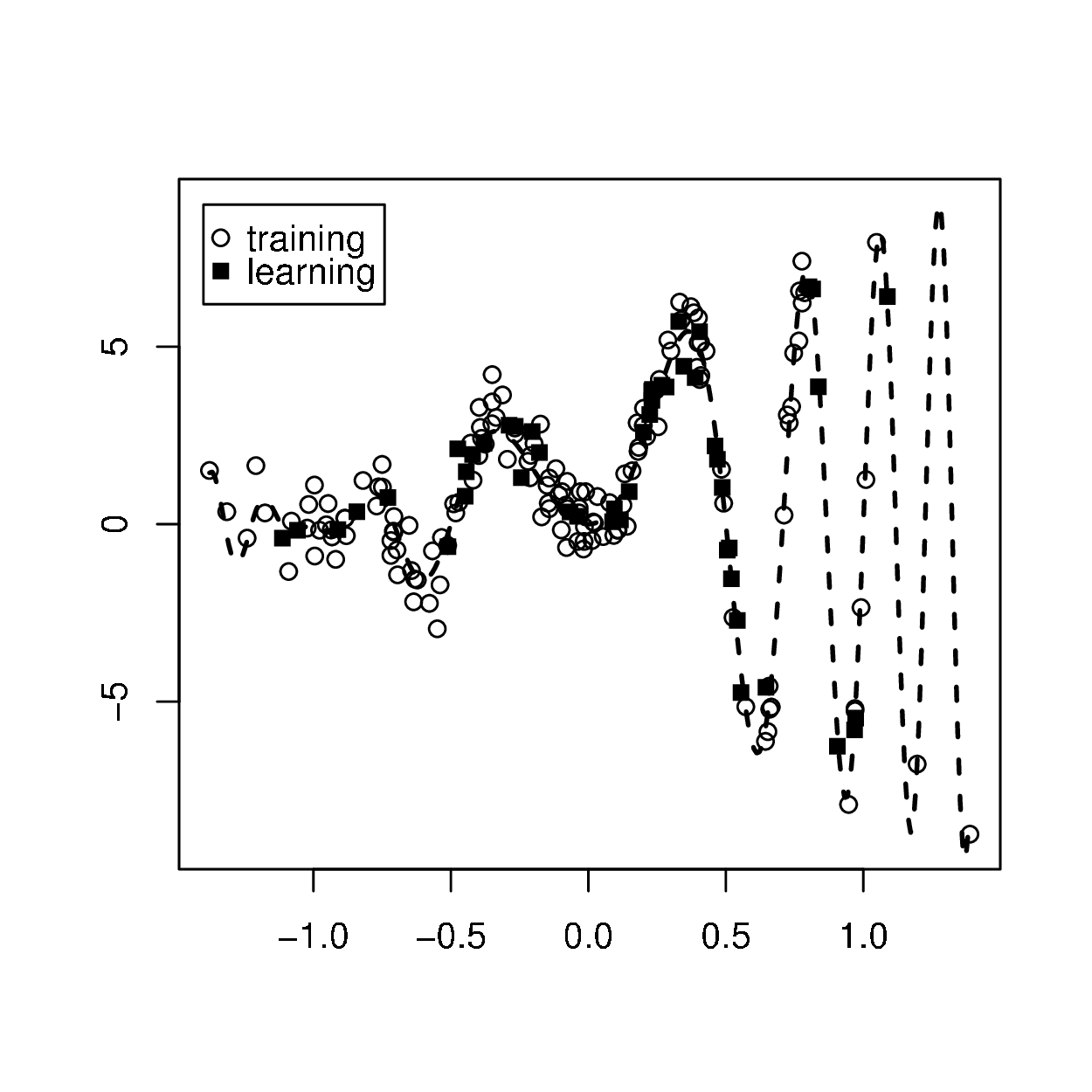}
  \includegraphics[width = 6cm]{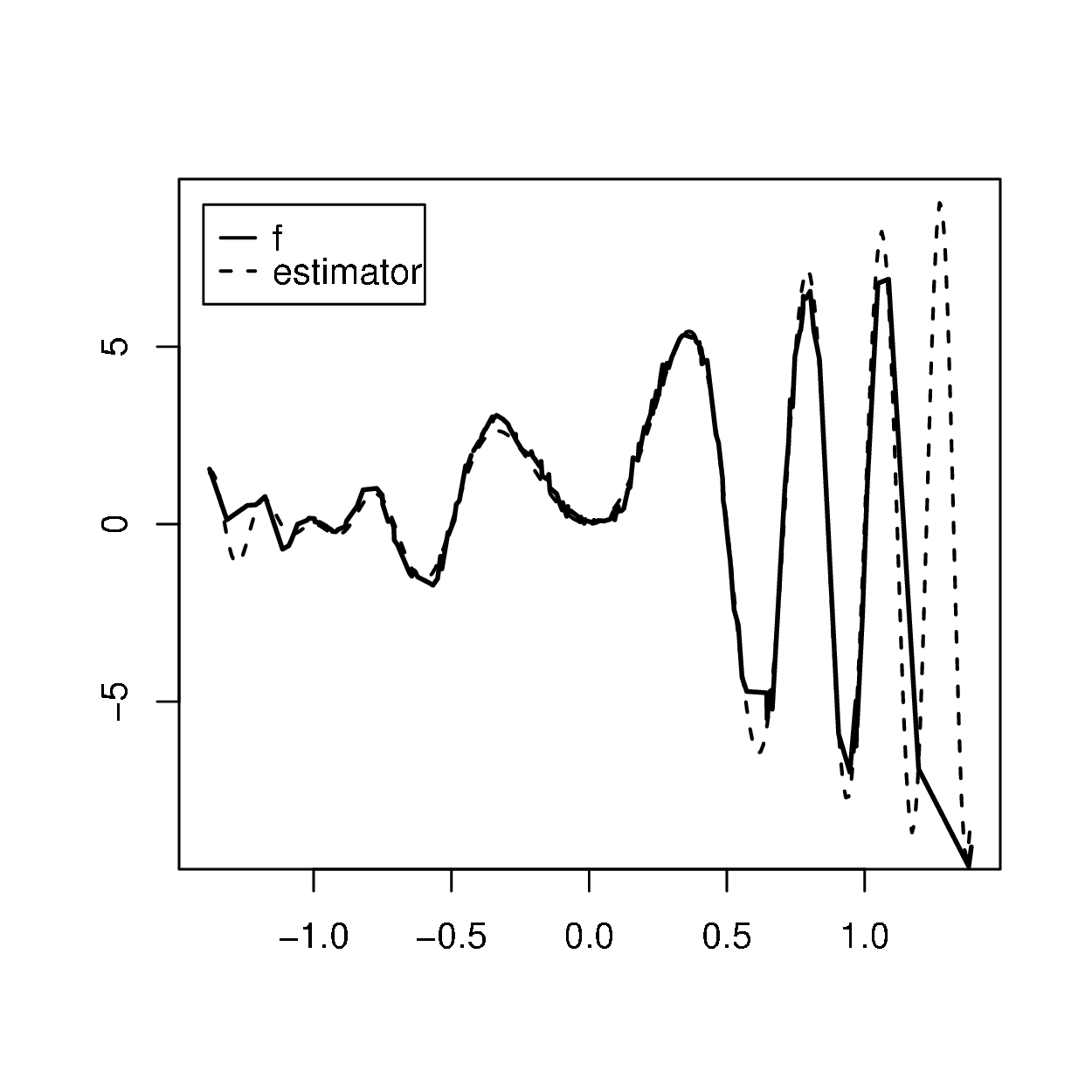}
  \includegraphics[width = 6cm]{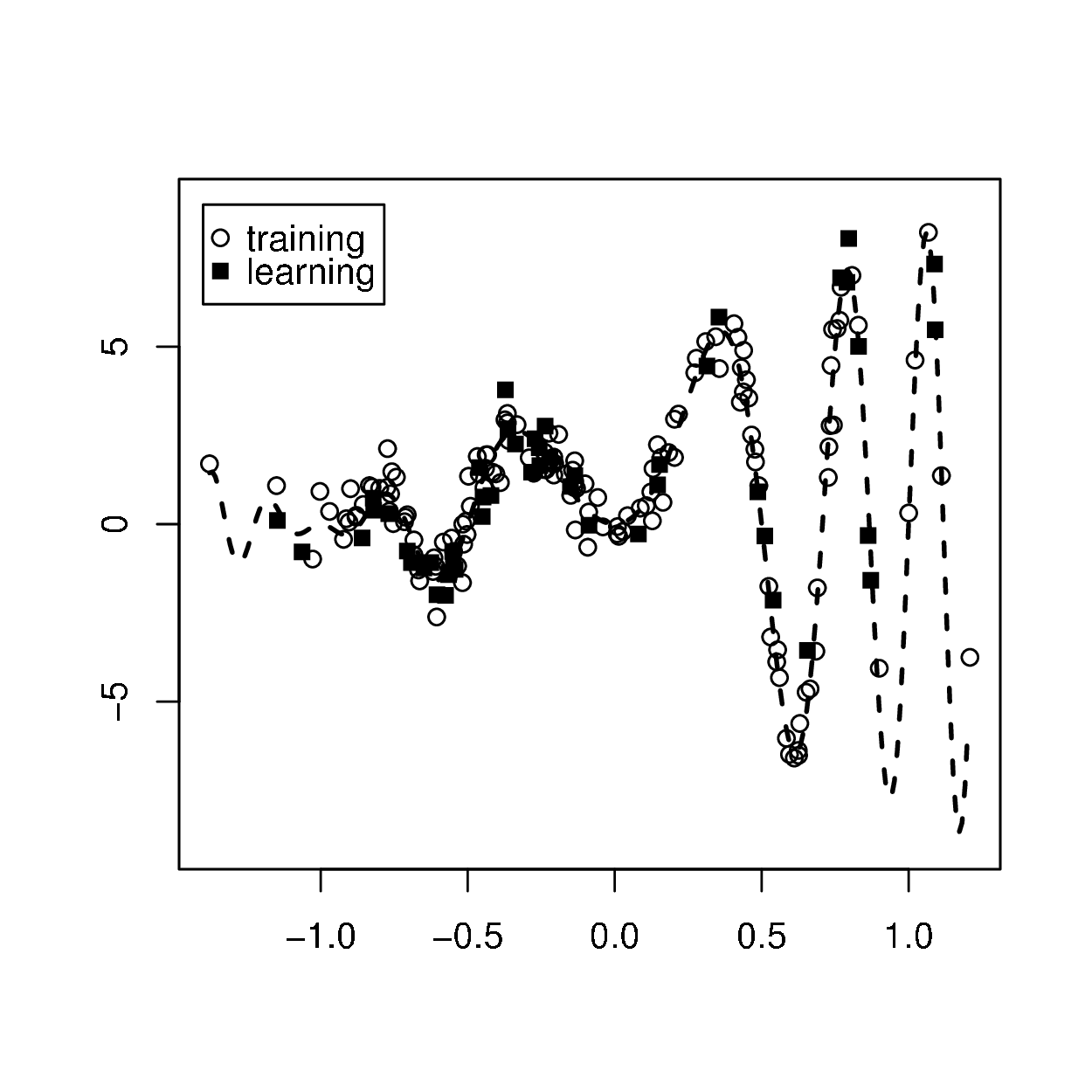}
  \includegraphics[width = 6cm]{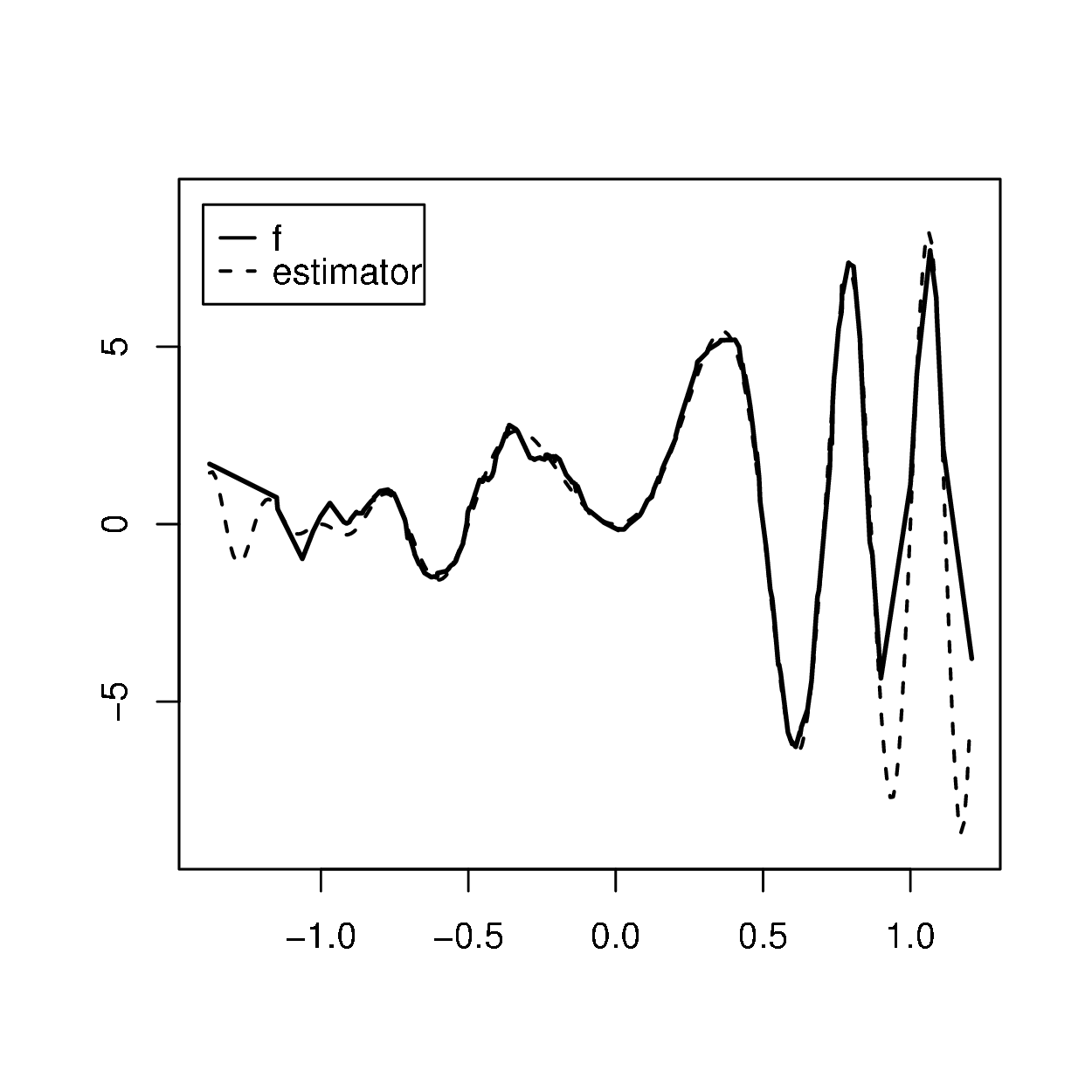}
  \caption{Simulated datasets and aggregated estimators with
    cross-validated temperature for $f = \mathtt{oscsine}$, $n=200$,
    and indexes $\vartheta = (1 / \sqrt{2}, 1 / \sqrt{2})$, $\vartheta
    = (2 / \sqrt{14}, 1/\sqrt{14}, 3/\sqrt{14})$,
    $\vartheta=(1/\sqrt{21}, -2/\sqrt{21}, 0, 4/\sqrt{21} )$ from top
    to bottom.}
  \label{fig:n200oscsine}
\end{figure}

The aim of Figures~\ref{fig:weightsd2} and~\ref{fig:weightsd3} is to
give an illustration of the aggregation phenomenon. In these figures,
we show the points
\begin{equation}
  \label{eq:weight_points}
  \big\{ (1 + w(\bar g^{(\lambda)})) \vartheta \text{ for } \lambda =
  (\vartheta, s, L) \in \Lambda = \bar S_\Delta^{d-1} \times \{ 3 \}
  \times \{ 1 \} \big\}
\end{equation}
obtained for a single run (that is, we take $s = 3$ and $L = 1$ in the
bandwidth choice~\eqref{eq:bandwidth} and we do not use the reduction
of complexity algorithm). These figures motivates the use of the
complexity reduction algorithm, since only the weights corresponding
to a point of $\bar S_\Delta^{d-1}$ which is close to the true index
are significant (at least numerically). Moreover, these weights
provide information about the true index: the direction $v \in \bar
S_\Delta^{d-1}$ corresponding to the largest coefficient $w(\bar
g^{(\lambda)})$ for $\lambda = (v, s, L)$ is an accurate estimator of
the index, see Figures~\ref{fig:weightsd2}
and~\ref{fig:weightsd3}. Finally, we show typical realisations for
several index functions, indexes and sample sizes in
Figures~\ref{fig:n100hardsine}, \ref{fig:n100oscsine},
\ref{fig:n200hardsine}, \ref{fig:n200oscsine}.

\section{Proofs}

\subsection*{Proof of Theorem~\ref{thm:UB}}

The functions $\bar g^{(\lambda)}$ are given by~\eqref{eq:weak}. They
are computed based on the training (or ``frozen'') sample $D_m$, which
is independent of the learning sample $D_{(m)}$. If $E^{(m)}$ denotes
the integration with respect to the joint law of $D_{(m)}$, we obtain
using Theorem~\ref{thm:oracle}:
\begin{align*}
  E^{(m)} \norm{\hat g - g}_{L^2(P_X)}^2 &\leq (1 + a) \min_{\lambda
    \in \Lambda} \norm{\bar g^{(\lambda)} - g}_{L^2(P_X)}^2 + \frac{C
    \log |\Lambda| (\log |D_{(m)}|)^{1/2}}{|D_{(m)}|} \\
  &\leq (1 + a) \norm{\bar g^{(\bar \lambda)} - g}_{L^2(P_X)}^2 +
  o(n^{-2s / (2s + 1)}),
\end{align*}
since $\log |\Lambda| (\log |D_{(m)}|)^{1/2} / |D_{(m)}| \leq d (\log
n)^{3/2 + \gamma} / (2 \smin n)$ (see~\eqref{eq:learning_size} and
\eqref{eq:lattice_step}), and where $\bar \lambda = (\bar \vartheta,
\bar s) \in \Lambda$ is such that $\norm{\bar \vartheta - \vartheta}_2
\leq \Delta$ and $\lfloor \bar s \rfloor = \lfloor s \rfloor$ with $s
\in [\bar s, \bar s + (\log n)^{-1}]$. By integration with respect to
$P^m$, we obtain
\begin{equation}
  \label{eq:step1}
  E^n \norm{\hat g - g}_{L^2(P_X)}^2 \leq (1 + a) E^m \norm{\bar
    g^{(\bar \lambda)} - g}_{L^2(P_X)}^2 + o(n^{-2s / (2s + 1)}).
\end{equation}
The choice of $\bar \lambda$ entails $H^Q(s, L) \subset H^Q(\bar s,
L)$ and
\begin{equation*}
  n^{-2 \bar s / (2 \bar s + 1)} \leq e^{1/2} n^{-2s / (2s + 1)}.
\end{equation*}
Thus, together with~\eqref{eq:learning_size} and~\eqref{eq:step1}, the
Theorem follows if we prove that
\begin{equation}
  \label{eq:step2}
  \sup_{f \in H^Q(\bar s, L)} E^m \norm{\bar g^{(\bar \lambda)} -
    g}_{L^2(P_X)}^2 \leq C m^{-2 \bar s / (2 \bar s + 1)}.
\end{equation}
for $n$ large enough, where $C > 0$. We cannot use directly
Theorem~\ref{thm:LPE} to prove this, since the weak estimator $\bar
g^{(\bar \lambda)}$ works based on data $D_m(\bar \vartheta)$
(see~\eqref{eq:proj_sample}) while the true index is $\vartheta$. In
order to clarify the proof, we write $\bar g^{(\bar \vartheta)}$
instead of $\bar g^{(\bar \lambda)}$ since in~\eqref{eq:step2}, the
estimator uses the ``correct'' smoothness parameter $\bar s$. We have
\begin{equation*}
  \| \bar g^{(\bar \vartheta)} - g \|_{L^2(P_X)}^2 \leq 2 \big( \norm{
    \bar g^{(\bar \vartheta)}(\cdot) - f(\bar \vartheta\T \cdot)
  }_{L^2(P_X)}^2 + \norm{ f(\bar \vartheta\T \cdot) - f(\vartheta\T
    \cdot)}_{L^2(P_X)}^2 \big)
\end{equation*}
and using together~\eqref{eq:lattice_step} and $f \in H^Q(s, L)$ for
$s \geq \smin$, we obtain
\begin{equation*}
  \norm{ f(\bar \vartheta\T \cdot) - f(\vartheta\T
    \cdot)}_{L^2(P_X)}^2 \leq L^2 \int \norm{x}_2^{2 \smin}
  P_X(dx) \Delta^{2 \smin} \leq C (n \log n)^{-1}.
\end{equation*}
Let us denote by $Q_{\vartheta}(\cdot | X_1^m)$ the joint law of
$(X_i, Y_i)_{1 \leq i \leq m}$ from model~\eqref{eq:model} (when the
index is $\vartheta$) conditional on the $(X_i)_{1 \leq i \leq m}$,
which is given by
\begin{equation*}
  Q_{\vartheta}(dy_1^m | x_1^m) := \prod_{i = 1}^m
  \frac{1}{(\sigma(x_i) (2 \pi)^{1/2})} \exp\Big( -\frac{ (y_i -
    f(\vartheta\T x_i))^2}{2 \sigma(x_i)^2} \Big) d y_i.
\end{equation*}
Under $Q_{\bar\vartheta}(\cdot | X_1^m)$, we have
\begin{align*}
  &L_X(\vartheta, \bar \vartheta) := \frac{d Q_{\vartheta}( \cdot |
    X_1^m)}{d Q_{\bar \vartheta}(\cdot | X_1^m)} \\
  &\overset{(\text{law})}{=} \exp\Big( -\sumim \frac{\epsilon_i
    (f(\bar \vartheta\T X_i) - f(\vartheta\T X_i) ) }{ \sigma(X_i)} -
  \frac{1}{2} \sumim \frac{(f(\bar \vartheta\T X_i) - f(\vartheta\T
    X_i))^2}{\sigma(X_i)^2} \Big).
\end{align*}
Hence, if $P_X^m$ denotes the joint law of $(X_1, \ldots, X_m$),
\begin{align}
  \label{eq:step3}
  \nonumber E^m \norm{ \bar g^{(\bar \vartheta)}(&\cdot) - f(\bar
    \vartheta\T \cdot) }_{L^2(P_X)}^2 \\
  \nonumber &= \int \int \norm{ \bar g^{(\bar \vartheta)}(\cdot) -
    f(\bar \vartheta\T \cdot) }_{L^2(P_X)}^2 L_X(\vartheta,
  \bar\vartheta) d Q_{\bar\vartheta}(\cdot | X_1^m) d P_X^m \\
  &\leq C \int \int \norm{ \bar f^{(\bar \vartheta)}(\bar \vartheta\T
    \cdot) - f(\bar \vartheta\T \cdot) }_{L^2(P_X)}^2 d
  Q_{\bar\vartheta}(\cdot | X_1^m) d P_X^m \\
  \nonumber &+ 4 Q^2 \int \int L_X(\vartheta, \bar\vartheta)
  \ind{\{L_X(\vartheta, \bar\vartheta) \geq C \}} d
  Q_{\bar\vartheta}(\cdot | X_1^m) d P_X^m,
\end{align}
where we decomposed the integrand over $\{ L_X(\vartheta,
\bar\vartheta) \geq C \}$ and $\{ L_X(\vartheta, \bar\vartheta) \leq C
\}$ for some constant $C \geq 3$, and where we used the fact that
$\norm{\bar g^{(\bar \vartheta)}}_{\infty} ,\norm{f}_{\infty} \leq
Q$. Under $Q_{\bar\vartheta}(\cdot | X_1^m)$, the $(X_i,Y_i)$ have the
same law as $(X, Y)$ from model~\eqref{eq:model} where the index is
$\bar\vartheta$. Moreover, we assumed that $P_{\bar \vartheta\T X}$
satisfies Assumption~\ref{assd}.  Hence, Theorem~\ref{thm:LPE} entails that,
uniformly for $f \in H^Q(\bar s, L)$,
\begin{equation*}
  \int \int \norm{ \bar f^{(\bar \vartheta)}(\bar \vartheta\T
    \cdot) - f(\bar \vartheta\T \cdot) }_{L^2(P_X)}^2 d
  Q_{\bar\vartheta}(\cdot | X_1^m) d P_X^m \leq C' m^{-2\bar s / (2
    \bar s + 1)}.
\end{equation*}
Moreover, the second term in the right hand side of~\eqref{eq:step3}
is smaller than
\begin{align*}
  4 Q^2 \int \Big( &\int L_X(\vartheta, \bar\vartheta)^2 d
  Q_{\bar\vartheta}(\cdot | X_1^m) \Big)^{1/2} Q_{\bar\vartheta} \big[
  L_X(\vartheta, \bar\vartheta) \geq C | X_1^m \big]^{1/2} dP_X^m.
\end{align*}
Since $f \in H^Q(s, L)$ for $s \geq \smin$, since $P_X$ is compactly
supported and since $\sigma(X) > \sigma_0$ a.s., we obtain
using~\eqref{eq:lattice_step}:
\begin{align*}
  \int L_X(\vartheta, \bar\vartheta)^2 d Q_{\bar\vartheta}(\cdot |
  X_1^m) \leq \exp\Big( \frac{1}{2} \sumim \frac{(f(\bar \vartheta\T
    X_i) - f(\vartheta\T X_i))^2}{\sigma(X_i)^2} \Big) \leq 1
\end{align*}
$P_X^m$-a.s. when $m$ is large enough. Moreover, with the same
arguments we have
\begin{align*}
  Q_{\bar\vartheta} \big[L_X(\vartheta, \bar\vartheta) \geq C | X_1^m
  \big] \leq m^{-(\log C)^2 / 2} \leq m^{-4 \bar s / (2 \bar s + 1) }
\end{align*}
for $C$ large enough, where we use the standard Gaussian deviation $P[
N(0, b^2) \geq a ] \leq \exp( -a^2 / (2 b^2) )$. This concludes the
proof of Theorem~\ref{thm:UB}. \hfill $\square$

\subsection*{Proof of Theorem~\ref{thm:LB}}

We want to bound the minimax risk
\begin{equation}
  \label{eq:LBproof1}
  \inf_{\tilde g} \sup_{f \in H^Q(s, L)} E^n \int \big( \tilde g(x) -
  f( \vartheta\T x) \big)^2 P_X( dx )
\end{equation}
from below, where the infimum is taken among all estimators $\mathbb
R^d \rightarrow \mathbb R$ based on data from
model~\eqref{eq:model},\eqref{eq:SI}. We recall that $\vartheta\T X$
satisfies Assumption~\ref{assd}. We consider $\vartheta^{(2)}, \ldots,
\vartheta^{(d)}$ in $\mathbb R^d$ such that $(\vartheta,
\vartheta^{(2)}, \ldots, \vartheta^{(d)})$ is an orthogonal basis of
$\mathbb R^d$. We denote by $\mathbf O$ the matrix with columns
$\vartheta, \vartheta^{(2)}, \ldots, \vartheta^{(d)}$. We define $Y :=
\mathbf O X = (Y^{(1)}, \ldots, Y^{(d)})$ and $Y_2^d := ( Y^{(2)},
\ldots, Y^{(d)} )$. By the change of variable $y = \mathbf O x$, we
obtain
\begin{align*}
  \int_{\mathbb R^d} \big( &\tilde g(x) - f( \vartheta\T x)
  \big)^2 P_X( dx ) \\
  &= \int_{\mathbb R^d} \big( \tilde g( \mathbf O^{-1} y) - f( y^{(1)}
  ) \big)^2 P_Y( dy ) \\
  &= \int_{\mathbb R} \int_{\mathbb R^{d-1}} \big( \tilde g( \mathbf
  O^{-1} y) - f( y^{(1)} ) \big)^2 P_{ Y_2^d | Y^{(1)} } ( d y_2^d |
  y^{(1)} ) P_{Y^{(1)}} (d y^{(1)} ) \\
  &\geq \int_{\mathbb R} \big( \tilde f( y^{(1)}) - f( y^{(1)} )
  \big)^2 P_{ \vartheta\T X} ( d y^{(1)} ),
\end{align*}
where $\tilde f( y^{(1)} ) := \int \tilde g( \mathbf O^{-1} y) P_{
  Y_2^d | Y^{(1)} } ( d y_2^d | y^{(1)} )$. Hence, if $Z :=
\vartheta\T X$,~\eqref{eq:LBproof1} is larger than
\begin{equation}
  \label{eq:LBproof2}
  \inf_{\tilde f} \sup_{f \in H^Q(s, L)} E^n \int \big( \tilde f(z) -
  f(z) \big)^2 P_{Z}( dz ),
\end{equation}
where the infimum is taken among all estimators $\mathbb R \rightarrow
\mathbb R$ based on data from model~\eqref{eq:model} with $d=1$
(univariate regression). In order to bound~\eqref{eq:LBproof2} from
below, we use the following Theorem, from \cite{tsybakov03}, which is
a standard tool for the proof of such a lower bound. We say that
$\partial$ is a \emph{semi-distance} on some set $\Theta$ if it is
symmetric, if it satisfies the triangle inequality and if
$\partial(\theta, \theta) = 0$ for any $\theta \in \Theta$. We
consider $K(P | Q) := \int \log (\frac{dP}{dQ}) dP$ the
Kullback-Leibler divergence between probability measures $P$ and $Q$.

\begin{theorem}
  \label{thm:tsybokov}
  Let $(\Theta, \partial)$ be a set endowed with a semi-distance
  $\partial$. We suppose that $\{ P_\theta ; \theta \in \Theta \}$ is
  a family of probability measures on a measurable space $(\cX, \cA)$
  and that $(v_n)_{n \in \mathbb{N} }$ is a sequence of positive
  numbers. If there exist $\{ \theta_0, \ldots, \theta_M \} \subset
  \Theta$, with $M \geq 2$, such that
  \begin{itemize}
  \item $\partial(\theta_j, \theta_k) \geq 2 v_n \quad \forall 0 \leq
    j < k \leq M$
  \item $P_{\theta_j} \ll P_{\theta_0} \quad \forall 1 \leq j \leq M,$
  \item $\frac{1}{M} \sum_{j=1}^M K(P_{\theta_j}^{n} |
    P_{\theta_0}^{n})\leq \alpha \log M$ for some $\alpha \in (0,
    1/8)$,
  \end{itemize}
  then
  \begin{equation*}
    \inf_{ \tilde \theta_n }\sup_{\theta \in \Theta}E_\theta^n
    [  ( v_n^{-1} \partial( \tilde \theta_n, \theta) )^2 ] \geq
    \frac{\sqrt{M}}{1 + \sqrt{M}} \bigg( 1 - 2 \alpha - 2
    \sqrt{ \frac{\alpha}{\log M}} \bigg),
  \end{equation*}
  where the infimum is taken among all estimators based on a sample of
  size~$n$.
\end{theorem}

Let us define $m := \lfloor c_0 n^{1 / (2s + 1)} \rfloor $, the
largest integer smaller than $c_0 n^{1 / (2s + 1)}$, where $c_0 >
0$. Let $\varphi : \mathbb R \rightarrow [0, +\infty)$ be a function
in $H^Q(s, 1/2; \mathbb R)$ with support in $[-1/2, 1/2]$. We take
$h_n := m^{-1}$ and $z_k := (k - 1/2) / m$ for $k \in \{ 1, \ldots, m
\}$. For $\omega \in \Omega := \{ 0, 1 \}^m$, we consider the
functions
\begin{equation*}
  f(\cdot; \omega) := \sum_{k=1}^m \omega_k \varphi_k(\cdot) \; \text{
    where } \; \varphi_k(\cdot) := L h_n^s \varphi \Big( \frac{\cdot -
    z_k}{h_n} \Big).
\end{equation*}
We have
\begin{align*}
  \norm{ f(\cdot; \omega) - f(\cdot; \omega') }_{L^2(P_Z)} &= \Big(
  \sum_{k=1}^m ( \omega_k - \omega_{k'} )^2 \int \varphi_k(z)^2
  P_Z(dz) \Big)^{1/2} \\
  &\geq \mu_0^{1/2} \rho(\omega, \omega') L^2 h_n^{2s + 1}
  \int_{S_{\mu}} \varphi(u)^2 du,
\end{align*}
where $S_{\mu} := \supp P_Z - \cup_{z} [a_{z}, b_{z}]$ (the union is
over the $z$ such that $\mu(z) = 0$, see Assumption~\ref{assd}), where $\mu_0
:= \min_{z \in S_\mu} \mu(z) > 0$ and where
\begin{equation*}
  \rho(\omega, \omega') := \sum_{k=1}^m \ind{ \omega_k \neq
    \omega'_{k}}
\end{equation*}
is the Hamming distance on $\Omega$. Using a result of
Varshamov-Gilbert (see~\cite{tsybakov03}) we can find a subset $\{
\omega^{(0)}, \ldots, \omega^{(M)} \}$ of $\Omega$ such that
$\omega^{(0)} = ( 0, \ldots, 0)$, $\rho(\omega^{(j)}, \omega^{(k)})
\geq m/8$ for any $0 \leq j < k \leq M$ and $M \geq 2^{m/8}$. Hence,
we have
\begin{equation*}
  \norm{ f(\cdot; \omega^{(j)}) - f(\cdot; \omega^{(k)}) }_{L^2(P_Z)}
  \geq D n^{-s / (2s + 1)},
\end{equation*}
where $D = \mu_0^{1/2} \int_{S_{\mu}} \varphi(u)^2 du / (8 c_0^{2s})
\geq 2$ for $c_0$ small enough. Moreover,
\begin{align*}
  \frac{1}{M} \sum_{k=1}^M K( P_{ f(\cdot, \omega^{(0)})}^n | P_{
    f(\cdot, \omega^{(k)})}^n ) &\leq \frac{n}{2 M \sigma_0^2}
  \sum_{k=1}^M \norm{ f(\cdot; \omega^{(0)}) - f(\cdot; \omega^{(k)})
  }_{L^2(P_Z)}^2 \\
  &\leq \frac{n}{2 \sigma_0^2}L^2h_n^{2s+1} \norm{\varphi}_2^2 m \leq
  \alpha \log M,
\end{align*}
where $\alpha := ( L^2 \norm{\varphi}_2^2 ) / (\sigma^2 c_0^{2s + 1}
\log 2) \in (0, 1/8)$ for $c_0$ small enough. The conclusion follows
from Theorem~\ref{thm:tsybokov}. $\hfill \square$

\subsection*{Proof of Theorem~\ref{thm:LPE}}

We recall that $r = \lfloor s \rfloor$ is the largest integer smaller
than $s$, and that $\lambda(M)$ stands for the smallest eigenvalue of
a matrix $M$.

\subsubsection*{Proof of (\ref{eq:LPE_point})}

First, we prove a bias-variance decomposition of the LPE at a fixed
point $z \in \supp P_Z$. This kind of result is commonplace, see for
instance \cite{fan_gijbels95, fan_gijbels96}. We introduce the
following weighted pseudo-inner product, for fixed $z \in \mathbb R$
and $h > 0$:
\begin{equation*}
  \prodsca{f}{g}_h := \frac{1}{m \bar P_Z[I(z, h)]} \sumim
  f(Z_i) g(Z_i) \ind{Z_i \in I(z, h)},
\end{equation*}
where we recall that $I(z, h) = [z - h, z + h]$, and that $\bar P_Z$
is given by~\eqref{eq:empirical_design}. We consider the associated
pseudo-norm $\norm{g}^2_h := \prodsca{g}{g}_h$. We introduce the power
functions $\varphi_a(\cdot) := ((\cdot - z) / h)^a$ for $a \in \{ 0,
\ldots, r \}$, which satisfy $\norm{\varphi_a}_h \leq 1$.

Note that the entries of the matrix $\bar {\mathbf Z}_m = \bar
{\mathbf Z}_m(z, h)$ (see~\eqref{eq:defMatrixZh}) satisfy $(\bar
{\mathbf Z}_m(z, h))_{a, b} := \prodsca{\varphi_a}{\varphi_b}_h$ for
$(a, b) \in \{ 0, \ldots, r \}^2$. Hence,~\eqref{eq:defLPE} is
equivalent to find $\bar P \in \text{Pol}_r$ such that
\begin{equation}
  \label{eq:LPE_variational}
  \prodsca{\bar P}{\varphi_a}_h = \prodsca{Y}{\varphi_a}_h
\end{equation}
for any $a \in \{0, \ldots, r \}$, where $ \prodsca{Y}{\varphi}_h := (
m \bar P_Z[I(z, h)])^{-1} \sumim Y_i \varphi(Z_i) \ind{Z_i \in I(z,
  h)}$. In other words, $\bar P$ is the projection of $Y$ onto
$\text{Pol}_r$ with respect to the inner product
$\prodsca{\cdot}{\cdot}_h$. For $e_1 := (1, 0, \ldots, 0) \in \mathbb
R^{r+1}$, we have
\begin{align*}
  \bar f(z) - f(z) = e_1^{\top} \bar {\mathbf Z}_m^{-1} \bar {\mathbf
    Z}_m (\bar \theta - \theta)
\end{align*}
whenever $\lambda(\bar {\mathbf Z}_m) > 0$, where $\bar \theta$ is the
coefficient vector of $\bar P$ and $\theta$ is the coefficient vector
of the Taylor polynomial $P$ of $f$ at $z$ with degree $r$. In view
of~\eqref{eq:LPE_variational}:
\begin{align*}
  ( \bar {\mathbf Z}_m (\bar \theta - \theta))_a = \prodsca{\bar P -
    P}{\varphi_a}_h = \prodsca{Y - P}{\varphi_a}_h,
\end{align*}
thus $\bar {\mathbf Z}_m (\bar \theta - \theta)) = B + V$ where $(B)_a
:= \prodsca{f - P}{\varphi_a}_h$ and $(V)_a := \prodsca{\sigma(\cdot)
  \xi}{\varphi_a}_h$. The bias term satisfies $|e_1^{\top} \bar
{\mathbf Z}_m^{-1} B| \leq (r+1)^{1/2} \norm{ \bar {\mathbf Z}_m^{-1}}
\norm{B}_{\infty}$ where for any $a \in \{ 0, \ldots, r\}$
\begin{equation*}
  |(B)_a| \leq \norm{f - P}_h \leq L h^s / r!.
\end{equation*}
Let $\bar{\mathbf Z}_m^{\sigma}$ be the matrix with entries
$(\bar{\mathbf Z}_m^{\sigma})_{a, b} := \prodsca{\sigma(\cdot)
  \varphi_a}{\sigma(\cdot) \varphi_b}_{h}$.  Since $V$ is,
conditionally on $Z_1^m = (Z_1, \ldots, Z_m)$, centered Gaussian with
covariance matrix $(m \bar P_Z[I(z, h)])^{-1} \bar{\mathbf
  Z}_m^{\sigma}$, we have that $e_1^{\top} \bar {\mathbf Z}_m^{-1} V$
is centered Gaussian with variance smaller than
\begin{equation*}
  (m \bar P_Z[I(z, h)])^{-1} e_1^{\top} \bar{\mathbf Z}_m^{-1}
  \bar{\mathbf Z}_m^{\sigma} \bar{\mathbf Z}_m^{-1} e_1 \leq \sigma_1^2
  (m \bar P_Z[I(z, h)])^{-1} \lambda (\bar{\mathbf Z}_m )^{-1}
\end{equation*}
where we used $\sigma(\cdot) \leq \sigma_1$. Hence, if $C_r :=
(r+1)^{1/2} / r!$, we obtain
\begin{equation*}
  E^m [ (\bar f(z) - f(z))^2 | Z_1^m] \leq \lambda( \bar{\mathbf
    Z}_m(z, h) )^{-2} \big( C_r L h^s + \sigma_1 (m \bar P_Z[I(z, h)]
  )^{-1/2} \big)^2
\end{equation*}
for any $z$, and the bandwidth choice~\eqref{eq:bandwidth}
entails~\eqref{eq:LPE_point}.

\subsection*{Proof of~\eqref{eq:LPE_L2}}

Let us consider the sequence of positive curves $h_m(\cdot)$ defined
as the point-by-point solution to
\begin{equation}
  \label{eq:def_h}
  L h_m(z)^s = \frac{\sigma_1}{(m P_Z[I(z, h_m(z))])^{1/2}}
\end{equation}
for all $z \in \supp P_Z$, where we recall $I(z, h) = [z-h, z+h]$, and
let us define
\begin{equation*}
  r_m(z) := L h_m(z)^s.
\end{equation*}
The sequence $h_m(\cdot)$ is the deterministic equivalent to the
bandwidth $H_m(\cdot)$ given by~\eqref{eq:bandwidth}. Indeed, with a
large probability, $H_m(\cdot)$ and $h_m(\cdot)$ are close to each
other in view of Lemma~\ref{lem:controlH} below. Under Assumption~\ref{assd}
we have $P_Z[I] \geq \gamma |I|^{\beta + 1}$, which entails together
with~\eqref{eq:def_h} that
\begin{equation}
  \label{eq:up_bound_h}
  h_m(z) \leq D_1 m^{-1 / (1 + 2s + \beta)}
\end{equation}
uniformly for $z \in \supp P_Z$, where $D_1 = (\sigma_1 / L)^{2 / (1 +
  2s + \beta)} (\gamma 2^{\beta + 1})^{-1/(1 + 2s +
  \beta)}$. Moreover, since $P_Z$ has a continuous density $\mu$ with
respect to the Lebesgue measure, we have
\begin{equation}
  \label{eq:lo_bound_h}
  h_m(z) \geq D_2 m^{-1 / (1 + 2s)}
\end{equation}
uniformly for $z \in \supp P_Z$, where $D_2 = (\sigma_1 / L)^{2/(1 +
  2s)} (2 \mu_{\infty})^{-1 / (2s + 1)}$. We recall that $P_Z^m$
stands for the joint law of $(Z_1, \ldots, Z_m)$.

\begin{lem}
  \label{lem:controlH}
  If $Z$ satisfies Assumption~\ref{assd}, we have for any $\epsilon \in (0,
  1/2)$
  \begin{equation*}
    P_Z^m \Big[ \sup_{z \in \supp(P_Z)} \Big| \frac{H_m(z)}{h_m(z)}
    - 1 \Big| > \epsilon \Big] \leq \exp( - D \epsilon^2 m^{\alpha})
  \end{equation*}
  for $m$ large enough, where $\alpha := 2s / (1 + 2s + \beta) $ and
  $D$ is a constant depending on $\sigma_1$ and $L$.
\end{lem}

The next lemma provides an uniform control on the smallest eigenvalue
of $\bar {\mathbf Z}_m(z) := \bar {\mathbf Z}_m(z, H_m(z))$ under
Assumption~\ref{assd}.
\begin{lem}
  \label{lem:VP}
  If $Z$ satisfies Assumption~\ref{assd}, there exists $\lambda_0 > 0$
  depending on $\beta$ and $s$ only such that
  \begin{equation*}
    P_Z^m \big[ \inf_{z \in \supp P_Z} \lambda(\bar{\mathbf Z}_m(z)) \leq
    \lambda_0 \big] \leq \exp( - D m^{\alpha} ),
  \end{equation*}
  for $m$ large enough, where $\alpha = 2s / (1 + 2s + \beta)$, and
  $D$ is a constant depending on $\gamma, \beta, s, L, \sigma_1$.
\end{lem}
The proofs Lemmas~\ref{lem:controlH} and~\ref{lem:VP} are given in
Section~\ref{sec:lemma_hell}. We consider the event
\begin{equation*}
  \Omega_m(\epsilon) := \big\{  \inf_{z \in \supp P_Z}
  \lambda(\bar{\mathbf Z}_m(z)) > \lambda_0 \big\} \cap \big\{ \sup_{z \in
    \supp P_Z} | H_m(z) / h_m(z) - 1| \leq \epsilon \big\},
\end{equation*}
where $\epsilon \in (0, 1/2)$. We have for any $f \in H^Q(s, L)$
\begin{equation*}
  E^m [ \norm{\tau_Q(\bar f) - f}_{L^2(P_Z)}^2
  \ind{\Omega_m(\epsilon)}] \leq \lambda_0^{-2} (1 + \epsilon)^{2s}
  \frac{\sigma_1^2}{m} \int \frac{P_Z(dz)}{\int_{z - h_m(z)}^{z +
      h_m(z)} P_Z(dt)},
\end{equation*}
where we used together the definition of
$\Omega_m(\epsilon)$,~\eqref{eq:LPE_point} and~\eqref{eq:def_h}. Let
us denote $I := \supp P_Z$ and let $I_{z*}$ be the intervals from
Assumption~\ref{assd}. Using together the fact that $\min_{z \in I -
  \cup_{z*} I_{z*}} \mu(z) > 0$ and~\eqref{eq:lo_bound_h}, we obtain
\begin{equation*}
  \frac{\sigma_1^2}{m} \int_{I - \cup_{z*} I_{z*}}
  \frac{P_Z(dz)}{\int_{z-h_m(z)}^{z + h_m(z)} P_Z(dt)} \leq C m^{-2s /
    (2s + 1)}.
\end{equation*}
Using the monoticity constraints from Assumption~\ref{assd}, we obtain
\begin{align*}
  \frac{\sigma_1^2}{m} &\int_{I_{z*}} \frac{P(dz)}{\int_{z -
      h_m(z)}^{z + h_m(z)} P_Z(dt)} \\
  &\leq \frac{\sigma_1^2}{m} \Big( \int_{z* - a_{z*}}^{z*}
  \frac{\mu(z) dz}{\int_{z - h_m(z)}^{z} \mu(t) dt} + \int_{z*}^{z* +
    b_{z*}} \frac{\mu(z) dz}{\int_{z}^{z + h_m(z)}
    \mu(t) dt} \Big) \\
  &\leq \frac{\sigma_1^2}{m} \int_{I_{z*}} h_m(z)^{-1} dz \leq C
  m^{-2s / (2s + 1)},
\end{align*}
hence $E^m [ \norm{\tau_Q(\bar f) - f}_{L^2(P_Z)}^2
\ind{\Omega_m(\epsilon)}] \leq C m^{-2s / (2s + 1)}$ uniformly for $f
\in H^Q(s, L)$. Using together Lemmas~\ref{lem:controlH}
and~\ref{lem:VP}, we obtain $E^m [ \norm{\tau_Q(\bar f) -
  f}_{L^2(P_Z)}^2 \ind{\Omega_m(\epsilon)^{\complement}} ] = o(n^{-2s
  / (2s + 1)})$, and~\eqref{eq:LPE_L2} follows. \hfill $\square$

\subsection*{Proof of Theorem~\ref{thm:oracle}}

In model~\eqref{eq:model}, when the noise $\epsilon$ is centered and
such that $E(\epsilon^2) = 1$, the risk of a function $\bar g :
\mathbb R^d \rightarrow \mathbb R$ is given by
\begin{equation*}
  A(\bar g) := E[ (Y - \bar g(X) )^2 ] = E [ \sigma(X) ^2 ] +
  \norm{\bar g - g}_{L^2(P_X)}^2,
\end{equation*}
where $g$ is the regression function. Therefore, the excess risk
satisfies
\begin{equation*}
  A(\bar g) - A = \norm{\bar g - g}_{L^2(P_X)}^2,
\end{equation*}
where $A := A(g) = E [ \sigma(X) ^2 ]$. Let us introduce $n := |D|$
the size of the learning sample, and $M := |\Lambda|$ the size of the
dictionary of functions $\{ \bar g^{(\lambda)} ; \lambda \in \Lambda
\}$. The least squares of $\bar g$ over the learning sample is given
by
\begin{equation*}
  A_n(\bar g) := \frac{1}{n} \sum_{i = 1}^{n} ( Y_i - \bar g(X_i)
  )^2.
\end{equation*}
We begin with a linearization of these risks. We consider the convex
set
\begin{equation*}
  \cC := \Big\{ (\theta_{\lambda})_{\lambda \in \Lambda} \text{ such
    that } \theta_{\lambda} \geq 0 \text{ and } \sum_{\lambda \in
    \Lambda} \theta_\lambda = 1 \Big\},
\end{equation*}
and define the linearized risks on $\cC$ as
\begin{equation*}
  \tilde{A}(\theta) := \sum_{\lambda \in \Lambda} \theta_\lambda
  A(\bar g^{(\lambda)}), \quad
  \tilde{A}_n(\theta) := \sum_{\lambda \in \Lambda} \theta_\lambda
  A_n(\bar g^{(\lambda)}),
\end{equation*}
which are linear versions of the risk $A$ and its empirical version
$A_n$. The exponential weights $w = (w_\lambda)_{\lambda \in \Lambda}
:= ( w( \bar g^{(\lambda)} ) )_{\lambda \in \Lambda}$ are actually the
unique solution of the minimization problem
\begin{equation}
  \label{eq:oracle_minimization}
  \min \Big( \tilde{A}_n(\theta) + \frac{1}{T n} \sum_{\lambda \in
    \Lambda} \theta_\lambda \log \theta_\lambda \;\big|\;  (\theta_\lambda)
  \in \cC \Big),
\end{equation}
where $T > 0$ is the temperature parameter in the
weights~\eqref{eq:weights}, and where we use the convention $0 \log 0
= 0$. Let $\hat \lambda \in \Lambda$ be such that $A_n ( \bar g^{(\hat
  \lambda)} ) = \min_{\lambda \in \Lambda} A_n(\bar
g^{(\lambda)})$. Since $\sum_{\lambda \in \Lambda} w_\lambda \log
\big( \frac{w_\lambda}{1 / M } \big) = K(w|u) \geq 0$ where $K(w|u)$
denotes the Kullback-Leibler divergence between the weights $w$ and
the uniform weights $u := (1 / M )_{\lambda \in \Lambda}$, we have
together with~\eqref{eq:oracle_minimization}:
\begin{align*}
  \tilde A_n(w) &\leq \tilde A_n(w) + \frac{1}{T n} K(w | u) \\
  &= \tilde A_n(w) + \frac{1}{T n} \sum_{\lambda \in \Lambda}
  w_\lambda \log w_\lambda + \frac{\log M}{T n} \\
  &\leq \tilde A_n( e_{\hat \lambda} ) + \frac{\log M}{T n},
\end{align*}
where $e_\lambda \in \cC$ is the vector with $1$ for the $\lambda$-th
coordinate and $0$ elsewhere. Let $a>0$ and $A_n := A_n(g)$. For any
$\lambda \in \Lambda$, we have
\begin{align*}
  \tilde{A}(w) - A &= (1 + a) ( \tilde{A}_n(w) - A_n ) + \tilde{A}(w)
  - A - (1 + a) (\tilde{A}_n(w) - A_n) \\
  &\leq (1 + a) ( \tilde A_n(e_\lambda) - A_n) + (1 + a) \frac{\log
    M}{T n} \\
  &+ \tilde{A}(w) - A - (1 + a) (\tilde{A}_n(w) - A_n).
\end{align*}
Let us denote by $E_K$ the expectation with respect to $P_K$, the
joint law of the learning sample for a noise $\epsilon$ which is
bounded almost surely by $K > 0$. We have
\begin{align*}
  E_K \big[ \tilde{A}(w) - A \big] \leq (1 &+ a) \min_{\lambda \in
    \Lambda} ( \tilde A_n(e_\lambda) - A_n) + (1 + a) \frac{\log
    M}{T n} \\
  &+ E_K \big[ \tilde{A}(w) - A - (1 + a) (\tilde{A}_n(w) - A_n)
  \big].
\end{align*}
Using the linearity of $\tilde A$ on $\cC$, we obtain
\begin{equation*}
  \tilde A(w) - A - (1 + a) (\tilde A_n(w) - A_n) \leq \max_{g \in
    \mathcal G_{\Lambda}} \big( A(g) - A - (1 + a) (A_n(g) - A_n)
  \big),
\end{equation*}
where $\mathcal G_\Lambda := \{ \bar g^{(\lambda)} \;;\; \lambda \in
\Lambda \}$. Then, using Bernstein inequality, we obtain for all
$\delta > 0$
\begin{align*}
  P_K \big[ &\tilde A(w) - A - (1+a) ( \tilde A_n(w) - A_n)\geq \delta
  \big] \\
  &\leq \sum_{g \in \mathcal G_{\Lambda}} P_K \Big[ A(g) - A - (
  A_n(g) - A_n ) \geq \frac{ \delta + a( A(g) - A ) }{1 + a} \Big] \\
  &\leq \sum_{g \in \mathcal G_{\Lambda}} \exp \Big( -\frac{n ( \delta
    + a(A(g) - A)) ^2 (1 + a)^{-1} }{8 Q^2 (1 + a) (A(g) - A) + 2 ( 6
    Q^2 + 2 \sigma K ) (\delta + a(A(g) - A)) / 3} \Big).
\end{align*}
Moreover, we have for any $\delta > 0$ and $g \in \mathcal
G_{\Lambda}$,
\begin{equation*}
  \frac{( \delta + a(A(g) - A))^2 (1 + a)^{-1} }{8 Q^2 (A(g) - A) + 2
    ( 6 Q^2 (1 + a)  + 2 \sigma K ) (\delta + a(A(g) - A)) / 3} \geq
  C(a, K) \delta,
\end{equation*}
where $C(a, K) := \big( 8 Q^2 (1 + a)^2 / a + 4 ( 6 Q^2 + 2 \sigma K)
(1 + a) / 3 \big)^{-1}$, thus
\begin{equation*}
  E_K \big[ \tilde A(w) - A - (1 + a) ( \tilde A_n(w) - A_n ) \big]
  \leq 2u + M \frac{\exp(-n C(a,K) u)}{ n C(a,K) }.
\end{equation*}
If we denote by $\gamma_A$ the unique solution of $\gamma = A
\exp(-\gamma)$, where $A > 0$, we have $(\log A) / 2 \leq \gamma_A
\leq \log A$. Thus, if we take $u = \gamma_M / (n C(a,K))$, we obtain
\begin{equation*}
  E_K \big[ \tilde A(w) - A - (1 + a) ( \tilde A_n(w) - A_n ) \big]
  \leq \frac{3 \log M}{C(a, K) n}.
\end{equation*}
By convexity of the risk, we have
\begin{equation*}
  \tilde A(w) - A \geq A( \hat g) - A,
\end{equation*}
thus
\begin{equation*}
  E_K \big[ \norm{ \hat g - g }_{L^2(P_X)}^2 \big] \leq (1 + a)
  \min_{\lambda \in \Lambda} \norm{ \bar g^{(\lambda)} -g}_{L²(P_X)}^2
  + C_1 \frac{\log M}{n},
\end{equation*}
where $C_1 := (1 + a) (T^{-1} + 3 C(a, K)^{-1})$. It remains to prove
the result when the noise is Gaussian. Let us denote
$\epsilon_\infty^n := \max_{1 \leq i \leq n} |\epsilon_i|$. For any $K
> 0$, we have
\begin{align*}
  E \big[ \norm{\hat g - g}_{L^2(P_X)}^2 \big] &= E \big[ \norm{\hat g
    - g}_{L^2(P_X)}^2 \ind{ \epsilon_\infty^n \leq K } \big] + E \big[
  \norm{\hat g - g}_{L^2(P_X)}^2 \ind{ \epsilon_\infty^n > K }
  \big] \\
  &\leq E_K \big[ \norm{\hat g - g}_{L^2(P_X)}^2 \big] + 2 Q^2 P[
  \epsilon_\infty^n > K ].
\end{align*}
For $K = K_n := 2 (2 \log n)^{1 / 2 }$, we obtain using standard
results about the maximum of Gaussian vectors that $P[
\epsilon_\infty^n > K_n] \leq P[ \epsilon_\infty^n - E
[\epsilon_\infty^n] > (2 \log n)^{1 / 2}] \leq 1/n $, which concludes
the proof of the Theorem. $\hfill \square$

\section{Proof of the lemmas}
\label{sec:lemma_hell}

\subsection*{Proof of Lemma~\ref{lem:controlH}}

Using together~\eqref{eq:bandwidth} and~\eqref{eq:def_h}, if
$I_m^{\epsilon}(z) := [z - (1 + \epsilon) h_m(z), z + (1 + \epsilon)
h_m(z)]$ and $I_m(z) := I_m^{0}(z)$, we obtain for any $\epsilon \in
(0, 1/2)$:
\begin{align*}
  \{ H_m(z) \leq (1 + \epsilon) h_m(z) \} &= \big\{ (1 +
  \epsilon)^{2s} \bar P_Z[ I_m^{\epsilon}(z) ] \geq P_Z[I_m(z)] \big\} \\
  &\supset \big\{ (1 + \epsilon)^{2s} \bar P_Z[ I_m(z) ] \geq
  P_Z[I_m(z)] \big\},
\end{align*}
where we used the fact that $\epsilon \mapsto P_Z[ I_m^{\epsilon}(z)]$
is nondecreasing. Similarly, we have on the other side
\begin{equation*}
  \{ H_m(z) > (1 - \epsilon) h_m(z) \} \supset \big\{ (1 -
  \epsilon)^{2s} \bar P_Z[ I_m(z) ] \leq P_Z [ I_m(z) ] \big\}.
\end{equation*}
Thus, if we consider the set of intervals
\begin{equation*}
  \cI_m := \bigcup_{z \in \supp{P_Z}} \big\{ I_m(z) \big\},
\end{equation*}
we obtain
\begin{equation*}
  \Big\{ \sup_{z \in \supp P_Z} \Big| \frac{H_m(z)}{h_m(z)} - 1 \Big|
  \geq \epsilon \Big\} \subset \Big\{ \sup_{I \in \cI_m} \Big|
  \frac{\bar P_Z[I]}{P_Z[I]} - 1 \Big| \geq \epsilon / 2 \Big\}.
\end{equation*}
Using together~\eqref{eq:def_h} and~\eqref{eq:up_bound_h}, we obtain
\begin{equation}
  \label{eq:minorproba}
  P_Z[ I_m(z) ] = \sigma_1^2 / (m L^2 h_m(z)^{2s})  \geq D m^{-(\beta +
    1) / (1 + 2s + \beta)} =: \alpha_m.
\end{equation}
Hence, if $\epsilon' := \epsilon (1 + \epsilon / 2) / (\epsilon + 2)$,
we have
\begin{align*}
  \Big\{ \sup_{I \in \cI_m} \Big| \frac{\bar P_Z[I]}{P_Z[I]} - 1 \Big|
  \geq \epsilon / 2 \Big\} \subset &\Big\{ \sup_{I \in \cI_m}
  \frac{\bar P_Z[I] - P_Z[I] }{\sqrt{\bar P_Z[I]}} \geq \epsilon'
  \alpha_m^{1/2} \Big\} \\
  & \cup \Big\{ \sup_{I \in \cI_m} \frac{P_Z[I] - \bar P_Z[I]
  }{\sqrt{P_Z[I]}} \geq \epsilon \alpha_m^{1/2} / 2 \Big\}.
\end{align*}
Then, Theorem~\ref{thm:vapnik74} (see Appendix) and the fact that the
shatter coefficient satisfies $\cS( \cI_m, m) \leq m(m+1) / 2$ entails
the Lemma. $\hfill \square$

\subsection*{Proof of Lemma~\ref{lem:VP}}

Let us denote $\bar {\mathbf Z}_m(z) := \bar {\mathbf Z}_m(z, H_m(z))$
where $\bar {\mathbf Z}_m(z, h)$ is given by~\eqref{eq:defMatrixZh}
and where $H_m(z)$ is given by~\eqref{eq:bandwidth}. Let us define the
matrix $\tilde {\mathbf Z}_m(z) := \tilde {\mathbf Z}_m(z, h_m(z))$
where
\begin{equation*}
  (\tilde {\mathbf Z}_m(z, h))_{a, b} := \frac{1}{m P_Z[I(z, h)]}
  \sumim \Big( \frac{Z_i - z}{h} \Big)^{a + b} \ind{Z_i \in I(z, h)}.
\end{equation*}

\noindent
\emph{Step~1.} Let us define for $\epsilon \in (0, 1)$ the event
\begin{equation*}
  \Omega_1(\epsilon) := \Big\{ \sup_{z \in \supp P_Z} \Big|
  \frac{H_m(z)}{h_m(z)} - 1 \Big| \leq \epsilon \Big\} \cap \Big\{
  \sup_{z \in \supp P_Z} \Big| \frac{\bar P_Z[I(z, H_m(z))]}{P_Z[ I(z,
    h_m(z))]} - 1 \Big| \leq \epsilon \Big\}.
\end{equation*}
For a matrix $A$, we denote $\norm{A}_{\infty}:= \max_{a, b} |(A)_{a,
  b}|$. We can prove that on $\Omega_1(\epsilon)$, we have
\begin{equation*}
  \norm{ \bar {\mathbf Z}_m(z) - \tilde {\mathbf Z}_m(z)}_{\infty}
  \leq \epsilon.
\end{equation*}
Moreover, using Lemma~\ref{lem:controlH}, we have $P_Z^m[
\Omega_1(\epsilon)^{\complement} ] \leq C \exp( -D \epsilon^2
m^{\alpha})$.  Hence, on $\Omega_1(\epsilon)$, we have for any $v \in
\mathbb R^d$, $\norm{v}_2 = 1$
\begin{equation*}
  v\T \bar {\mathbf Z}_m(z) v \geq v\T \tilde {\mathbf Z}_m(z) v -
  \epsilon
\end{equation*}
uniformly for $z \in \supp P_Z$.

\noindent
\emph{Step~2.} We define the deterministic matrix $\mathbf Z(z) :=
\mathbf Z(z, h_m(z))$ where
\begin{equation*}
  ({\mathbf Z}(z, h))_{a, b} := \frac{1}{P_Z[I(z, h)]}
  \int_{I(z, h)} \Big( \frac{t - z}{h} \Big)^{a + b} P_Z(dt),
\end{equation*}
and
\begin{equation*}
  \lambda_0 := \limInf_m \inf_{z \supp P_Z} \lambda\big( {\mathbf Z}(z,
  h_m(z)) \big).
\end{equation*}
We prove that $\lambda_0 > 0$. Two cases can occur: either $\mu(z) =
0$ or $\mu(z) > 0$. We show that in both cases, the $\limInf$ is
positive. If $\mu(z) > 0$, the entries $(\mathbf Z(z, h_m(z)))_{a, b}$
have limit $(1 + (-1)^{a+b})/(2(a + b + 1))$, which defines a positive
definite matrix. If $\mu(z) = 0$, we know that the density
$\mu(\cdot)$ of $P_Z$ behaves as the power function $|\cdot -
z|^{\beta(z)}$ around $z$ for $\beta(z) \in (0, \beta)$. In this case,
$(\mathbf Z(z, h_m(z)))_{a, b}$ has limit $(1 + (-1)^{a+b})(\beta(z) +
1) / [ 2(1 + a + b + \beta(z)) ]$, which defines also a definite
positive matrix.

\noindent
\emph{Step~3.} We prove that
\begin{equation*}
  P_Z^m [ \sup_{z \in \supp P_Z} \norm{ \tilde {\mathbf Z}_m(z) -
    {\mathbf Z}(z) }_{\infty} > \epsilon] \leq \exp( -D \epsilon^2
  m^{\alpha}).
\end{equation*}
We consider the sets of nonnegative functions (we recall that $I(z, h)
= [z-h, z+h]$)
\begin{equation*}
  F^{(\text{even})} := \bigcup_{\substack{z \in \supp P_Z \\ a \text{
        even and } 0 \leq a \leq 2 r }} \Big\{ \Big( \frac{\cdot -
    z}{h_m(z)}\Big)^a \ind{I(z, h_m(z))}(\cdot) \Big\},
\end{equation*}
\begin{equation*}
  F_{+}^{(\text{odd})} := \bigcup_{\substack{z \in \supp P_Z \\ a \text{
        odd and } 0 \leq a \leq 2 r}} \Big\{ \Big( \frac{\cdot -
    z}{h_m(z)}\Big)^a \ind{ [z, z + h_m(z)]}(\cdot) \Big\},
\end{equation*}
\begin{equation*}
  F_{-}^{(\text{odd})} := \bigcup_{\substack{z \in \supp P_Z \\ a \text{
        odd and } 0 \leq a \leq 2 r}} \Big\{ \Big( \frac{z -
    \cdot}{h_m(z)}\Big)^a \ind{ [z - h_m(z), z]}(\cdot) \Big\}.
\end{equation*}
Writing $I(z, h_m(z)) = [z - h_m(z), z) \cup [z, z + h_m(z)]$ when
$a+b$ is odd, and since
\begin{equation*}
  P_Z[ I(z, h_m(z)) ] \geq E f(Z_1)
\end{equation*}
for any $f \in F := F^{(\text{even})} \cup F_+^{(\text{odd})} \cup
F_-^{(\text{odd})}$, we obtain
\begin{equation*}
  \norm{ \tilde {\mathbf Z}_m(z) - {\mathbf Z}(z) }_{\infty} \leq
  \sup_{f \in F} \frac{| \frac{1}{m} \sum_{i=1}^m f(Z_i) - E f(Z_1)| }{
    E f(Z_1) }.
\end{equation*}
Hence, since $x \mapsto x / (x + \alpha)$ is increasing for any
$\alpha > 0$, and since $\alpha := E f(Z_1) \geq D m^{-(\beta + 1) /
  (1 + 2s + \beta)} =: \alpha_m$ (see~\eqref{eq:minorproba}), we
obtain
\begin{align*}
  \big\{ \sup_{z \in \supp P_Z} &\norm{ \tilde {\mathbf Z}_m(z) -
    {\mathbf Z}(z) }_{\infty} > \epsilon \big\} \\
  &\subset \Big\{ \sup_{f \in F} \frac{| \frac{1}{m} \sum_{i=1}^m
    f(Z_i) - E f(Z_1)| }{ \alpha_m + \frac{1}{m} \sum_{i=1}^m f(Z_i) +
    E f(Z_1) } > \epsilon / 2 \Big\}.
\end{align*}
Then, using Theorem~\ref{thm:haussler} (note that any $f \in F$ is
non-negative), we obtain
\begin{align*}
  P_Z^m [ \sup_{z \in \supp P_Z} &\norm{ \tilde {\mathbf Z}_m(z) -
    {\mathbf Z}(z) }_{\infty} > \epsilon] \\
  &\leq 4 E[ \cN_1(\alpha_m \epsilon/8, F, Z_1^m)] \exp \big( - D
  \epsilon^2 m^{2s / (1 + 2s + \beta)} \big).
\end{align*}
Together with the inequality
\begin{equation}
  \label{eq:covering}
  E[ \cN_1(\alpha_m \epsilon/8, F, Z_1^m)] \leq D (\alpha_m
  \epsilon)^{-1} m^{1 / (2s + 1) + (\beta-1) / (2s + \beta)},
\end{equation}
(see the proof below), this entails the Lemma. $\hfill \square$

\subsection*{Proof of~(\ref{eq:covering})}

It suffices to prove the inequality for $F^{(\text{even})}$ and a
fixed $a \in \{ 0, \ldots, 2 r \}$, since the proof is the same for
$F_+^{(\text{odd})}$ and $F_-^{(\text{odd})}$. We denote $f_z(\cdot)
:= ( (\cdot - z) / h_m(z) )^a \ind{I(z, h_m(z))}(\cdot)$. We prove the
following statement:
\begin{equation*}
  \cN ( \epsilon, F, \norm{\cdot}_\infty) \leq D \epsilon^{-1} m^{1 /
    (2s + 1) + (\beta-1) / (2s + \beta)},
\end{equation*}
which is stronger than~\eqref{eq:covering}, where
$\norm{\cdot}_\infty$ is the uniform norm over the support of
$P_Z$. Let $z, z_1, z_2 \in \supp P_Z$. We have
\begin{equation*}
  |f_{z_1}(z) - f_{z_2}(z)| \leq \max(a, 1) \Big|\frac{z - z_1}{h_1} -
  \frac{z - z_2 }{h_2} \Big| \ind{I_1 \cup I_2},
\end{equation*}
where $h_j := h_m(z_j)$ and $I_j := [ z_j - h_j, z_j + h_j]$ for $j =
1,2$. Hence,
\begin{equation*}
  |f_{z_1}(z) - f_{z_2}(z)| \leq \frac{ |h_1 -h_2| + |z_1 -
    z_2|}{\min(h_1 , h_2)}.
\end{equation*}
Using~\eqref{eq:def_h} together with a differentiation of $z \mapsto
h_m(z)^{2s} P_Z[ I(z, h_m(z)) ]$, we obtain that
\begin{align*}
  &|h_m(z_1) - h_m(z_2)| \\
  &\leq \sup_{z_1 \leq z \leq z_2} \Big| \frac{ h_m(z)^{2s+1} (\mu( z
    - h_m(z)) - \mu( z + h_m(z)) ) } { (2s \sigma_1^2) / (m L) +
    h_m(z)^{2s+1} (\mu(z - h_m(z)) + \mu(z + h_m(z))) } \Big| |z_1 -
  z_2|,
\end{align*}
for any $z_1 < z_2$ in $\supp \mu$. This entails together with
Assumption~\ref{assd},~\eqref{eq:up_bound_h} and~\eqref{eq:lo_bound_h}:
\begin{equation*}
  |h_m(z_1) - h_m(z_2)| \leq \frac{\mu_{\infty}}{2 s (\gamma
    L)^{ (2s + 1) / (2s + \beta + 1)}} \Big( \frac{m}{\sigma_1^2}
  \Big)^{\frac{\beta}{2s+\beta+1}} |z_1 - z_2|,
\end{equation*}
for any $z_1 < z_2$ in $\supp \mu$. Hence,
\begin{equation*}
  |f_{z_1}(z) - f_{z_2}(z)| \leq D m^{\frac{1}{2s + 1} + \frac{\beta -
      1}{2s + \beta}} |z_1 - z_2|,
\end{equation*}
which concludes the proof of~\eqref{eq:covering}. $\hfill \square$

\appendix

\section{\texorpdfstring{Some tools from empirical process theory}{Appendix A: Some tools from empirical process theory}}

Let $\cA$ be a set of Borelean subsets of $\mathbb{R}$. If $x_1^n :=
(x_1, \ldots, x_n) \in \mathbb R^n$, we define
\begin{equation*}
  N( \cA, x_1^n ) := \big| \big\{ \{ x_1, \ldots, x_n \}
  \cap A | A \in \cA \big\} \big|
\end{equation*}
and we define the \emph{shatter} coefficient
\begin{equation}
  \label{eq:shatter_coefficient}
  S(\cA, n) := \max_{x_1^n \in \mathbb R^n} N(\cA, (x_1,\ldots,x_n)).
\end{equation}
For instance, if $\cA$ is the set of all the intervals $[a,b]$ with
$-\infty \leq a < b \leq +\infty$, we have $S(\cA, n) = n(n + 1) / 2$.

Let $X_1, \ldots, X_n$ be i.i.d. random variables with values in
$\mathbb R$, and let us define $\mu[A] := P(X_1 \in A)$ and $\bar
\mu_n[A] := n^{-1} \sum_{i=1}^n \ind{X_i \in A}$. The following
inequalities for relative deviations are due to Vapnik and
Chervonenkis (1974), see for instance in \cite{vapnik98}.

\begin{theorem}[Vapnik and Chervonenkis (1974)]
  \label{thm:vapnik74}
  We have
  \begin{equation*}
    P\Big[ \sup_{A \in \mathcal A} \frac{\mu(A) - \bar
      \mu_n(A)}{\sqrt{\mu(A)}} > \epsilon \Big] \leq 4 \mathcal
    S(\mathcal A, 2 n) \exp(-n \epsilon^2 / 4)
  \end{equation*}
  and
  \begin{equation*}
    P\Big[ \sup_{A \in \mathcal A} \frac{\bar \mu_n(A) -
      \mu(A)}{\sqrt{\bar \mu_n(A)}} > \epsilon \Big] \leq 4 \mathcal
    S(\mathcal A, 2 n) \exp(-n \epsilon^2 / 4)
  \end{equation*}
  where $\mathcal S_{\mathcal A}(2 n)$ is the shatter coefficient of
  $\mathcal A$ defined by~\eqref{eq:shatter_coefficient}.
\end{theorem}

Let $(\cX, \tau)$ be a measured space and $\cF$ be a class of
functions $f : \cX \rightarrow [-K,K]$. Let us fix $p \geq 1$ and
$z_1^n \in \cX^n$. Define the semi-distance $d_p(f,g)$ between $f$ and
$g$ by
\begin{equation*}
  d_p(f,g) := \Big( \frac{1}{n} \sum_{i=1}^n |f(z_i) - g(z_i)|^p
  \Big)^{1/p}
\end{equation*}
and denote by $B^p(f,\epsilon)$ the $d_p$-ball with center $f$ and
radius $\epsilon$. The $\epsilon-$covering number of $\cF$ w.r.t $d_p$
is defined as
\begin{equation*}
  \cN_p(\epsilon, \cF, z_1^n) := \min\big( N \,|\, \exists
  f_1, \ldots, f_N  \text{ s.t. } \cF
  \subseteq \cup_{j=1}^M B^p (f_j,\epsilon) \big).
\end{equation*}

\begin{theorem}[Haussler (1992)]
  \label{thm:haussler}
  If $\cF$ consists of functions $f : \cX \rightarrow [0, K]$, we have
  \begin{equation*}
    P \Big[ \sup_{f\in\cF} \frac{\big| E[f(X_1)] - \frac{1}{n}
      \sum_{i=1}^nf(X_i) \big|}{\alpha +  E[f(X_1)] +
      \frac{1}{n} \sum_{i=1}^nf(X_i)} \geq \epsilon \Big]\leq 4 E[
    \cN_p(\alpha \epsilon / 8, \cF, X_1^n)] \exp \Big( -\frac{n \alpha
      \epsilon^2}{16 K^2} \Big).
  \end{equation*}
\end{theorem}

\bibliographystyle{ims}

\footnotesize%

\begin{thebibliography}{27}
\expandafter\ifx\csname natexlab\endcsname\relax\def\natexlab#1{#1}\fi
\expandafter\ifx\csname url\endcsname\relax
  \def\url#1{\texttt{#1}}\fi
\expandafter\ifx\csname urlprefix\endcsname\relax\def\urlprefix{URL }\fi
\providecommand{\eprint}[2][]{\url{#2}}

\bibitem[{Audibert and Tsybakov(2007)}]{audibert_tsybakov07}
\textsc{Audibert, J.-Y.} and \textsc{Tsybakov, A.} (2007).
\newblock Fast learning rates for plug-in estimators under the margin
  condition.
\newblock \textit{The Annals of Statistics}, \textbf{35}.
\MR{2336861}

\bibitem[{Birg{\'e}(2005)}]{b:04}
\textsc{Birg{\'e}, L.} (2005).
\newblock Model selection via testing: an alternative to (penalized) maximum
  likelihood estimators.
\newblock Available at \url{http://www.proba.jussieu.fr/mathdoc/textes/PMA-862.pdf}.
\MR{2219712}

\bibitem[{Catoni(2001)}]{catbook:01}
\textsc{Catoni, O.} (2001).
\newblock \textit{Statistical Learning Theory and Stochastic Optimization}.
\newblock Ecole d'{\'e}t{\'e} de Probabilit{\'e}s de Saint-Flour 2001, Lecture
  Notes in Mathematics, Springer, N.Y.
\MR{2163920}

\bibitem[{Delecroix et~al.(2003)Delecroix, H{\"a}rdle and
  Hristache}]{delecroix_etal03}
\textsc{Delecroix, M.}, \textsc{H{\"a}rdle, W.} and \textsc{Hristache, M.}
  (2003).
\newblock Efficient estimation in conditional single-index regression.
\newblock \textit{J. Multivariate Anal.}, \textbf{86} 213--226.
\MR{1997761}

\bibitem[{Delecroix et~al.(2006)Delecroix, Hristache and
  Patilea}]{delecroix_etal06}
\textsc{Delecroix, M.}, \textsc{Hristache, M.} and \textsc{Patilea, V.} (2006).
\newblock On semiparametric {$M$}-estimation in single-index regression.
\newblock \textit{J. Statist. Plann. Inference}, \textbf{136} 730--769.
\MR{2181975}

\bibitem[{Fan and Gijbels(1995)}]{fan_gijbels95}
\textsc{Fan, J.} and \textsc{Gijbels, I.} (1995).
\newblock Data-driven bandwidth selection in local polynomial fitting: variable
  bandwidth and spatial adaptation.
\newblock \textit{Journal of the Royal Statistical Society. Series B.
  Methodological}, \textbf{57} 371--394.
\MR{1323345}

\bibitem[{Fan and Gijbels(1996)}]{fan_gijbels96}
\textsc{Fan, J.} and \textsc{Gijbels, I.} (1996).
\newblock \textit{Local polynomial modelling and its applications}.
\newblock Monographs on Statistics and Applied Probability, Chapman \& Hall,
  London.
\MR{1383587}

\bibitem[{Geenens and Delecroix(2005)}]{geenens_delecroix05}
\textsc{Geenens, G.} and \textsc{Delecroix, M.} (2005).
\newblock A survey about single-index models theory.
\newblock \urlprefix\url{http://www.stat.ucl.ac.be/ISpub/dp/2005/dp0508.pdf}.

\bibitem[{Gy{\"o}rfi et~al.(2002)Gy{\"o}rfi, Kohler, Krzy{\.z}ak and
  Walk}]{kohler02}
\textsc{Gy{\"o}rfi, L.}, \textsc{Kohler, M.}, \textsc{Krzy{\.z}ak, A.} and
  \textsc{Walk, H.} (2002).
\newblock \textit{A distribution-free theory of nonparametric regression}.
\newblock Springer Series in Statistics, Springer-Verlag, New York.

\bibitem[{Horowitz(1998)}]{horowitz98}
\textsc{Horowitz, J.~L.} (1998).
\newblock \textit{Semiparametric methods in econometrics}, vol. 131 of
  \textit{Lecture Notes in Statistics}.
\newblock Springer-Verlag, New York.
\MR{1624936}

\bibitem[{Hristache et~al.(2001)Hristache, Juditsky and Spokoiny}]{spok01}
\textsc{Hristache, M.}, \textsc{Juditsky, A.} and \textsc{Spokoiny, V.} (2001).
\newblock Direct estimation of the index coefficient in a single-index model.
\newblock \textit{Ann. Statist.}, \textbf{29} 595--623.
\MR{1865333}

\bibitem[{Juditsky and Nemirovski(2000)}]{jn:00}
\textsc{Juditsky, A.} and \textsc{Nemirovski, A.} (2000).
\newblock Functional aggregation for nonparametric estimation.
\newblock \textit{Ann. Statist.}, \textbf{28(3)} 681--712.
\MR{1792783}

\bibitem[{Juditsky et~al.(2005{\natexlab{a}})Juditsky, Rigollet and
  Tsybakov}]{juditsky_etal05}
\textsc{Juditsky, A.}, \textsc{Rigollet, P.} and \textsc{Tsybakov, A.}
  (2005{\natexlab{a}}).
\newblock Learning by mirror averaging.
\newblock \urlprefix\url{http://arxiv.org/abs/math/0511468}.

\bibitem[{Juditsky et~al.(2005{\natexlab{b}})Juditsky, Nazin, Tsybakov and
  Vayatis}]{juditsky_nazin05}
\textsc{Juditsky, A.~B.}, \textsc{Nazin, A.~V.}, \textsc{Tsybakov, A.~B.} and
  \textsc{Vayatis, N.} (2005{\natexlab{b}}).
\newblock Recursive aggregation of estimators by the mirror descent method with
  averaging.
\newblock \textit{Problemy Peredachi Informatsii}, \textbf{41} 78--96.
\MR{2198228}

\bibitem[{Lecu{\'e}(2007)}]{lec9:07}
\textsc{Lecu{\'e}, G.} (2007).
\newblock Suboptimality of penalized empirical risk minimization in
  classification.
\newblock 20th Annual Conference On Learning Theory, COLT07, Springer,
  142--156.

\bibitem[{Leung and Barron(2006)}]{leung_barron06}
\textsc{Leung, G.} and \textsc{Barron, A.~R.} (2006).
\newblock Information theory and mixing least-squares regressions.
\newblock \textit{IEEE Trans. Inform. Theory}, \textbf{52} 3396--3410.
\MR{2242356}

\bibitem[{Nemirovski(2000)}]{n:00}
\textsc{Nemirovski, A.} (2000).
\newblock \textit{Topics in {N}on-parametric {S}tatistics}, vol. 1738 of
  \textit{Ecole d'{\'e}t{\'e} de Probabilit{\'e}s de Saint-Flour 1998, Lecture
  Notes in Mathematics}.
\newblock Springer, N.Y.
\MR{1775640}

\bibitem[{Stone(1982)}]{stone82}
\textsc{Stone, C.~J.} (1982).
\newblock Optimal global rates of convergence for nonparametric regression.
\newblock \textit{The Annals of Statistics}, \textbf{10} 1040--1053.
\MR{0673642}

\bibitem[{Stute and Zhu(2005)}]{stute_zhu05}
\textsc{Stute, W.} and \textsc{Zhu, L.-X.} (2005).
\newblock Nonparametric checks for single-index models.
\newblock \textit{Ann. Statist.}, \textbf{33} 1048--1083.
\MR{2195628}

\bibitem[{Tsybakov(2003{\natexlab{a}})}]{tsybakov03}
\textsc{Tsybakov, A.} (2003{\natexlab{a}}).
\newblock \textit{Introduction l'estimation non-paramerique}.
\newblock Springer.

\bibitem[{Tsybakov(2003{\natexlab{b}})}]{tsy:03}
\textsc{Tsybakov, A.~B.} (2003{\natexlab{b}}).
\newblock Optimal rates of aggregation.
\newblock \textit{Computational Learning Theory and Kernel Machines.
  B.Sch{\"o}lkopf and M.Warmuth, eds. Lecture Notes in Artificial
  Intelligence}, \textbf{2777} 303--313.
\newblock Springer, Heidelberg.

\bibitem[{Vapnik(1998)}]{vapnik98}
\textsc{Vapnik, V.~N.} (1998).
\newblock \textit{Statistical learning theory}.
\newblock Adaptive and Learning Systems for Signal Processing, Communications,
  and Control, John Wiley \& Sons Inc., New York.
\newblock , A Wiley-Interscience Publication.
\MR{1641250}

\bibitem[{Xia and H\"ardle(2006)}]{hardle_xia06}
\textsc{Xia, Y.} and \textsc{H\"ardle, W.} (2006).
\newblock Semi-parametric estimation of partially linear single-index models
  1162--1184.
\MR{2276153}

\bibitem[{Yang(2000{\natexlab{a}})}]{yang00}
\textsc{Yang, Y.} (2000{\natexlab{a}}).
\newblock Combining different procedures for adaptive regression.
\newblock \textit{J. Multivariate Anal.}, \textbf{74} 135--161.
\MR{1790617}

\bibitem[{Yang(2000{\natexlab{b}})}]{yang:00}
\textsc{Yang, Y.} (2000{\natexlab{b}}).
\newblock Mixing strategies for density estimation.
\newblock \textit{Ann. Statist.}, \textbf{28} 75--87.
\MR{1762904}

\bibitem[{Yang(2004)}]{yang04}
\textsc{Yang, Y.} (2004).
\newblock Aggregating regression procedures to improve performance.
\newblock \textit{Bernoulli}, \textbf{10} 25--47.
\MR{2044592}

\bibitem[{Yang and Barron(1999)}]{yang_barron99}
\textsc{Yang, Y.} and \textsc{Barron, A.} (1999).
\newblock Information-theoretic determination of minimax rates of convergence.
\newblock \textit{Ann. Statist.}, \textbf{27} 1564--1599.
\MR{1742500}

\end{thebibliography}

\end{document}